\documentclass{amsart}
\usepackage[title]{appendix}

\usepackage[square,numbers]{natbib}
\usepackage[english]{babel}

\usepackage{hyperref}

\usepackage{comment}

\usepackage[T1]{fontenc}

\usepackage{palatino}

\usepackage[autostyle]{csquotes}

\usepackage{amsmath}
\usepackage{amsfonts,amssymb,amsmath,latexsym,mathrsfs,bbm,stmaryrd}
\usepackage[all]{xy}
\usepackage[usenames]{color}
\usepackage{epsfig}

\usepackage{graphicx}
\usepackage{subcaption}
\usepackage{float}

\usepackage{amsthm}
\usepackage{thmtools}
\usepackage{thm-restate}

\usepackage{tikz}
\usetikzlibrary{arrows.meta}
\usetikzlibrary{matrix}
\usepackage{caption}
\usepackage{subcaption}

\usepackage{tikz}
\usepackage{xparse}

\NewDocumentCommand\DownArrow{O{2.0ex} O{black}}{%
	\mathrel{\tikz[baseline] \draw [<-, line width=0.5pt, #2] (0,0) -- ++(0,#1);}
}
\usetikzlibrary{shapes.arrows}
\newcommand{\arrowdown}[1]{
	\tikz[baseline=-1ex]{\node [draw, fill=black, single arrow, 
		minimum height=6.5ex, single arrow head extend=1ex,
		rotate=-90](arr) {}; \node[anchor=west] at ([xshift=0.5em]arr.90) {#1};}
}

\declaretheorem[numberwithin=section]{theorem}
\declaretheorem[sibling=theorem]{proposition}
\declaretheorem[sibling=theorem]{definition}
\declaretheorem[sibling=theorem]{corollary}
\declaretheorem[sibling=theorem]{lemma}

\declaretheorem[sibling=theorem]{notation}

\declaretheorem[sibling=theorem,style=remark]{remark}
\declaretheorem[sibling=theorem,style=remark]{example}

\newtheorem{thmy}{Theorem}
\newenvironment{thmx}{\stepcounter{theorem}\begin{thmy}}{\end{thmy}}

\newtheorem*{condition-non}{Condition A}
\newtheorem{condition-test}{Condition A}

\numberwithin{equation}{section}

\usepackage[mathcal]{euscript}
\usepackage{times}

\def\R{\mathbb R}
\def\C{\mathbb C}

\def\E{\mathbb E}

\def\1{\mathbbm{1}}

\newcommand{\x}{X}

\def\C{{\mathbb C}}
\def\B{{\mathbb B}}
\def\E{{\mathbb E}}

\def\CA{{\mathcal A}}
\def\CB{{\mathcal B}}
\def\CM{{\mathcal M}}

\def\i{{\rm i}}

\def\supp{{\text{\rm supp}}}
\def\Prob{{\text{\rm Prob}}}

\def\unit{{}}

\begin{document}

\title[Brown measure of the sum of an elliptic operator]{Brown measure of the sum of an elliptic operator and a free random variable in a finite von Neumann algebra}

\author{
	PING ZHONG
}

\address{
Department of Mathematics\\
University of Houston\\
Houston, Texas 77204-3008\\
{pzhong@central.uh.edu}
}

\maketitle

\begin{abstract}
Given an $n\times n$ random matrix $X_n$ with i.i.d. entries of unit variance, the circular law says that the empirical spectral distribution (ESD) of $X_n/\sqrt{n}$ converges to the uniform measure on the unit disk. Let $M_n$ be a deterministic matrix that converges in $*$-moments to an operator ${x}$. It is known from the work by \'{S}niady and Tao--Vu that the ESD of $X_n/\sqrt{n}+M_n$ converges to the Brown measure of ${x}+c$, where $c$ is Voiculescu's circular operator. We obtain a formula for the Brown measure of ${x}+c$ which provides a description of the limit distribution. This answers a question of Biane--Lehner for arbitrary operator ${x}$. 

Generalizing the case of circular and semi-circular operators, we also consider a family of twisted elliptic operators that are $*$-free from ${x}$. For an arbitrary twisted elliptic operator $g$, possible degeneracy then prevents a direct calculation of the Brown measure of ${x}+g$. We instead show that the whole family of Brown measures are the push-forward measures of the Brown measure of ${x}+c$ under a family of self-maps of the plane, which could possibly be singular. 
We calculate explicit formula for the case ${x}$ is self-adjoint. In addition, we prove that the Brown measure of the sum of an $R$-diagonal operator and a twisted elliptic element is supported in a deformed ring where the inner boundary is a circle and  the outer boundary is an ellipse. 

These results generalize some known results about free additive Brownian motions where the free random variable ${x}$ is assumed to be self-adjoint.  The approach is based on a Hermitian reduction and subordination functions. 

\end{abstract}


\tableofcontents

\section{Introduction}
\subsection{Brown measure of free random variables and random matrices}
Let $\CM$ be a von Neumann algebra with a faithful, normal, tracial state $\phi$. The Fuglede-Kadison determinant of $x\in \CM$ is defined by
\[
  \Delta(x)=\exp \left( \int_0^\infty \log t\, d\mu_{|x|}(t) \right),
\]
where $\mu_{|x|}$ is the spectral measure of $|x|$ with respect to $\phi$. 
Brown \cite{Brown1986} proved that the function  $\log\Delta(x-\lambda\unit)$ of $\lambda$ is a subharmonic function  whose Riesz measure is the unique, compactly supported probability measure $\mu_x$ on $\mathbb{C}$ with the property that
\[
  \log\Delta(x-\lambda\unit)=\int_\mathbb{C}\log |z-\lambda| d\mu_x(z), \qquad \lambda \in \mathbb{C}. 
\]
The measure $\mu_x$ is called the \emph{Brown measure} of $a$. 
In other words, $\mu_x$ is the distributional Laplacian of the function $\log\Delta(x-\lambda\unit)$ up to a constant.

Voiculescu's free probability theory is a suitable framework to describe the limits of the joint distribution of a family of random matrix models \cite{Voiculescu1991}. The convergence of $*$-mixed moments of suitable random matrix models have been studied well. The $*$-mixed moments of free random variables can be described using either analytic or combinatorial tools. For non-normal random matrix models, very little is known about the limit of the empirical spectral distribution (ESD) of a polynomial of independent random matrices, even for the sum or product of two random matrices. The Brown measure of an operator in $\CM$ is an analogue of eigenvalue distribution of a finite dimensional matrix. The Brown measures of the sum or product or a polynomial of free random variables are natural candidates for the limits of the ESD of the sum or product or a polynomial of suitable random matrix models as the size of the matrices tends to infinity.

Let $X_n$ be an $n\times n$ random matrix whose entries are independent identically distributed copies of a complex random variable with zero mean and unit variance. The circular law  says that the ESD of $X_n/\sqrt{n}$ converges weakly to the uniform measure on the unit disc which is also the Brown measure of Voiculescu's circular operator, denoted by $c$. The circular law was established in the 1960s by Ginibre \cite{Ginibre1965} for Gaussian distributed entries and was proved by Tao and Vu \cite{TaoVu2010} under the minimal assumptions after a long list of partial progresses (see \cite{Bordenave-Chafai-circular} and references therein). In fact, Tao and Vu proved results stronger than circular law. In particular, they showed the existence of the limit of the summation $X_n/\sqrt{n}+M_n$ where $M_n$ is a deterministic $n\times n$ matrix satisfying some technical conditions. In \cite{Sniady2002}, \'{S}niady showed that the ESD of $X_n/\sqrt{n}+M_n$ converges to the Brown measure of ${x}+c$ provided that $X_n$ is a Ginibre ensemble and $M_n$ converges in $*$-moments to ${x}$.  Hence, by combining Tao--Vu's replacement principle \cite[Section 2]{TaoVu2010}, we conclude that the ESD of $X_n/\sqrt{n}+M_n$ converges to the Brown measure of ${x}+c$ under the minimal requirements on $X_n$ when $M_n$ converges in $*$-moments to ${x}$.

In the above-cited paper, Tao and Vu did not pursue what the limit actually is (see \cite[Theorem 1.17]{TaoVu2010}) and they mentioned that the limit distribution ESD of $X_n/\sqrt{n}+M_n$  was established in a work of Krishnapur--Vu for the case where $M_n$ is a diagonal matrix  (equivalently, ${x}$ is a normal operator as the limit) and $X_n$ is a Ginibre ensemble. But that work has not appeared even as a preprint (communicated with Krishnapur). 
The case when ${x}$ is self-adjoint was known in the author's joint work with Ho \cite{HoZhong2020Brown} using PDE methods. The density formula when ${x}$ is a Gaussian distributed normal operator was established in a work of Bordenave-Caputo-Chafa\"{\i} \cite[Theorem 1.4]{Bordenave-Caputo-Chafai-cpam2014} using random matrix techniques. Their method can be extended to all normal operators \cite{Bordenave-Capitaine-cpam2016}.
Here we study the Brown measure of the sum of free circular operator and a $*$-free random variable ${x}$ with an arbitrary distribution, not necessarily normal, which finishes the description of the limit ESD of $X_n/\sqrt{n}+M_n$. 
The present work answers an earlier question of Biane--Lehner \cite[Section 5]{BianeLehner2001} in general case.

The twisted elliptic operators generalize circular operator, semi-circular operator and elliptic operator. Let $c_t$ be a circular operator with variance $t$, and let $g_{t,\gamma}$ be a twisted elliptic operator, and let ${x}$ be an operator $*$-free from $\{c_t, g_{t,\gamma}\}$. The calculation of the Brown measure of ${x}+g_{t,\gamma}$ is more involved than ${x}+c_t$ because there is possibly degeneracy. We show that there is a remarkable connection between the Brown measure of the addition of a free circular operator with a free random variable and the Brown measure of the sum of an elliptic operator with the same free random variable, in the following sense: the Brown measure of ${x}+g_{t,\gamma}$ is the push-forward measure of the Brown measures of ${x}+c_t$ under a natural self-map of the complex plane, which can be constructed explicitly. The push-forward connection between ${x}+c_t$ and ${x}+g_{t,\gamma}$ was first proved in some very special cases in \cite{HoHall2020Brown, HoZhong2020Brown} under the assumption that ${x}$ is self-adjoint. Our work extends this connection to full generality without any restriction on the distribution of ${x}$. The multiplicative analogue of the push-forward connection have been studied in a related context in our joint work with Ho \cite{HoZhong2020Brown} and by Hall-Ho \cite{HallHo2021Brown}. 

We calculate explicit density Brown measure formulas for the case where ${x}$ is self-adjoint. In addition, we describe the Brown measure of the sum of an $R$-diagonal operator and a twisted elliptic element. In this case, the Brown measure is supported in a deformed ring where the inner boundary is a circle and  the outer boundary is an ellipse. This can be viewed as a deformation of the limit distribution in the single ring theorem \cite{GuionnetKZ-single-ring2011} in random matrix theory. 

The present work extends previous results \cite{BianeLehner2001, HoHall2020Brown, Ho2020Brown, HoZhong2020Brown} for the sum of a self-adjoint operator with a circular operator or a (non-twisted) elliptic operator. All these work rely on some PDE methods and did not explain why subordination functions appeared in the Brown measure formulas. We use a completely different approach based on Hermitian reduction and subordination functions. The new method provides a conceptual explanation how subordination functions play a key role. The method and subordination results developed in this paper are likely to be useful in the study of non-normal random matrices. Our main results provide potential applications to unify various important deformed random matrix models (see \cite{Bordenave-Caputo-Chafai-cpam2014, CapitaineP2016_fullrank, Tao2013_circular_finiterank, TaoVu2010} for example).

\subsection{Statements of the results}
Let $g_{t,\gamma}$ be a twisted elliptic operator with parameters $t>0$ and $\gamma\in\mathbb{C}$ such that $|\gamma|\leq t$ (see Section \ref{section:elliptic} for the definition). Such operator has the same distribution as an operator of the form $e^{{\i}\theta} (s_{t_1}+{\i}s_{t_2})$ where $\theta\in[0,2\pi]$, and $s_{t_1}, s_{t_2}$ are semicircular operators with variances $t_1, t_2$ respectively that are freely independent in the sense of Voiculescu. The case when $\gamma=t$, the operator $g_{t,\gamma}$ is the semicircular operator $s_t$ with mean zero and variance $t$, and the case when $\gamma\in[-t,t]$ the operator $g_{t,\gamma}$ is an elliptic operator. 

Let ${x}\in\mathcal{M}$ be a random variable that is free from $\{c_t, g_{t,\gamma}\}$. We show that the Brown measure of ${x}+c_t$ can be calculated directly using subordination functions in free probability. The Brown measure of ${x}+g_{t,\gamma}$ is not calculated directly. Instead, we show that there is a natural push forward map between the Brown measure of ${x}+c_t$ and the Brown measure of ${x}+g_{t,\gamma}$. Our main results extends previous work \cite{HoHall2020Brown, Ho2020Brown, HoZhong2020Brown} and are also applicable to non-self-adjoint operators. 

Fix $t>0$ and ${x}\in \mathcal{M}$. Consider the open set (see Proposition \ref{prop:Xi-t-open-set}) 
\begin{equation}
\label{defn:Xi-t}
 \Xi_t=\left\{ \lambda\in\mathbb{C} : \phi \left( |\lambda-x|^{-2}\right)>\frac{1}{t}  \right\},
\end{equation}
where $|\lambda-x|=\big((\lambda-x)^*(\lambda-x)\big)^{1/2}$ by functional calculus and 
\[
 \phi \left( |\lambda-x|^{-2}\right)
   =\int_\mathbb{R}\frac{1}{u^2}d\mu_{|\lambda-x|}(u).
\]
For each $\lambda\in \Xi_t$, let $w(0;\lambda,t)$ be a positive function of $\lambda$  such that
\[
  \phi ( (\lambda-x)^*(\lambda-{x})+w(0;\lambda,t)^2\unit )^{-1} )=\frac{1}{t},
\]
and let $w(0;\lambda,t)=0$ for $\lambda\in\mathbb{C}\backslash\Xi_t$. 
For each $\lambda\in\mathbb{C}$ and $\varepsilon>0$, denote by $w(\varepsilon;\lambda,t)$ the imaginary part of a subordination function (valued at ${\i}\varepsilon$) for the free convolution of the symmetrization of $\mu_{|x-\lambda|}$ and the semicircular distribution with variance $t$. 
We will show that $w(0;\lambda,t)=\lim_{\varepsilon\downarrow 0} w(\varepsilon;\lambda,t)$ for any $\lambda\in\mathbb{C}$ (see Proposition \ref{prop:sub-identical-symmetric-measure}). 

We define the map $\Phi^{(\varepsilon)}_{t,\gamma}:\mathbb{C}\rightarrow\mathbb{C}$ by
\begin{align}
 \label{defn:Phi-t-gamma-epsitlon-v00}
 \Phi^{(\varepsilon)}_{t,\gamma} (\lambda)
 &=\lambda+\gamma\cdot \phi\left( (\lambda-x-c_t)^*(((\lambda-x-c_t)(\lambda-x-c_t)^*+\varepsilon^2)^{-1}) \right).
\end{align}
We can rewrite it as (see Theorem \ref{thm:sub-X0-gt-formula})
\[
  \Phi^{(\varepsilon)}_{t,\gamma} (\lambda)
    =\lambda+\gamma\cdot\phi\left( (\lambda-x)^*(((\lambda-x)(\lambda-x)^*+w(\varepsilon;\lambda,t))^{-1}) \right).
\]
We then denote
\[
\Phi_{t,\gamma} (\lambda)= \lambda+\gamma\cdot p_\lambda^{(0)}( w(0;\lambda,t)  ),
\]
where 
  \begin{equation}\label{eqn:1.2-intro}
    p_\lambda^{(0)}( w(0;\lambda,t) ) =
     \phi\bigg[ (\lambda-{x})^*\big( (\lambda\unit-{x})(\lambda\unit-{x})^*+w(0;\lambda,t)^2 \big)^{-1}\bigg].
  \end{equation}
For $\lambda\in\mathbb{C}\backslash\Xi_t$, the displayed formula in the right hand side of \eqref{eqn:1.2-intro} means $p_\lambda^{(0)}( 0 )=\lim_{\varepsilon\rightarrow 0}p_\lambda^{(0)}( \varepsilon )$ which is a finite value guaranteed by the fact that $\phi(|x-\lambda|^{-2})$ is finite for such $\lambda$.
We will show in Lemma \ref{lemma:regularization-Phi-uniform-convergence} that $\lim_{\varepsilon\rightarrow {0}}\Phi^{(\varepsilon)}_{t,\gamma}=\Phi_{t,\gamma}$
uniformly (see Lemma \ref{lemma:convergence-p-lambda}). 
In addition, $w(0;\lambda,t)$ is a real analytic function of $\lambda\in\Xi_t$ and hence the function $\Phi_{t,\gamma}$ is a real analytic function of $\lambda$ in the set $\Xi_t$.

\begin{thmx}
[See Theorem \ref{thm:main-FK-det-ct-0} and Theorem \ref{thm:main-FK-det-x-t-gamma-0}]
For every $\lambda\in \Xi_t$, we have 
 \begin{equation*}
  \Delta \big({x}+c_t-\lambda \big)^2
    ={\Delta\big( ({x}-\lambda\unit)^*({x}-\lambda\unit)+w^2 \unit \big)} \left[{\exp(-w^2/{t})}\right],
\end{equation*}
and, if the map $\lambda\mapsto \Phi_{t,\gamma}(\lambda)$ is non-singular at $\lambda\in\Xi_t$, then 
\begin{equation*}
{\Delta  ({x}+g_{t,\gamma}-z)^2}
  ={\Delta  ({x}+c_t-\lambda\unit)^2}
   \exp\left(H(\lambda)\right),
\end{equation*}
where $w=w(0;\lambda,t)$, $z=\Phi_{t,\gamma}(\lambda)$ and $H(\lambda)=\Re \bigg( \gamma( p_\lambda^{(0)}(w)) ^2\bigg)$.
\end{thmx}
 
We point out that we also have $\Delta({x}+c_t-\lambda\unit)=\Delta({x}-\lambda\unit)$ for any $\lambda\in \mathbb{C}\backslash \Xi_t$ (see Theorem \ref{thm:main-FK-det-ct-0}) which implies that $\mu_{{x}+c_t}$ and $\mu_{{x}}$ coincide in the interior of $\mathbb{C}\backslash  \Xi_t$. Moreover, we obtain a general result on the support of the Brown measure (see Theorem \ref{thm:support-Brown-general}) which allows us to deduce that the interior of $\mathbb{C}\backslash  \Xi_t$ is not in the support of $\mu_{{x}}$. Hence, our focus is the Brown measure of $\mu_{{x}+c_t}$ within the set $\Xi_t$. 
The above Fuglede-Kadison formulas are fundamental in our study which allows us to calculate the Brown measure formulas.
The Brown measure of ${x}+c_t$ can be described as follows.

	\begin{thmx}
	[See Theorem \ref{thm:BrownFormula-x0-ct-general} and Theorem \ref{thm:support-x0-c_t-general}]
	 The Brown measure of ${x}+c_t$ has no atom and is supported in $\overline{\Xi_t}$. It is absolutely continuous with respect to Lebesgue measure in the open set $\Xi_t$, and
	 the density of the Brown measure at any $\lambda\in \Xi_t$ is strictly positive which can be expressed as
	 \begin{equation}
	 \label{eqn:density-1.1-intro}
	  \frac{1}{\pi} \left( \frac{1}{t}-  \frac{\partial}{\partial \overline{\lambda}} \bigg( \phi \big( {x}^* ( ({x}-\lambda\unit)({x}-\lambda\unit)^*+w(0;\lambda,t)^2\unit )^{-1} \big) \bigg) \right).
	\end{equation}
	\end{thmx}
	After this work was done, in a followup joint work with Belinschi and Yin \cite[Section 7]{BelinschiYinZhong2021Brown}, we proved that the Brown measure of ${x}+c_t$ is absolutely continuous with respect to Lebesgue measure on $\mathbb{C}$ and we also showed that the density function is bounded by $1/\pi{t}$. 
	We can also rewrite the above density formula in some form without involving derivative by implicit differentiation (see \eqref{eqn:Brown-density-circular-positive}), which shows that the density \eqref{eqn:density-1.1-intro} is strictly positive in $\Xi_t$. 
The subordination function $w(0;\lambda,t)$ can be calculated explicitly for a large family of operators that include all self-adjoint operators (see Section \ref{section:addition-selfadj}) and $R$-diagonal operators (see Section \ref{section:addition-R-diag}). This generalizes a result in our earlier work with Ho \cite{HoZhong2020Brown} in which ${x}$ is self-adjoint.
	Unlike \cite{HoZhong2020Brown} and its generalizations for semicircular operators and elliptic operators \cite{HoHall2020Brown, Ho2020Brown}, we do not use PDE methods.

	\begin{thmx}
	[See Theorem \ref{thm:main-push-forward-general}]
	\label{thm:main-push-forward-introduction}
The Brown measure of ${x}+g_{t,\gamma}$ is the push-forward measure of the Brown measure of ${{x}+c_t}$ by the map $\lambda\mapsto \Phi_{t,\gamma}(\lambda)$. That is, for an arbitrary Borel measurable set $E$ in $\mathbb{C}$, we have 
	\begin{equation*}
	  \mu_{{x}+g_{t,\gamma}}(E)=\mu_{{x}+c_t}(\Phi_{t,\gamma}^{-1}(E))).
	\end{equation*}
The Brown measure of $x+g_{t,\gamma}$ is supported on $\Phi_{t,\gamma}(\overline{\Xi_t})$. 
	\end{thmx} 
	The above push-forward connection between two different models was known in earlier works \cite{HoHall2020Brown, Ho2020Brown, HoZhong2020Brown} under additional assumption that $x$ is self-adjoint for some certain elliptic operators.
	Theorem \ref{thm:main-push-forward-introduction} generalizes these results
	in two directions. Theorem \ref{thm:main-push-forward-introduction} is applicable for operators ${x}$ not necessarily self-adjoint. In addition, the twisted elliptic operator include semicircular operators and elliptic operators as special cases. What is more, the proof of Theorem \ref{thm:main-push-forward-introduction} provides a conceptual explanation about why such push-forward map exists. 
	
The push-forward map $\Phi_{t,\gamma}$ may be rewritten in terms of the Brown measure $\mu_{x+c_t}$ solely as follows. 
	\begin{thmx}
 [See Theorem \ref{thm:Phit_new} and Lemma \ref{lemma_continuityCauchy}]
 The Cauchy transform of $\mu_{x+c_t}$ is a continuous function defined on $\mathbb{C}$. 
 The map $\Phi_{t,\gamma}$ can be rewritten as
  \[
    \Phi_{t,\gamma}(\lambda)=\lambda+\gamma \cdot G_{\mu_{x+c_t}}(\lambda), \qquad \lambda\in\mathbb{C},
  \]
  where 
  \[
  G_{\mu_{x+c_t}}(\lambda)=\int_{\mathbb{C}}\frac{1}{\lambda-z}d{\mu_{x+c_t}}(z).
\] 
\end{thmx}

\begin{example}\label{example:circular-semicircular} [From circular law to elliptic/semicircular law]
Take $t=1$ and $|\gamma|\leq {1}$. Let ${x}=0$ and $\mu=\delta_0$. In this case, the formula for the set $\Xi_t$ is simplified as 
\[
  \Xi_t=\{\lambda\in\mathbb{C}: |\lambda|<1\},
\]
which is also the support of the standard circular operator $c$. Then the condition determining $w(0;\lambda,t)$ is written as
\[
  \frac{1}{|\lambda|^2+w(0;\lambda,t)^2}=1, \quad |\lambda|<1.
\]
Hence $w(0;\lambda,t)=\sqrt{1-|\lambda|^2}$. Then the push forward map is 
\begin{align*}
  \Phi_{t,\gamma}(\lambda)=
    \begin{cases}
      \lambda+\gamma\overline{\lambda}, \qquad & \text{for}\qquad 
          |\lambda|<1;\\
         \lambda+\frac{\gamma}{\lambda}  \qquad &\text{for}\qquad 
         |\lambda|\geq 1.
    \end{cases}
\end{align*}
Note that this map is a homeomorphism of $\mathbb{C}$ for $|\gamma|<1$, but when $|\gamma|=1$ it fails to be injective. 
\end{example}

The push-forward map $\Phi_{t,\gamma}$ not only connects different operator model and their corresponding random matrix models, but also appears in other related contexts.  For instance, this map was used in a recent preprint of Hall-Ho \cite[Section 2.2]{HallHo2022Conjecture}, where they proposed conjectures concerning zeros of characteristic polynomials of deformed i.i.d. random matrix model and deformed GUE or deformed elliptic random matrix model. 
 See also \cite{HallHJK2024zeros} for another work on zeros of certain random polynomials evolving under the heat flow. The transport map in \cite{HallHJK2024zeros} is the same map as above in the case of Weyl polynomials. 
	
	We conjecture that $\Phi_{t,\gamma}$ is a self-homeomorphism of $\mathbb{C}$ provided that $|\gamma|<t$, and we verify that this is indeed the case for self-adjoint and $R$-diagonal operators ${x}$ (see Theorem \ref{thm:main-self-adjoint-introduction} and Theorem \ref{thm:main-Rdiag-introduction}). 
	However, the map $\Phi_{t,\gamma}$ could be singular in general, as shown in Example \ref{example:circular-semicircular}.
	
Let us outline the strategy for the proof of Theorem \ref{thm:main-push-forward-introduction}. The possible degeneracy prevents a direct calculation of the Brown measure of ${x}+g_{t,\gamma}$. 
We will show that the map $\Phi^{(\varepsilon)}_{t,\gamma}$ defined in \eqref{defn:Phi-t-gamma-epsitlon-v00} is a self-homeomorphism of $\mathbb{C}$.
We prove that the regularized Brown measure $\mu_{{x}+g_{t,\gamma},\varepsilon}$ is the push-forward measure of $\mu_{{x}+c_t,\varepsilon}$ under the regularized push-forward map $\Phi^{(\varepsilon)}_{t,\gamma}$. After we establish the convergence $\Phi^{(\varepsilon)}_{t,\gamma} \rightarrow\Phi_{t,\gamma}$ as $\varepsilon$ tends to zero, Theorem \ref{thm:main-push-forward-introduction} would follow by push-forward connection between regularized Brown measures and the fact that regularized Brown measures converge to Brown measure weakly. We hence have the following commutative diagram. 

\begin{center}
\begin{tikzpicture}
\matrix (m) [matrix of math nodes,row sep=3em,column sep=4em,minimum width=2em]
{
	\mu_{{x}+c_t,\varepsilon} & \mu_{{x}+g_{t,\gamma},\varepsilon} \\
	\mu_{{x}+c_t} & \mu_{{x}+g_{t,\gamma}} \\};
\path[-stealth]
(m-1-1) edge node [left] {$\varepsilon\rightarrow 0$} (m-2-1)
edge [double] node [below] {$\Phi_{t,\gamma}^{(\varepsilon)}$} (m-1-2)
(m-2-1.east|-m-2-2) edge [double] node [below] {{$\Phi_{t,\gamma}$}}
(m-2-2)
(m-1-2) edge node [right] {$\varepsilon\rightarrow 0$} (m-2-2)
(m-2-1);
\end{tikzpicture}
\end{center}

	\begin{thmx}
	[See Theorem \ref{thm:dentisty-x0-x-t-gamma-main-self-adjoint} and Example \ref{example-Phi-elliptic}]
	\label{thm:main-self-adjoint-introduction}
	Let ${x}$ be a self-adjoint operator that is $*$-free from $g_{t,\gamma}$. 
	For any $|\gamma|\leq t$ with $\gamma\neq t$, the map $\Phi_{t,\gamma}$ is one-to-one and non-singular in $\Xi_t$. The Brown measure of ${x}+g_{t,\gamma}$ is the push-forward measure of the Brown measure of ${x}+c_t$ under the map $\Phi_{t,\gamma}$. 
	
	Moreover, the Brown measure $\mu_{{x}+g_{t,\gamma}}$ is concentrated on $\Phi_{t,\gamma}(\Xi_t)$ and the density is given by
	\begin{equation}
	\label{eqn:density-Brown-self-adjoint-introduction}
	   d\mu_{{x}+g_{t,\gamma}}(z)=\frac{1}{2\pi \tau_1}\frac{d\psi_t(a)}{d\delta(a)} dz_1 dz_2, \qquad z\in 
	   \Phi_{t,\gamma}(\Xi_t)
	\end{equation}
	where  $z=z_1+{\i}z_2=\Phi_{t,\gamma}(a+{\i}b)$, $\tau_1=t-\Re(\gamma)$, and $\psi_t, \delta$ are two increasing homeomorphisms of $\mathbb{R}$ onto $\mathbb{R}$. 
	In particular, if $\gamma\in \mathbb{R}$ (equivalently, $g_{t,\gamma}$ is an elliptic operator), then $z_1=\delta(a)$ depending only on $a$, in which case the Brown measure is constant along vertical lines. 
	\end{thmx}
	
	The density formula could be understood as follows. For any self-adjoint operator ${x}$, it is known \cite{Biane1997} that there is a continuous function $v_t$ such that 
	\[
	  \Xi_t=\{a+{\i}b\in\mathbb{C}: |b|<v_t(a) \}.
	\]
	See Section \ref{section-sub-x0-self-adjoint} for a review. It is shown \cite{HoZhong2020Brown} that
	the density of ${x}+c_t$ is constant along vertical segments in ${\Xi_t}$. In this case, for any fixed $a\in\mathbb{R}$ such that the vertical line through $a$ intersects the set $\Xi_t$, the map $b\mapsto\Phi_{t, \gamma}(a+{\i}b)$ is an affine transformation of $b$. Hence, the Brown measure of ${x}+g_{t,\gamma}$ is expected to be constant along the trajectory of $\Phi_{t, \gamma}(a+{\i}b)$ as $b$ varies in $(-v_t(a), v_t(a))$ such that $a+{\i}b$ changes within $\Xi_t$. The formula \eqref{eqn:density-Brown-self-adjoint-introduction} describe precisely this observation. Indeed, as $\Phi_{t,\gamma}$ is one-to-one under the assumption of Theorem \ref{thm:main-self-adjoint-introduction}, the set $\Phi_{t,\gamma}(\Xi_t)$ can be parametrized by $a+{\i}b\in\Xi_t$ under the push-forward map $\Phi_{t,\gamma}$. Hence, we can say that the density formula \eqref{eqn:density-Brown-self-adjoint-introduction} depends on only one parameter and is constant in one direction. Theorem \ref{thm:main-self-adjoint-introduction} can be viewed as an analogue result for the free additive convolution in a recent work of Hall--Ho \cite{HallHo2021Brown} concerning free multiplicative Brownian motions. See Remark \ref{remark:twisted-coordinate-x0-self-adjoint} for details. 
	
	 We demonstrate the application of Theorem \ref{thm:main-push-forward-introduction} to $R$-diagonal operators. Let $T$ be an $R$-diagonal that is $*$-free from $\{c_t, g_{t,\gamma}\}$. It is known that the sum of two freely independent $R$-diagonal operators is again an $R$-diagonal operator \cite{NicaSpeicher-Rdiag}. The Brown measure of any $R$-diagonal operator is supported in a single ring \cite{HaagerupLarsen2000}. Let $\lambda_1,\lambda_2$ be inner and outer radii of the support of the Brown measure of $T+c_t$. 
	
	\begin{thmx}
		[See Theorem \ref{thm:Brown-addition-R-diag-elliptic-pushforward}] 
		\label{thm:main-Rdiag-introduction}
		If $T$ is $R$-diagonal, then the support of the Brown measure of $T+g_{t,\gamma}$ is  the deformed single ring where the inner boundary is the circle centered at the origin with radius $\lambda_1$, and the outer boundary is an ellipse rotated by some angle $\alpha$ determined by $\gamma$, centered at the origin, with semi-axes $\lambda_2-\frac{|\gamma|}{\lambda_2}$, and
		$\lambda_2+\frac{|\gamma|}{\lambda_2}$.
		The Brown measure is absolutely continuous and its density is strictly positive in the support. 
		
		The Brown measure of $T+g_{t,\gamma}$ is the push-forward measure of the Brown measure $T+c_t$ by the map $\Phi_{t,\gamma}$. The map sends the circle centered at the origin with radius $r$ to the ellipse whose semi-axes are given by 
		\[
		a(r)=r-|\gamma| m(r)/r, \qquad b(r)=r+|\gamma|m(r)/r,
		\] 
			 where 
		\[
		m(r)={\mu_{T+c_t}(\{z\in\mathbb{C} : |z|\leq r \})}.
		\]
	If $|\gamma|<t$, then both $a(r)$ and $b(r)$ are strictly increasing and the map $\Phi_{t,\gamma}$ is a homeomorphism. 
	\end{thmx}
	

	\subsection{Discussions on methodologies}

	Our approach is based on a Hermitian reduction method and subordination functions. The Hermitian reduction method was already used for the calculation of Brown measure of quasi-nilpotent DT operators in Aagaard--Haagerup's work \cite{HaagerupAagaard2004}. In physics literature, to our knowledge, the method was first used in two independent work \cite{FeinbergZee1997-a} and \cite{Nowak1997-non-Hermitian}. In the work by Jarosz and Nowak \cite{NowakJ2004-preprint, NowakJ2006-addition}, the authors used Hermitian reduction approach to study the support of the Brown measure of  ${x}+{\i}g_t$ for self-adjoint ${x}$, where the method is not mathematically rigorous as written. 
	The Hermitian reduction also appeared in Voiculescu's earlier work \cite{DVV-operator-valued-1992} which serves a motivation to introduce free probability theory with amalgamation. 
	The idea is to study the Brown measure of a non-normal free random variable $x$ by considering the Hermitian matrix
	\[
	    X=  \begin{bmatrix}
	      0  & x\\
	      x^* & 0
	   \end{bmatrix}
	\]
	and 
	 the $2\times 2$ matrix-valued Cauchy transform 
	\begin{equation}\nonumber
	    G_{X}\bigg ( \begin{bmatrix}
	        {\i}\varepsilon & \lambda\\
	        \overline{\lambda} & {\i}\varepsilon
	\end{bmatrix}      \bigg)=\mathbb{E} \bigg( 
	     \begin{bmatrix}
	        {\i}\varepsilon & \lambda-x\\
	        (\lambda-x)^* & {\i}\varepsilon
	\end{bmatrix} ^{-1}    \bigg)
	\end{equation}
	where $\mathbb{E}$ is the entry-wise conditional expectation $\mathbb{E}: M_2(\CM)\rightarrow M_2(\mathbb{C})$. 
	The entries of the  matrix-valued Cauchy transform $G_X$ carry important information for the calculation of the Brown measure. 
	
	The Hermitian reduction method becomes more powerful by combining with subordination functions and this approach was outlined in Belinschi--\'{S}niady--Speicher's work \cite{BSS2018}. In particular, it is showed \cite{BSS2018} that one can iterate certain fixed point equation for subordination functions to approximate boundary values of subordination functions, and then get approximation of Brown measure of an arbitrary polynomial of free random variables by some numerical schemes.
	For the sums ${x}+c_t$ or ${x}+g_{t,\gamma}$, we consider their Hermitian reductions and treat them as the summation of self-adjoint free random variables in the framework of operator-valued free probability. We then use the subordination functions in operator-valued free probability theory to study the Brown measure of ${x}+c_t$ or ${x}+g_{t,\gamma}$. It turns out there are nice formulas for subordination functions and this allows us to obtain explicit Brown measure formulas using subordination functions. 
	
	Our approach extends techniques from Aagaardd--Haagerup \cite{HaagerupAagaard2004}, and Haagerup--Schultz \cite{HaagerupSchultz2007}, and Belinschi--\'{S}niady--Speicher \cite{BSS2018}. We can view the method in Section \ref{subsection:operator-sub} as an operator-valued version of Biane's method used in \cite{Biane1997}.  To our best knowledge,  the Hermitian reduction methods have not be used broadly enough to calculate the explicit Brown measure formula, and existing results mainly focus on a single operator \cite{HaagerupAagaard2004, BSS2018, HaagerupSchultz2007}. The main technical issues when applying this approach are explicit formulas for subordination functions and analyticity of subordination functions and Cauchy transforms on their domains. The present work overcomes these issues and demonstrates that they are also applicable to study the explicit formula of the Brown measure of the sum of two free random variables ${x}+g_{t,\gamma}$. We expect that more applications of these methods are possible.  
	
	The density formula of ${x}+c_t$ for a (unbounded) Gaussian distributed normal operator ${x}$ was established earlier by Bordenave-Caputo-Chafa\"{\i} \cite[Theorem 1.4]{Bordenave-Caputo-Chafai-cpam2014}, where they used Hermitian reduction method and some non-trivial random matrix results. Their techniques can be extended for an arbitrary normal operator.  In particular, Proposition 4.3 in \cite{Bordenave-Caputo-Chafai-cpam2014} is a version of our subordination result in Theorem \ref{thm:sub-X0-gt-formula} for such operators. Some further properties of the Brown measure formula of ${x}+c_t$ was studied in \cite[Section 2]{Bordenave-Capitaine-cpam2016}. Our proof is different from the method used in \cite{Bordenave-Caputo-Chafai-cpam2014} and this allows us to work on non-normal operators. We emphasize that we are interested in the push-forward property between circular case and elliptic case for an arbitrary operator ${x}$. The main results in the present work can be extended to unbounded operators and this was done in a joint work with Belinschi and Yin \cite{BelinschiYinZhong2021Brown}.

	The paper has seven sections. After the Introduction and Preliminaries, in Section 3 we study the sum of a circular operator and a free random variable. We obtain a formula for the Fuglede-Kadison determinant $\Delta( {x}+g_{t,\gamma}-z\unit)$ and $\Delta( {x}+c_t-\lambda\unit)$ using subordination functions. In Section 4, we study the Brown measures of ${x}+c_t$. In Section 5, we show that the Brown measure $\mu_{{x}+g_{t,\gamma}}$ is the push-forward measure of $\mu_{{x}+c_t}$ under the map $\Phi_{t,\gamma}$. In Section 6, we calculate the Brown measure of the sum of a twisted elliptic operator and a self-adjoint operator. In Section 7, we calculate the push-forward map and the Brown measure for the case that ${x}$ is an $R$-diagonal operator. Finally, we calculate explicit formulas for some non-self-adjoint examples in Section 8.

\section{Preliminaries}
\label{Prelim}
\subsection{Free probability and subordination functions}
\label{FreeProb}

We recall the definition of freeness with amalgamation over a subalgebra \cite{SpeicherAMS1998, DVV-operator-valued-1992}. 
An \emph{operator-valued $W^*$-probability space} $(\CA, \mathbb{E}, \CB)$ consists of a von Neumann algebra $\CA$, a unital $*$-subalgebra $\CB\subset \CA$, and  a \emph{conditional expectation}
$\mathbb{E}:\CA\rightarrow \mathcal{B}$. Thus, $\mathbb{E}$ is a linear, unital linear positive map satisfying:
(1) $\mathbb{E}(b)=b$ for all $b\in\CB$, and
(2) $\mathbb{E}(b_1xb_2)=b_1\mathbb{E}(x)b_2$ for all $x\in\CA$, $b_1, b_2\in\CB$. 
Let $(\CA_i)_{i\in I}$ be a family of sublagebras $\CB\subset \CA_i\subset \CA$. We say that $(\CA_i)_{i\in I}$ are \emph{free with amalgamation} over $\CB$ with respect to the conditional expectation $\mathbb{E}$ \emph{(}or free with amalgamation in $(\CA, \mathbb{E}, \CB)) $ if
\[
    \mathbb{E}(x_1x_2\cdots x_n)=0
\]
whenever there are $n\geq 1$ and indexes $i_1, i_2,\cdots, i_n \in I$ such that $i_1\neq i_2$, $i_2\neq i_3$, $\cdots$, $i_{n-1}\neq i_n$, and for $j=1,2,\cdots, n$, we have $x_j\in \CA_{i_j}$ and $\mathbb{E}(x_j)=0$. 

Let $(\CA, \mathbb{E}, \CB)$ be an operator-valued $W^*$-probability space. The elements in $\CA$ are called (noncommuntative) random variables.  We call 
\[
   \mathbb{H}^+(\CB)=\{ b\in\CB: \exists \varepsilon>0, \Im (b)\geq \varepsilon \unit \}
\]
the Siegel upper half-plane of $\CB$, where we use the notation $\Im (b)=\frac{1}{2{\i}}(b-b^*)$.
We set $\mathbb{H}^-(\CB)=\{-b: b\in  \mathbb{H}^+(\CB)\}$. The $B$-valued Cauchy transform $G_X$ of any self-adjoint operator $X\in \CA$ is defined by
\[
  G_X(b)=\mathbb{E}[ (b-X)^{-1} ],  \quad  b\in  \mathbb{H}^+(\CB).
\] 
The $B$-valued Cauchy transform $G_X$ is a map from $\mathbb{H}^+(\CB)$ to $\mathbb{H}^-(\CB)$. 
The Cauchy transform is one-to-one in $\{b\in \mathbb{H}^+(\CB): ||b^{-1}||<\varepsilon \}$ for $\varepsilon$ sufficiently small, and Voiculescu's amalgamated $R$-transform is now defined for $X\in\CA$ by
\[
   R_X(b)=G_X^{\langle -1 \rangle}(b)-b^{-1}
\]
for $b$ being invertible element of $\CB$ suitably close to zero. 

Let $X, Y$ be two self-adjoint operators that are free with amalgamation in $(\mathcal{A}, \mathbb{E},\mathcal{B})$.  The $R$-transform linearizes the free convolution in the sense that if $X, Y$ are free with amalgamation in $(\CA, \mathbb{E}, \CB)$, then 
\[
  R_{X+Y}(b)=R_X(b)+R_Y(b)
\]
for $b$ in some suitable domain. There exist two analytic self-maps $\Omega_1,\Omega_2$
of the upper half-plane $\mathbb{H}^+(\mathcal{B})$
so that
\begin{equation}
(\Omega_1(b)+\Omega_2(b)-b)^{-1}=G_X(\Omega_1(b))=G_Y(\Omega_2(b))=G_{X+Y}(b),
\end{equation}
for all $b\in \mathbb{H}^+(\mathcal {B})$. We refer the reader to \cite{BelinschiTR2018-sub-operator-valued,  Biane1998, DVV-general} for details.

When $\CM$ is a von Neumann algebra, 
 $\CB=\mathbb{C}\unit$ consists of scalar multiples of identity, and $\phi$ is a normal, faithful tracial state on $\CM$. Then the pair $(\CM, \phi)$ replaces the triple $(\CM, \mathbb{C}\unit, \phi)$. We say $(\CM,\phi)$ is a tracial \emph{$W^*$-probability space}. 
For any self-adjoint element $x\in\CM$, let $\mu=\mu_x$ be its spectral measure in $(\CM, \phi)$ determined by
\[
   \phi(f(x))=\int_\mathbb{R}f(u)d\mu_x(u)
\]
for all $f\in C(\sigma(x)$.
The Cauchy transform of $\mu$ (or the Cauchy transform of $x$) can be written as 
$$G_\mu(z)=G_{x}(z) = \phi((z-x)^{-1})= \int_\R\frac{1}{z-u}\:d\mu(u),  \quad z\in\C^+.$$
We also set $F_\mu(z)=F_x(z)=1/G_\mu(z)$. 
The reciprocal Cauchy transform $F_\mu$ maps the upper half plane 
$\mathbb{C}^+$ into itself. 
The $R$-transform of $\mu$ is now an analytic function 
\begin{equation}
\label{Rtransform}
  \mathcal{R}_\mu(z)=G_\mu^{\langle -1 \rangle}(z)-\frac{1}{z}
\end{equation}
where $G_\mu^{\langle -1\rangle}$ denotes the inverse function to $G_\mu$, that is defined in a truncated Stolz angle $\{z\in\C\ : \Im z>\beta, \left\vert \Re z\right\vert<\alpha (\Im z\})$ for some $\alpha, \beta>0$.

Suppose that the self-adjoint random variables $x,y\in \CM$ are freely independent. 
Denote by $\mu_1 $ the spectral measure of $x$,  and $\mu_2$ the spectral measure of $y$, and $\mu_1\boxplus \mu_2$ the \emph{free additive convolution} of $\mu_1$ and $\mu_2$ in the sense that 
$\mu_1\boxplus \mu_2:=\mu_{x+y}$.
The $\mathcal{R}$-transform~\eqref{Rtransform} also linearizes the free additive convolution \cite{BercoviciVoiculescu1993} such that $\mathcal{R}_{\mu\boxplus \nu}(z)=\mathcal{R}_\mu(z)+\mathcal{R}_\nu(z)$ in the domain where all the three $\mathcal{R}$-transforms are defined. 
In this scalar case, there exists a unique pair of analytic maps $\omega_1, \omega_2: \mathbb{C}^+ \to\mathbb{C}^+$ such that
\[
F_{\mu_1\boxplus\mu_2}(z) =  F_{\mu_1}(\omega_1(z)) = F_{\mu_2}(\omega_2(z)) =\omega_1(z)+\omega_2(z)-z
\]
for all $z\in \mathbb{C}^+$.
The above subordination relations can also be written in terms of Cauchy transform. That is, 
$G_{\mu_1\boxplus\mu_2}(z) =  G_{\mu_1}(\omega_1(z)) = G_{\mu_2}(\omega_2(z))$.
The existence of subordination functions leads to many regularity results (see \cite{Belinschi2008, BercoviciVoiculescu1998} and the survey paper \cite[Chapter 6]{survey-2013}). 
The regularity of subordination functions is important in our approach. See Lemma \ref{lemma:s-epsilon-t-limit} for example. 
 For a probability measure $\mu$ on $\mathbb{R}$, denote $H_\mu(z)=F_\mu(z)-z$. In \cite{BB2007new}, Belinschi--Bercovici showed that $\omega_1, \omega_2$ can be obtained from the following fixed point equations
\[
   \omega_1(z)=z+H_{\mu_2}(z+H_{\mu_1}(\omega_1(z))), \qquad
    \omega_1(z)=z+H_{\mu_1}(z+H_{\mu_2}(\omega_1(z))).
\]
Although the subordination functions, in general, cannot be computed explicitly, 
they play a key role in our study by adopting this fixed point approach. See Subsection \ref{section-subordination-scalar} for details.

\subsection{The Brown Measure}
The spectral theorem does not apply to non-normal operators. The Brown measure of an operator in $\mathcal{M}$ was introduced by Brown \cite{Brown1986} and is a natural replacement of the spectral distribution of a non-normal operator.
Given $x\in\CM$, the \emph{Fuglede--Kadison determinant} $\Delta(x)$ \cite{FugledeKadison1952}  of $x$ is defined as 
\[
\Delta(x)=\exp[\phi(\log(|x|))]\in [0,\infty).
\]
The \emph{Brown measure} \cite{Brown1986} of $x$ is then defined to be the distributional Laplacian of the subharmonic function $\log\Delta(x-\lambda)$. That is,
\begin{equation}
\label{eq:brown} \mu_x = \frac{1}{2\pi}\nabla^2_\lambda \log\Delta(x-\lambda)
= \frac{2}{\pi} \frac{\partial}{\partial
\lambda}  \frac{\partial}{\partial \bar{\lambda}}  \log
\Delta(x-\lambda\unit).
\end{equation}
In fact, $\mu_x$ is a
probability measure supported on a subset of the spectrum of
$x$. When $\mathcal{M}=M_n(\C)$ and $\phi$ is the normalized trace on $M_n(\mathbb{C})$,  for $x\in M_n(\C)$, we have
$$\log\Delta(x-\lambda) = \log |\det (x-\lambda I)|^{1/n}=\frac{1}{n}\sum_{i=1}^n\log |\lambda-\lambda_i|,$$
where $\lambda_1,\cdots, \lambda_n$ are the eigenvalues of $x$. 
Hence the eigenvalue distribution of $x$ can be recovered by taking the distributional Lapalacian
\[
  \frac{1}{n}\left( \delta_{\lambda_1}+\cdots \delta_{\lambda_n}\right)=\frac{1}{2\pi}\nabla^2_\lambda \log |\det (x-\lambda I)|^{1/n}.
\]

It is useful to consider the regularized function
\[
    S(x,\lambda,\varepsilon)=\phi(\log((x-\lambda\unit)^*(x-\lambda\unit)+\varepsilon^2\unit)), 
     \qquad \varepsilon>0.
\]
Then, by the tracial property of $\phi$, we have $\log\Delta(x-\lambda)=\frac{1}{2}\lim_{\varepsilon\rightarrow 0} S(x,\lambda,\varepsilon)$, and the Brown measure is calculated as
\[
   \mu_x=\frac{1}{4\pi}\nabla^2_\lambda \left(\lim_{\varepsilon\rightarrow 0} \phi(\log((x-\lambda\unit)^*(x-\lambda\unit)+\varepsilon^2\unit))\right).
\]
It is known \cite{HaagerupSchultz2007, Larsen1999-thesis} that, for any $\varepsilon>0$, the function $\lambda\mapsto S(x,\lambda,\varepsilon)$ is subharmonic and its Riesz measure is a probability measure, defined by
\[
    \mu_{x,\varepsilon}=\frac{1}{4\pi}\nabla^2_\lambda S(x,\lambda,\varepsilon).
\]
Moreover, $\mu_{x,\varepsilon}$ converges to $\mu_x$ weakly as $\varepsilon$ tends to zero. 
This regularization process makes the calculation of general operator in $\mathcal{M}$ more tractable. For the summation ${x}+c_t$ and ${x}+g_{t,\gamma}$, we are able to identify the domain $\Xi_t$, so that $ \log\Delta({{x}+c_t}-\lambda)$ is real analytic for any $\lambda\in\Xi_t$, and $\log\Delta({{x}+g_{t,\gamma}}-z)$ is real analytic for any $z\in \Phi_{t,\gamma}(\Xi_t)$ if $\Phi_{t, \gamma}$ is non-singular. Hence, the Brown measures in this paper can be calculated in classic sense. See Theorem \ref{thm:main-FK-det-ct-0} for details.

\subsection{Hermitian reduction method for the sum of two free random variables}
\label{subsection:2.3Hermitian}
Let $(\CM, \phi)$ be a tracial $W^*$-probability space. 
We equip
the algebra $M_2(\CM)$, the $2\times 2$ matrices with entries from $\CM$, with the conditional expectation
$\E:M_2(\CM)\rightarrow M_2(\C)$ given by \begin{equation}
\label{eq:conditional} \E
\begin{bmatrix} a_{11} &
a_{12} \\ a_{21} & a_{22} \end{bmatrix} = \begin{bmatrix} \phi(a_{11}) & \phi(a_{12}) \\
\phi(a_{21}) & \phi(a_{22}) \end{bmatrix}. 
\end{equation}
Then the triple  $(M_2(\CM), \mathbb{E}, M_2(\mathbb{C}))$ is a operator-valued $W^*$-probability space. 
Given $x\in\CM$,
let
\begin{equation}
\label{eq:X} \x=
\begin{bmatrix} 0 & x \\ x^\ast & 0 \end{bmatrix} \in M_2(\CM),
\end{equation}
which a self-adjoint element in $M_2(\CM).$
For $\lambda\in\mathbb{C}$ and $\delta\in\mathbb{C}$ with $\Im\delta>0$, the element
\begin{equation}
\label{defn:Theta}
\Theta
(\lambda, \delta) =  \begin{bmatrix} \delta & \lambda \\
\bar{\lambda}  & \delta \end{bmatrix} \in  M_2(\C) 
\end{equation}
belongs to the domain of the $M_2(\C)$-valued Cauchy $G_X$.
We now record that
\begin{align}
&\label{inverse}
 \begin{bmatrix}
  a & b \\
  c & d
 \end{bmatrix}^{-1}
 =\begin{bmatrix}
  d(ad-bc)^{-1} & -b(ad-cb)^{-1}\\
  -c(ad-bc)^{-1} & a(ad-cb)^{-1}
 \end{bmatrix}
\end{align}
where $a,d\in\mathbb{C}$ and $b,c\in\mathcal{M}$ such that $ad-bc$ is invertible (which is equivalent to $ad-cb$ is invertible).
We then have
\begin{equation}
	\begin{aligned}
	\label{eqn:inverse-lambda-X}
	&(\Theta
	(\lambda,{\i}\varepsilon) -\x)^{-1}\\
	& =  \begin{bmatrix}
	-{\i}\varepsilon \big((\lambda-x)(\lambda-x)^\ast+\varepsilon^2\big)^{-1}
	& (\lambda-x)\big((\lambda-x)^\ast (\lambda-x)+\varepsilon^2\big)^{-1} \\
	(\lambda-x)^\ast \big(
	(\lambda-x)(\lambda-x)^\ast+\varepsilon^2\big)^{-1}  & -{\i} \varepsilon
	\big( (\lambda-x)^\ast(\lambda-x)+\varepsilon^2 \big)^{-1}
	\end{bmatrix}.
	\end{aligned}
\end{equation}
and 
\begin{equation} 
\label{eqn-Cauchy-X-Theta} 
G_X
(\Theta
(\lambda,{\i}\varepsilon) ) 
=
\E \big( (\Theta
(\lambda, {\i}\varepsilon) - \x)^{-1}\big)
=
\begin{bmatrix}
\displaystyle g_{X,11}(\lambda, \varepsilon) &
\displaystyle g_{X,12}(\lambda, \varepsilon)\\
\displaystyle g_{X,21}(\lambda, \varepsilon) &
\displaystyle  g_{X,22}(\lambda, \varepsilon)
\end{bmatrix}.
\end{equation} 
where 
\begin{eqnarray*}
g_{X,11}(\lambda, \varepsilon) & = & -{\i}\varepsilon\phi\left(\big((\lambda-x)(\lambda-x)^*+\varepsilon^2
\big)^{-1}\right)\\
g_{X,12}(\lambda, \varepsilon) &=& \phi\left((\lambda-x)\big((\lambda-x)^* (\lambda-x)+\varepsilon^2
\big)^{-1}\right) \\
g_{X,21}(\lambda, \varepsilon)&=&\phi\left((\lambda-x)^*\big((\lambda-x)(\lambda-x)^*+
\varepsilon^2\big)^{-1}\right)\\
g_{X,22}(\lambda, \varepsilon) & = & -{\i}\varepsilon\phi\left(\big( (\lambda-x)^*(\lambda-x)+\varepsilon^2 
\big)^{-1}\right).
\end{eqnarray*}

We note that by the tracial property of $\phi$, we have 
\begin{equation}
 \label{eqn-G-entries-symmetry}
  g_{X+Y,11}(\lambda, \varepsilon)=g_{X+Y,22}(\lambda, \varepsilon),
  \qquad 
  g_{X,21}(\lambda, \varepsilon)=\overline{g_{X,12}(\lambda, \varepsilon)}
\end{equation}
We observe that entries of the Cauchy transform \eqref{eqn-Cauchy-X-Theta} have symmetry similar to the matrix $\Theta(\lambda,\varepsilon)$. This can be explained as follows. Define the map $J: M_2(\mathbb{C})\rightarrow M_2(\mathbb{C})$ by 
\[
    b\mapsto J(b)=-b^*.
\]
Then we have $G_X(J(b))=J(G_X(b))$. Notice that $J(\Theta(\lambda,{\i}\varepsilon))=\Theta(\lambda,{\i}\varepsilon)$ and hence $G_X(\Theta(\lambda,{\i}\varepsilon))$ has symmetric property \eqref{eqn-G-entries-symmetry}.

Equations \eqref{eqn-Cauchy-X-Theta} show that two
diagonal entries of $G_X
(\Theta
(\lambda, {\i}\varepsilon) ) $ carry important
information to calculate the Brown measure of $x$. 
Let $x$ and $y$ be two $*$-free random variables.  We have
to understand the $\mathcal M_2(\mathbb C)$-valued distribution of
$$\begin{bmatrix} 0 & x+y \\ (x+y)^* & 0 \end{bmatrix}=X+Y=
\begin{bmatrix} 0 & x \\ x^* & 0 \end{bmatrix}+
\begin{bmatrix} 0 & y \\ y^* & 0 \end{bmatrix}$$
in terms of the $\mathcal M_2(\mathbb C)$-valued distributions of $X$ and of $Y$. Note that $X$ and $Y$ are free over $\mathcal M_2(\mathbb C)$.
The subordination functions in this context are two analytic self-maps  $\Omega_1,\Omega_2$
of the upper half-plane $\mathbb{H}^+(\mathcal M_2(\mathbb C))$
so that
\begin{equation}\label{eqn:subord}
(\Omega_1(b)+\Omega_2(b)-b)^{-1}=G_X(\Omega_1(b))=G_Y(\Omega_2(b))=G_{X+Y}(b),
\end{equation}
for every $b\in \mathbb{H}^+(\mathcal M_2(\mathbb C))$. We shall be concerned with $b=\Theta(\lambda, {\i}\varepsilon)$. Indeed, 
we have, by \eqref{eqn-Cauchy-X-Theta}, 
$$
G_{X+Y}(\Theta(\lambda,{\i}\varepsilon))=\begin{bmatrix}
\displaystyle g_{X+Y,11}(\lambda, \varepsilon) &
\displaystyle g_{X+Y,12}(\lambda, \varepsilon)\\
\displaystyle g_{X+Y,21}(\lambda, \varepsilon) &
\displaystyle  g_{X+Y,22}(\lambda, \varepsilon)
\end{bmatrix}
$$
where 
\begin{eqnarray*}
g_{X+Y,11}(\lambda, \varepsilon) & = & -{\i}\varepsilon\phi\left(\big((\lambda-x-y)(\lambda-x-y)^*+\varepsilon^2
\big)^{-1}\right)\\
g_{X+Y,12}(\lambda, \varepsilon) &=& \phi\left((\lambda-x-y)\big((\lambda-x-y)^* (\lambda-x-y)+\varepsilon^2
\big)^{-1}\right) \\
g_{X+Y,21}(\lambda, \varepsilon)&=&\phi\left((\lambda-x-y)^*\big((\lambda-x-y)(\lambda-x-y)^*+
\varepsilon^2\big)^{-1}\right)\\
g_{X+Y,22}(\lambda, \varepsilon) & = & -{\i}\varepsilon\phi\left(\big( (\lambda-x-y)^*(\lambda-x-y)+\varepsilon^2 
\big)^{-1}\right).
\end{eqnarray*}
The idea of calculating the Brown measure of $x+y$ is to separate the information of $X$ and $Y$ in some tractable way. We achieve this by using subordination functions \eqref{eqn:subord}.

\subsection{The elliptic operator and the operator-valued semicircular element}
\label{section:elliptic}

In the tracial $W^*$-probability space $(\CM, \phi)$, 
for $t>0$, 
Voiculescu's circular operator with variance $t$, denoted by $c_t$, is defined as
\[
  c_t=\frac{1}{\sqrt{2}}(s_t+{\i}s'_t)
\]
where $\{ s_t, s_t' \}$ is a free semicircular family and each of them has variance $t$. 

Let $t>0$ and $\gamma\in \mathbb{C}$ such that $|\gamma|\leq t$. 
The \emph{twisted elliptic operator} operator $g_{t,\gamma}$ can be constructed as follows. 
Let $\{s_{t_1},  s_{t_2}\}$ be semicircular operators with zero expectation and variance $t_1, t_2$ respectively such that $\{s_{t_1},  s_{t_2}\}$ are freely independent.  For $\theta\in [0,2\pi]$, consider the operator $y_{t_1, t_2,\theta}=e^{{\i}\theta}(s_{t_1}+{\i}s_{t_2})$, 
by choosing $t_1, t_2$ such that $t_1+t_2=t$, $t_1-t_2=|\gamma|$ and $e^{{\i}2\theta}=\gamma/|\gamma|$, we can check directly that $g_{t,\gamma}$ and $y_{t_1,t_2, \theta}$ have the same $*$-distribution, whose only nonzero free cumulants are given by
\[
  \kappa(y, y^*)=\kappa(y^*,y)=t, \qquad \kappa(y, y)=\gamma, \qquad
    \kappa(x^*, x^*)=\overline{\gamma},
\]
where $y=g_{t,\gamma}$. The operator $g_{t,\gamma}$ include the following operators as special cases:
(i) if $\gamma=0$, $y$ is a circular operator with variance $t$; (ii) if $\gamma=t$, $y$ is a semicircular operator $s_t$ with variance $t$; (iii) if $\gamma=-t$, then $y$ has the distribution as ${\i}s_t$; (iv) if $\gamma\in[-t,t]$, then $y$ has the same distribution as an elliptic operator.

In the operator-valued $W^*$-probability space $(\mathcal{A}, \mathbb{E}, \mathcal{B})$,  
following Voiculescu \cite{DVV-operator-valued-1992} and Speicher \cite{SpeicherAMS1998}, we say $Y\in \CA$ is \emph{$\mathcal{B}$-Gaussian} or \emph{an operator-valued semicircular element} if and only if the $R$-transform has a particular simple form
\begin{equation}
\label{eqn:R-transform-B-Gaussian}
  R_Y(b)=\mathbb{E}_\mathcal{B}(YbY).
\end{equation}
Condition \eqref{eqn:R-transform-B-Gaussian} says that only $\mathcal{B}$-cumulants of length two survive. Note that a linear combination of two operator-valued semicircular elements in $(\mathcal{A}, \mathbb{E}, \mathcal{B})$ is again an
operator-valued semicircular element.

The following result is a special case of \cite[Example 19 in Section 9.4]{MingoSpeicherBook}. One can also deduce it from a general formula about a relation between matrix-valued and scalar-valued free cumulants in \cite[Theorem 6.2]{NicaSS2002-Rcyclic} (or \cite[Proposition 13 in Section 9.3]{MingoSpeicherBook}).
\begin{proposition}
\label{prop:semi-circular-operatoro-valued}
Let $g_{t,\gamma}$ be a twisted elliptic operator with parameters $t, \gamma$ in the tracial $W^*$-probability space $(\CM, \phi)$ and $\lambda\in\mathbb{C}$. Denote
\[
  Y=\begin{bmatrix}
  0 & g_{t,\gamma}\\
  g^*_{t,\gamma} &0
\end{bmatrix}\in M_2(\CM).
\]
Then $Y$ is an operator-valued semicircular element in the operator-valued $W^*$-probability space 
$(M_2(\CM), \mathbb{E}, M_2(\mathbb{C}))$. 
\end{proposition}


\section{The Fuglede-Kadison determinant and subordination functions}
\label{section-FK-det-subordination}
Given $t>0$ and $\gamma\in \mathbb{C}$ such that $|\gamma|\leq t$,
let $y=g_{t,\gamma}$ be a twisted elliptic operator and let ${x}$ be a random variable that is $*$-free from $y$ in the $W^*$-probability space $(\mathcal{M}, \phi)$.
In this section, we study the Fuglede-Kadison determinant $\Delta({x}+y-\lambda\unit)$
for $\lambda\in\mathbb{C}$.
We denote 
\begin{equation}
\label{eqn:X0-Ct}
X=
\begin{bmatrix}
0 & {x}\\
{x}^* & 0
\end{bmatrix}, 
\qquad 
Y=
\begin{bmatrix}
0 & y\\
y^* & 0
\end{bmatrix}.
\end{equation}
Note that $X$ and $Y$ are free over $\mathcal M_2(\mathbb C)$. 
There exist two analytic self-maps $\Omega_1, \Omega_2$
of upper half-plane $\mathbb{H}^+(\mathcal{M}_2(\mathbb{C}))$ of $\mathcal M_2(\mathbb C)$
such that
\begin{equation}\label{eqn:subord-X0-Ct}
(\Omega_1(b)+\Omega_2(b)-b)^{-1}=G_X(\Omega_1(b))=G_Y(\Omega_2(b))=G_{X+Y}(b),
\end{equation}
for all $b\in\mathcal M_2(\mathbb C)$ with $\Im b>0$.
We choose
\[
\Theta(\lambda,{\i}\varepsilon)=
\begin{bmatrix}
{\i}\varepsilon & \lambda\\
\overline{\lambda} & {\i}\varepsilon
\end{bmatrix}
\]
where $\varepsilon>0$ and $\lambda\in \mathbb{C}$. Our strategy is to find more explicit formulations for subordination functions $\Omega_1, \Omega_2$.

\subsection{Subordination functions in free convolution with a semicircular distribution}
\label{section-subordination-scalar}

Recall that 
$\lambda\in\Xi_t$ defined in \eqref{defn:Xi-t} if the following condition holds
\begin{equation}
\label{eqn:3.3000-in-proof}
\phi \left[\big( ({x}-\lambda\unit)^*({x}-\lambda\unit )\big)^{-1}\right]>\frac{1}{t}.
\end{equation}
For any $\lambda\in\Xi_t$, let $w=w(0;\lambda,t)$ be the unique positive function of $\lambda$ such that
\begin{equation}
\label{eqn:3.4000-in-proof}
\phi \left[\big( ({x}-\lambda\unit)^*({x}-\lambda\unit ) +w^2\unit\big)^{-1}\right]=\frac{1}{t}.
\end{equation}
In this section, we show that $w(0;\lambda,t)$ is the imaginary part of a subordination function and it is a real analytic function of $\lambda$ as long as \eqref{eqn:3.3000-in-proof} holds. 

\begin{proposition}
	\label{prop:Xi-t-open-set}
	The set $\Xi_t$ is bounded and open for any $t>0$. 
\end{proposition}
\begin{proof}
	For $|\lambda|$ large enough, $\lambda\notin \sigma({x})$ and 
	\[
	\lim_{\lambda\rightarrow\infty} \phi \left[\big( ({x}-\lambda\unit)^*({x}-\lambda\unit )\big)^{-1}\right]=0.
	\]
	Hence, the set $\Xi_t$ is bounded. 
	For any $\varepsilon\geq 0$, denote the function $f_\varepsilon$ of $\lambda$ by
	\[f_\varepsilon(\lambda)=\phi\left[\big( ({x}-\lambda\unit)^*({x}-\lambda\unit) +\varepsilon\unit \big)^{-1}\right].\] 
	For any $\varepsilon>0$, the function $f_\varepsilon(\lambda)$ is a continuous function of $\lambda$. Observe that $f_0$ is the limit of the increasing sequence of $f_\varepsilon$, hence it is lower semi-continuous. The set $\Xi_t$ can be rewritten as $\Xi_t=\{\lambda: f_0(\lambda)>1/t \}$ and therefore $\Xi_t$ is open for any $t>0$. 
\end{proof}

For $\mu\in\Prob(\R, \B)$ let $\tilde\mu$ denote the
{\it symmetrization} of $\mu$. That is, $\tilde\mu\in \Prob(\R, \B)$ is given
by
\[
\tilde\mu(B) = \textstyle{\frac12}(\mu(B)+\mu(-B)), \qquad (B\in\B).
\]
	For a probability measure $\mu$ on $\mathbb{R}$, let 
	\begin{equation}\label{defn:function-h}
	 h_\mu(s)=\int_\mathbb{R}\frac{s}{s^2+u^2}d\mu(u), \qquad s>0. 
	\end{equation}
\begin{proposition}
	\label{prop:sub-identical-symmetric-measure}
	Let $\mu_1=\tilde{\mu}_{|{x}-\lambda\unit|}$ and $\mu_2$ be the semicircular distribution with variance $t$. Denote $\mu=\mu_1\boxplus \mu_2$.
	Let $\omega_1^{(\lambda)}, \omega_2^{(\lambda)}$ be subordination functions such that
	\begin{equation}\label{eqn-subordination-F}
	F_{{\mu}}(z)= F_{{\mu}_{1}}(\omega_1^{(\lambda)}(z))=
	F_{\mu_2}(\omega_2^{(\lambda)}(z)).
	\end{equation}
	For any $\varepsilon>0$, set $W_1(\varepsilon)=\Im\omega_1^{(\lambda)}({\i}\varepsilon)$, then $W_1(\varepsilon)$ satisfies the identity
	\begin{equation}
	\label{eqn:condition-s-general-t}
	h_{\mu_1}(W_1(\varepsilon))=\int_0^\infty\frac{W_1(\varepsilon)}{W_1(\varepsilon)^2+u^2}d\mu_{|{x}-\lambda\unit|}(u)
	=\frac{W_1(\varepsilon)-\varepsilon}{t}.
	\end{equation}
\end{proposition}
\begin{proof}
	For simplicity, we denote $\omega_j=\omega_j^{(\lambda)}$ for $j=1, 2$. 
	We have $\omega_1(z)+\omega_2(z)=z+F_{{\mu}}(z)$.
	Notice that, by symmetry of $\mu_1, \mu_2$ and $\mu$, we have 
	\begin{equation}\label{eqn:G_tilde_h_lambda}
	G_\mu({\i}\varepsilon)=\int_\mathbb{R}\frac{1}{{\i}\varepsilon-u}d\mu(u)=
	\int_\mathbb{R}\frac{-{\i}\varepsilon}{\varepsilon^2+u^2}d\mu(u)=-{\i}h_\mu(\varepsilon).
	\end{equation}
	and similarly
	\begin{equation}\label{eqn:G_tilde_h}
	G_{\mu_1}({\i}\varepsilon)=-{\i}h_{\mu_1}(\varepsilon), \quad
	G_{\mu_2}({\i}\varepsilon)=-{\i}h_{\mu_2}(\varepsilon).
	\end{equation}
	Let $H_1(z)=F_{\mu_1}(z)-z$ and $H_2(z)=F_{\mu_2}(z)-z$.
	Then $\omega_1$ satisfy the following fixed point equation
	\[
	\omega_1(z)=z+H_2(z+H_1(\omega_1(z))).
	\]
	Since $\mu_2$ is the semicircle distribution with variance $t$, its Cauchy transform satisfies
	\[
	\frac{1}{G_{\mu_2}(z)}+t G_{\mu_2}(z)=z,
	\]
	Hence $H_2(z)=-tG_{\mu_2}(z)$ and  the fixed point equation reads
	\[
	\omega_1(z)-z=-tG_{\mu_2}(z+H_1(\omega_1(z))).
	\]
	It is clear that $\omega_1({\i}\varepsilon)$ is pure imaginary. We then have
	\begin{equation}
	\label{eqn:3.13-in-proof}
	{\i}W_1(\varepsilon)-{\i}\varepsilon=-tG_{\mu_2}({\i}\varepsilon+H_1({\i}W_1(\varepsilon)))
	\end{equation}
	Note that $z+H_1(\omega_1(z))=\omega_2(z)$. Hence 
	\[
	  G_{\mu_2}({\i}\varepsilon+H_1({\i}W_1(\varepsilon)))
	    =G_{\mu_1}(\i W_1(\varepsilon)).
	\]
	Using $G_{\mu_1}(\i \varepsilon)=-\i h_{\mu_1}(\varepsilon)$, then Equation \eqref{eqn:3.13-in-proof} implies 
	\begin{equation*}
	h_{\mu_1}(W_1(\varepsilon))=\frac{W_1(\varepsilon)-\varepsilon}{t}
	\end{equation*}
which yields \eqref{eqn:condition-s-general-t}.	
\end{proof}

Note that the defining identity \eqref{eqn:condition-s-general-t} for $W_1(0)$ reduces to the defining identity \eqref{eqn:3.4000-in-proof} for $w(0;\lambda,t)$ provided that $W_1(0)>0$. This suggess that $w(0;\lambda,t)$ must be the boundary value of 
the subordination function $\omega_1^{(\lambda)}$ parametrized by $\lambda\in\mathbb{C}$ that was defined in Proposition \ref{prop:sub-identical-symmetric-measure}.
We put
	\begin{equation}\label{defn:function-k}
	k(s,\varepsilon)=\frac{s-\varepsilon}{h_\mu(s)}, \qquad s>0, \quad \varepsilon>0.
	\end{equation}
\begin{lemma}
	\label{lemma:monotonicity-k-s-epsilon}
	Given a probability measure $\mu$ on $\mathbb{R}$, the function $k$ is an analytic function on $(0,\infty)\times(0,\infty)$. Moreover, for $\varepsilon>0$, the map
	$s\mapsto k(s,\varepsilon)$ is a strictly increasing bijection of $(\varepsilon,\infty)$ onto $(0, \infty)$, and for $\varepsilon=0$, the map $s\mapsto k(s,0)$ is a strictly increasing bijection of $(0,\infty)$ onto $\left(\lambda_1(\mu)^2,\infty\right)$,
	where 
	\[
	\lambda_1(\mu)^2=\left(\int_\mathbb{R} \frac{1}{u^2}d\mu(u)\right)^{-1}.
	\]
\end{lemma}
\begin{proof}
	It is clear that $k$ is analytic. Moreover, for $\varepsilon>0$,
	\begin{equation}
	\label{eqn:k-s-epsilon-product}
	k(s,\varepsilon)=\frac{s-\varepsilon}{s}\left(\int_\mathbb{R} \frac{1}{s^2+u^2}d\mu(u)\right)^{-1}, 
	\end{equation}
	which is a product two increasing and positive functions of $s$ on $(\varepsilon, \infty)$.
	The monotonicity properties of $s\mapsto k(s,\varepsilon)$ follows for $\varepsilon\geq 0$.
\end{proof}

\begin{definition}
	\label{defn:s-epsilon-t-00}
	For $\lambda\in\mathbb{C}$, set $\mu=\mu_{|{x}-\lambda\unit|}$. Let $h_\mu, k$ as in \eqref{defn:function-h} and \eqref{defn:function-k} respectively. For $\varepsilon, t\in (0,\infty)$, let $w(\varepsilon;\lambda,t)$ denote the unique solution $w\in (\varepsilon, \infty)$  to the equation $k(w,\varepsilon)=t$ following Lemma \ref{lemma:monotonicity-k-s-epsilon}. In other words, $w=w(\varepsilon;\lambda,t)\in (\varepsilon,\infty)$ is the unique solution of the equation
	\begin{equation}
	\label{eqn:condition-s-general-t-v2}
	\int_0^\infty\frac{w}{w^2+u^2}d\mu_{|{x}-\lambda\unit|}(u)
	=\frac{w-\varepsilon}{t}.
	\end{equation}
	
	Note that $\lambda\in\Xi_t$ if and only if $t\in (\lambda_1(\mu)^2,\infty)$. For $\lambda\in\Xi_t$, let $w(0;\lambda,t)$ be the unique solution $w\in (0,\infty)$ to the equation $k(w, 0)=t$, which is equivalent to 
	\begin{equation}
	\label{eqn:condition-s-0-t}
	\int_0^\infty \frac{1}{w^2+u^2}d\mu_{|{x}-\lambda\unit|}(u)=\frac{1}{t},
	\end{equation}
	and can be rewritten as
	\begin{equation}
	\label{eqn:w-epsiton-0-identity}
	\phi\left[\big( ({x}-\lambda\unit)^*({x}-\lambda\unit) +w(0; \lambda,t)^2\unit \big)^{-1}\right]=\frac{1}{t}.
	\end{equation}
	For $\lambda\in \mathbb{C}\backslash\Xi_t$, the operator $(x-\lambda)^*(x-\lambda)$ exists in $L^1(\mathcal{M},\phi)$, we set $w(0;\lambda,t)=0$. 
\end{definition}

\subsection{Convergence of subordination functions}
By Proposition \ref{prop:sub-identical-symmetric-measure} and Definition \ref{defn:s-epsilon-t-00}, we can view $w(\varepsilon;\lambda,t)$ as a family of scalar-valued subordination functions with parameter $\lambda\in\mathbb{C}$. 
In this section, we study the convergence of $w(\varepsilon;\lambda,t)$ to $w(0;\lambda,t)$ as $\varepsilon$ tends to zero.

\begin{lemma}
	\label{lemma:s-epsilon-t-limit}
	The function $(\varepsilon, t)\mapsto w(\varepsilon;\lambda,t)$ is real analytic in $(0,\infty)\times (0,\infty)$.
	The function $\lambda\mapsto w(0;\lambda,t)$ is a continuous function on $\mathbb{C}$. 
	Moreover, 
	\begin{align*}
	\lim_{\varepsilon\rightarrow 0^+}w(\varepsilon;\lambda,t)
	=w(0;\lambda,t).
	\end{align*}
	In addition, for $\lambda=a+{\i}b\in\Xi_t$, the function
	$(a,b)\mapsto w(0; \lambda,t)$ is real analytic.
\end{lemma}
\begin{proof}
	Let $\Omega=\{ (s,\varepsilon)\in\mathbb{R}^2: 0<\varepsilon<s \}$. Consider the analytic function
	\[
	F(s,\varepsilon)=( k(s,\varepsilon), \varepsilon), \quad (s,\varepsilon)\in \Omega. 
	\]
	By Lemma \ref{lemma:monotonicity-k-s-epsilon}, $F$ is one-to-one map of $\Omega$ onto $(0,\infty)\times (0,\infty)$. Moreover, its inverse function $F^{-1}: (0,\infty)\times (0,\infty)\rightarrow \Omega$ is given by 
	\[
	F^{-1}(t,\varepsilon)=(w(\varepsilon;\lambda,t), \varepsilon), \quad t, \varepsilon>0.
	\] 
	We now calculate the determinant of Jacobian of $F$ as 
	\[
	\det(\mathrm{J}(F))(s, \varepsilon)=\frac{\partial}{\partial s}k(s,\varepsilon)>0
	\]
	by \eqref{eqn:k-s-epsilon-product}.
	
		We denote $h(\varepsilon,\lambda)=h_{\mu_{|\lambda-{x}|}}(\varepsilon)$ and observe that
	\[
	h(\varepsilon,\lambda)=\frac{1}{2}\int_0^\infty \left( \frac{1}{\varepsilon+{\i}u}+\frac{1}{\varepsilon-{\i}u} \right) d\mu_{|\lambda-{x}|}(u)=-\Im G_{\tilde{\mu}_{|\lambda-{x}|}}(i\varepsilon).
	\]
	It follows that $h$ and also $k$ is analytic in $\varepsilon$. 
	 Hence, $F^{-1}$ is analytic on $(0,\infty)\times (0,\infty)$. Consequently, $w(\varepsilon;\lambda,t)$ is analytic on $(0,\infty)\times (0,\infty)$.
	
	Following notations in Lemma \ref{lemma:monotonicity-k-s-epsilon}, we set   $\mu=\mu_{|{x}-\lambda\unit|}$ and $\lambda_1(\mu)^2=\left(\int_0^\infty \frac{1}{u^2}d\mu(u)\right)^{-1}$. Recall that $\lambda\in\Xi_t$ if and only if $t>\lambda_1(\mu)^2$. 
	For $\lambda\in \Xi_t$, then $t\in (\lambda_1(\mu)^2,\infty)$ and $w(0;\lambda,t)>0$ by Definition \ref{defn:s-epsilon-t-00}.  Since 
	\[
	\frac{\partial}{\partial s}k(s,0)=\frac{\partial}{\partial s} \frac{s}{h_\mu(s)} >0
	\]
	for $s>0$, 
	where $h_\mu$ is defined in Lemma \ref{lemma:monotonicity-k-s-epsilon}.
	It follows that $F$ is also analytic in some neighborhood $U_0$ of $(w(0; \lambda,t),0)$
	and $F$ has an analytic inverse $F^{-1}$ in a neighborhood $V_0$ of $F(s,0)=(t,0)$. Now
	\begin{equation}
	\lim_{\varepsilon\rightarrow 0^+} F^{-1}(t, \varepsilon)=F^{-1}(t, 0)=(w(0; \lambda,t), 0). 
	\end{equation}
	Hence, for $\lambda\in\Xi_t$, 
	\[
	\lim_{\varepsilon\rightarrow 0^+}w(\varepsilon; \lambda,t)=w(0; \lambda,t). 
	\]
	Moreover, the function $w(0; \lambda,t)$ is a real analytic function of $(a, b)$ where $\lambda=a+{\i}b$. 
	
	We next study convergence for $\lambda\in\mathbb{C}\backslash\Xi_t$. 
	Note that for fixed $\varepsilon>0$, the map $t\mapsto w(\varepsilon;\lambda,t)$ is an increasing function of $t$, which can be verified directly from Lemma \ref{lemma:monotonicity-k-s-epsilon} and Definition \ref{defn:s-epsilon-t-00}. Hence, for $\lambda\in\mathbb{C}\backslash\Xi_t$, then $t\leq \lambda_1(\mu)^2$, and for 
	any $t'>\lambda_1(\mu)^2$ we have 
	\[
	\limsup_{\varepsilon\rightarrow 0^+}w(\varepsilon;\lambda,t)\leq 
	\limsup_{\varepsilon\rightarrow 0^+}w(\varepsilon; \lambda, t')=w(0;\lambda,t').
	\]
	But $t'\mapsto w(0;\lambda,t')$ is a bijection of $(\lambda_1(\mu)^2,\infty)$ onto $(0,\infty)$. It then follows that 
	\[
	\lim_{\varepsilon\rightarrow 0^+}w(\varepsilon;\lambda,t)=0
	\]
	whenever $\lambda\in\mathbb{C}\backslash\Xi_t$. 
	
	Recall that $w(0;\lambda,t)=0$ for $\lambda\in\mathbb{C}\backslash\Xi_t$ and $w(0;\lambda,t)>0$ in the open set $\Xi_t$. 
	Hence, to show that $\lambda\mapsto w(0;\lambda,t)$ is a continuous function in $\mathbb{C}$, it remains to show that for any $\lambda_0\in \mathbb{C}\backslash\Xi_t$ and a sequence $\{ \lambda_n \}\subset \Xi_t$ converging to $\lambda_0$, we have 
	\[
	\lim_{n\rightarrow \infty} w(0;\lambda_n, t)=0. 
	\]
	Suppose this is not true. By dropping to a subsequence if necessary, we may assume that there exists $\delta>0$ such that for all $n$, $w(0;\lambda_n,t)>\delta$. In this case, we have 
	\[
	\int_0^\infty\frac{1}{w(0;\lambda_n,t)^2 +u^2}d\mu_{|{x}-\lambda_n|}(u)=\frac{1}{t}
	\]
	which yields
	\begin{equation}
	\label{eqn:314-in-proof}
	\int_0^\infty\frac{1}{\delta^2 +u^2}d\mu_{|{x}-\lambda_n|}(u)>\frac{1}{t}
	\end{equation}
	Hence we have 
	\[
	\phi( ( |{x}-\lambda_0|^2+\delta^2)^{-1} )=\lim_{n\rightarrow\infty} \phi( ( |{x}-\lambda_n|^2+\delta^2)^{-1} )\geq\frac{1}{t},
	\]
	which implies that $w(0;\lambda_0,t)\geq \delta$. This contradicts to our choice $\lambda_0\in\mathbb{C}\backslash \Xi_t$. Therefore, $\lambda\mapsto w(0;\lambda,t)$ is a continuous function in $\mathbb{C}$. 
\end{proof}

\begin{lemma}
	\label{lemma:sub-w-uniform-convergence}
	Fix $t>0$, then
	the functions $\lambda\mapsto w(\varepsilon;\lambda,t)$ converge to the function $\lambda\mapsto w(0;\lambda,t)$ uniformly on $\mathbb{C}$ as $\varepsilon$ tends to zero. 
\end{lemma}
\begin{proof}
	We denote $w_1=w(\varepsilon_1; \lambda, t), w_2=w(\varepsilon_2; \lambda, t)$. 
	We observe that $\varepsilon_1< w_1<w_2$ for any $0<\varepsilon_1<\varepsilon_2$.
	Indeed, $w_1$ is the unique solution of 
	\[
	\int_0^\infty\frac{w_1}{w_1^2+u^2}d\mu_{|{x}-\lambda\unit|}(u)
	=\frac{w_1-\varepsilon_1}{t},
	\]
	which yields
	\[
	\int_0^\infty\frac{w_1}{w_1^2+u^2}d\mu_{|{x}-\lambda\unit|}(u)
	>\frac{w_1-\varepsilon_2}{t}.
	\]
	That is $k(w_1, \varepsilon_2)<t$. On the other hand, $k(w_2,\varepsilon_2)=t$. It follows from the monotonicity of $s\mapsto k(s,\varepsilon)$ in Lemma \ref{lemma:monotonicity-k-s-epsilon} that $w_1<w_2$. 
	
	We claim that  $w(\varepsilon;\lambda,t)<2\varepsilon$ uniformly as $\varepsilon$ tends to zero for $\mathbb{C}\backslash B$ where $B$ is a closed ball with large radius. Assume that $||{x}-\lambda||>M\gg 1$ for $\lambda\in \mathbb{C}\backslash B$. Since $w=w(\varepsilon;\lambda,t)$ is the unique solution of 
	\[
	\int_0^\infty\frac{w}{w^2+u^2}d\mu_{|{x}-\lambda\unit|}(u)
	=\frac{w-\varepsilon}{t},
	\]
	it follows that
	\begin{equation}
	\label{eqn:lemma3.7-inequality-in-proof}
	\frac{w-\varepsilon}{t}<\frac{w}{w^2+M^2}<\frac{w}{M^2}. 
	\end{equation}
	One can then verity that $w(\varepsilon;\lambda,t)<2\varepsilon$ for $\lambda\in \mathbb{C}\backslash B$ provided that the radius of $B$ is large so that $M$ is sufficiently large. 
	
	By Lemma \ref{lemma:s-epsilon-t-limit}, the function $\lambda\mapsto w(\varepsilon;\lambda,t)$ is a continuous function of $\lambda$ for any $\varepsilon\geq 0$.
	Let $B$ be a closed ball in $\mathbb{C}$ with large radius so that $\Xi_t\subset B$. Because $w(\varepsilon;\lambda,t)$ converges pointwise to the continuous function $w(0;\lambda,t)$ as $\varepsilon$ tends to zero by Lemma \ref{lemma:s-epsilon-t-limit}, and $\{w(\varepsilon;\lambda,t)\}_{\varepsilon>0}$ is a monotone sequence of continuous functions of $\lambda$, it follows by Dini's theorem that $w(\varepsilon;\lambda,t)$ converge to the function $w(0;\lambda,t)$ uniformly on the closed ball $B$ as $\varepsilon$ tends to zero. 
	
	The above discussions show that $w(\varepsilon;\lambda,t)$ converge to the function $w(0;\lambda,t)$ uniformly on $\mathbb{C}$ as $\varepsilon$ tends to zero. 
\end{proof}

\subsection{The operator-valued subordination functions}
\label{subsection:operator-sub}

Let $y=g_{t,\gamma}$ be a twisted elliptic operator in $(\CM, \phi)$ and let ${x} \in \CM$ be a random variable that is $*$-free from $y$. We recall that
$$
G_{X}(\Theta(\lambda,{\i}\varepsilon))=\begin{bmatrix}
\displaystyle g_{X,11}(\lambda, \varepsilon) &
\displaystyle g_{X,12}(\lambda, \varepsilon)\\
\displaystyle g_{X,21}(\lambda, \varepsilon) &
\displaystyle  g_{X,22}(\lambda, \varepsilon)
\end{bmatrix},
$$
and
$$
G_{X+Y}(\Theta(\lambda,{\i}\varepsilon))=\begin{bmatrix}
\displaystyle g_{X+Y,11}(\lambda, \varepsilon) &
\displaystyle g_{X+Y,12}(\lambda, \varepsilon)\\
\displaystyle g_{X+Y,21}(\lambda, \varepsilon) &
\displaystyle  g_{X+Y,22}(\lambda, \varepsilon)
\end{bmatrix}.
$$

The main result in this section is the following. 
\begin{theorem}
	\label{thm:sub-X0-gt-formula}
	Let $y=g_{t,\gamma}\in \mathcal{M}$ and ${x}$ be a random variable that is free from $y$. We have 
	\begin{equation}
	\label{eqn:Omega-1-formula-in-proof-0}
	\Omega_1(\Theta(z,{\i}\varepsilon))
	=\Theta(z,{\i}\varepsilon)-R_Y(G_{X+Y}(\Theta(z,{\i}\varepsilon))).
	\end{equation}
	For any $\varepsilon>0$ and $z\in \mathbb{C}$, using notations in \eqref{eqn:X0-Ct} and \eqref{eqn:subord-X0-Ct}, 
	and set
	\begin{equation}
	\label{eqn-z-lambda-subordination}
	\lambda=z-\gamma\cdot \phi\bigg( (z-{x}-y)^*\big((z-{x}-y)(z-{x}-y)^*+\varepsilon^2 \big)^{-1} \bigg)
	\end{equation}
	Then $\Omega_1(\Theta(z,{\i}\varepsilon))=\Theta(\lambda, {\i}w(\varepsilon;\lambda,t))$. That is,
	\begin{align}
	\Omega_1\left( 
	\begin{bmatrix}
	{\i}\varepsilon & z\\
	\overline{z} & {\i}\varepsilon
	\end{bmatrix}\right)
	=\begin{bmatrix}
	{\i} w(\varepsilon;\lambda,t)  & \lambda\\
	\overline{\lambda} & {\i} w(\varepsilon; \lambda,t)
	\end{bmatrix},
	\end{align}
	where $w(\varepsilon;\lambda,t)$ is defined in Proposition \ref{prop:sub-identical-symmetric-measure} and Definition \ref{defn:s-epsilon-t-00}. 
	In addition, $w(\varepsilon;\lambda,t)$ can be expressed as 
	\[
	  w(\varepsilon;\lambda,t)=\varepsilon+t \varepsilon \phi\bigg( \big( (z-{x}-y)^*(z-{x}-y) +\varepsilon^2 \big)^{-1} \bigg)
	\]
	with $\lambda$ given by \eqref{eqn-z-lambda-subordination}. 
	
	The subordination relation $G_{X+Y}(\Theta(z,{\i}\varepsilon))=G_X(\Omega_1(\Theta(z,{\i}\varepsilon)))$ is expressed in terms of $z$ and $\varepsilon$ as
	\begin{align*}
		\mathbb{E}\left( \begin{bmatrix}
		{\i}\varepsilon & z-({x}+y)\\
		\overline{z}-({x}+y)^* & {\i}\varepsilon
		\end{bmatrix}^{-1} \right)
		=\mathbb{E}\left( \begin{bmatrix}
		  {\i}w(\varepsilon;\lambda,t) & \lambda-{x}\\
		  \overline{\lambda}-{x}^* & {\i}w(\varepsilon;\lambda,t)
		\end{bmatrix}^{-1} \right),
	\end{align*}
	which is also equivalent to
	\begin{align}
	\label{eqn:subordination-Cauchy-entries}
	g_{X+Y, 11}(z,\varepsilon)=g_{X, 11}(\lambda, w(\varepsilon;\lambda,t)),
	\qquad  g_{X+Y, 12}(z,\varepsilon)=g_{X, 12}(\lambda, w(\varepsilon;\lambda,t)).
	\end{align}
In other words, we have
\begin{align}
  	\label{eqn:subordination-Cauchy-entries-1}
  &\varepsilon \phi\bigg( \big( (z-{x}-y)^*(z-{x}-y) +\varepsilon^2 \big)^{-1} \bigg)\nonumber\\
   &\qquad\qquad
    =w(\varepsilon;\lambda,t) \phi\bigg( \big( (\lambda-{x})^*(\lambda-{x}) +w(\varepsilon;\lambda,t)^2 \big)^{-1} \bigg)
\end{align}
and
\begin{align}
  	\label{eqn:subordination-Cauchy-entries-2}
  &\ \phi\bigg((z-x-y)^* \big( (z-{x}-y)(z-{x}-y)^* +\varepsilon^2 \big)^{-1} \bigg)\nonumber\\
   &\qquad\qquad
    = \phi\bigg((\lambda-x)^* \big( (\lambda-{x})(\lambda-{x})^* +w(\varepsilon;\lambda,t)^2 \big)^{-1} \bigg).
\end{align}
\end{theorem}

\begin{proof}
	The only nonzero free cumulants of $\{y,y^*\}$ are
	\[
	\kappa(y,y^*)=\kappa(y^*,y)=t, \qquad \kappa(y, y)=\gamma, \qquad
	\kappa(y^*, y^*)=\overline{\gamma}.
	\]
	The operator $     Y=
	\begin{bmatrix}
	0 & y\\
	y^* & 0
	\end{bmatrix}$ is an operator-valued semicircular element whose $R$-transform is explicitly given by
	\[
	R_Y(b)=\mathbb{E}(YbY), \qquad b\in M_2(\mathbb{C}).
	\]
	Hence, for $b=\begin{bmatrix}
	a_{11} & a_{12}\\
	a_{21} & a_{22}
	\end{bmatrix}$,
	\[
	R_Y\left(\begin{bmatrix}
	a_{11} & a_{12}\\
	a_{21} & a_{22}
	\end{bmatrix}\right)
	=\begin{bmatrix}
	a_{22} \kappa(y, y^*) & a_{21} \kappa(y, y)\\
	a_{12}\kappa(y^*, y^*) & a_{11} \kappa(y^*, y)
	\end{bmatrix}
	=\begin{bmatrix}
	a_{22}t & a_{21}\gamma\\
	a_{12}\overline{\gamma} & a_{11}t
	\end{bmatrix}.
	\]
	In particular, for any $\varepsilon>0$ and $z\in\mathbb{C}$, we have 
	\begin{equation}
	 \label{eqn:R-Y-Theta}
	 R_Y(\Theta(z, {\i}\varepsilon))=
	R_Y\left(\begin{bmatrix}
	{\i}\varepsilon & z\\
	\overline{z} & {\i}\varepsilon
	\end{bmatrix}\right)=\begin{bmatrix}
	{\i}\varepsilon t & \overline{z} \gamma\\
	{z} \overline{\gamma}& {\i}\varepsilon t
	\end{bmatrix}.
	\end{equation}

	Since $X, Y$ are free with amalgamation in the operator-valued $W^*$-probability space $(M_2(\CM), \mathbb{E}, M_2(\mathbb{C}))$, we have 
	\[
	R_{X+Y}(b)=R_X(b)+R_Y(b).
	\]
	Hence
	\[
	G_{X+Y}^{\langle -1 \rangle}(b)=G_{X}^{\langle -1 \rangle}(b)
	+R_Y(b).
	\]
	By replacing $b$ with $G_{X+Y}(\beta)$, we obtain a formula for the subordination function
	\begin{equation}
	\Omega_1(\beta)=  G_{X}^{\langle -1 \rangle}(G_{X+Y}(\beta))
	= \beta-R_Y(G_{X+Y}(\beta))
	\end{equation}
	for $\beta$ in a neighborhood of infinity. Hence, for and $z\in\mathbb{C}$ and $\varepsilon$ large, the identity \eqref{eqn:Omega-1-formula-in-proof-0} holds: 
	\begin{equation}\nonumber
	\Omega_1(\Theta(z,{\i}\varepsilon))
	=\Theta(z,{\i}\varepsilon)-R_Y(G_{X+Y}(\Theta(z,{\i}\varepsilon))).
	\end{equation}
	
	We next show that \eqref{eqn:Omega-1-formula-in-proof-0} holds for any $z\in\mathbb{C}$ and $\varepsilon>0$. Observing that the right hand side of \eqref{eqn:Omega-1-formula-in-proof-0}
	is defined for any $\varepsilon>0$, we only need to check that the right hand side of \eqref{eqn:Omega-1-formula-in-proof-0} is the correct expression of $\Omega_1(\Theta(z,{\i}\varepsilon))$. Recall that 
	\[
	G_{X+Y}(\Theta(z,{\i}\varepsilon))=
	\begin{bmatrix}
	g_{X+Y, 11}(z,\varepsilon) & g_{X+Y, 12}(z,\varepsilon)\\
	g_{X+Y, 21}(z,\varepsilon) & g_{X+Y, 22}(z,\varepsilon)
	\end{bmatrix}
	\]
	holds for any $z\in\mathbb{C}$ and $\varepsilon>0$. Hence,
	\[
	R_Y(G_{X+Y}(\Theta(z,{\i}\varepsilon)))
	=\begin{bmatrix}
	t\cdot g_{X+Y, 22}(z,\varepsilon) & \gamma \cdot g_{X+Y, 21}(z,\varepsilon)\\
	\overline{\gamma}\cdot g_{X+Y, 12}(z,\varepsilon)
	& t\cdot g_{X+Y, 11}(z,\varepsilon)
	\end{bmatrix},
	\]
	where we reminder the reader the notations of $g_{X+Y, ij}$ in Section \ref{subsection:2.3Hermitian} and 
	\[
	g_{X+Y, 11}(z,\varepsilon)=g_{X+Y, 22}(z,\varepsilon),
	\qquad 
	g_{X+Y, 12}(z,\varepsilon)=\overline{g_{X+Y, 21}(z,\varepsilon)}.
	\]
	Therefore, for $\varepsilon>0$ large enough, by \eqref{eqn:Omega-1-formula-in-proof-0}, we have 
	\begin{align*}
	\Omega_1(\Theta(z,{\i}\varepsilon))
	&=
	\begin{bmatrix}
	{\i}\varepsilon - t\cdot g_{X+Y, 22}(z,\varepsilon)
	&
	z-\gamma \cdot g_{X+Y, 21}(z,\varepsilon)\\
	\overline{z}-\overline{\gamma}\cdot g_{X+Y, 12}(z,\varepsilon)
	&
	{\i}\varepsilon-t\cdot g_{X+Y, 11}(z,\varepsilon)
	\end{bmatrix}.
	\end{align*}
	In particular, for $\varepsilon>0$ large enough, we have 
	\begin{align*}
	 &\text{the imaginary part of (1,1)-entry of }\, \Omega_1(\Theta(z,{\i}\varepsilon))\\
	 &=\varepsilon +{\i} t\cdot g_{X+Y, 22}(z,\varepsilon)\\
	 &=\varepsilon + t\varepsilon\phi\bigg( \big( (z-{x}-y)^*(z-{x}-y)+\varepsilon^2 \big)^{-1} \bigg)>\varepsilon.
	\end{align*}
	For any $\delta\in\mathbb{C}^+$, $\Im G_{X+Y}(\Theta(z,\delta))<0$ and hence
	$\Im R_Y(G_{X+Y}(\Theta(z,\delta)))<0$. Hence,
	$\Im \Omega_1(\Theta(z,\delta))>\Im \delta>0$ 
	for any $\delta\in\mathbb{C}^+$. The map $\delta\mapsto R_Y(G_{X+Y}(\Theta(z,\delta)))$
	is a complex analytic function of $\delta$ in $\mathbb{C}^+$.
	Since $\delta\mapsto \Omega_1(\Theta(z,\delta))$ is also a complex analytic function of $\delta$ in $\mathbb{C}^+$ and $\Im\Omega_1(\Theta(z,\delta))>\Im\delta$, we conclude that for any $\varepsilon>0$ we have 
	\[
	 \Omega_1(\Theta(z,{\i}\varepsilon))
	=\Theta(z,{\i}\varepsilon)-R_Y(G_{X+Y}(\Theta(z,{\i}\varepsilon)))
	\]
	 by the uniqueness of analytic functions.

	For any $z\in\mathbb{C}$ and $\varepsilon>0$, we denote 
	\begin{align}
	\label{eqn-defn-varepsilon0-z-2}
	\lambda&=z-\gamma \cdot g_{X+Y, 21}(z,\varepsilon)\\
	 &=z-\gamma\cdot \phi\bigg( (z-{x}-y)^*\big((z-{x}-y)(z-{x}-y)^*+\varepsilon^2 \big)^{-1} \bigg)\nonumber
	\end{align}
	and 
	\begin{align}
	\label{eqn-defn-varepsilon0-z-1}
	\varepsilon_0&= \varepsilon+{\i} t\cdot g_{X+Y, 22}(z,\varepsilon)\\
	 &=\varepsilon+t \varepsilon \phi\bigg( \big( (z-{x}-y)^*(z-{x}-y) +\varepsilon^2 \big)^{-1} \bigg).\nonumber
	\end{align}
	Then, for any $z\in\mathbb{C}$ and $\varepsilon>0$, we have 
	\begin{equation}
	\label{eqn:Omega-1-formula-in-proof}
	\Omega_1(\Theta({z},\varepsilon))=\Theta(\lambda,\varepsilon_0)= \begin{bmatrix}
	{\i}\varepsilon_0  &   \lambda\\
	\overline{\lambda} & {\i}\varepsilon_0
	\end{bmatrix}
	\end{equation}
	Hence, the Cauchy transform of $X$ at $\Omega_1(\Theta({z},\varepsilon))$ is given by
	\begin{equation}
	\nonumber
	G_X(\Omega_1(\Theta(z,\varepsilon)))= \mathbb{E}[ (\Omega_1(\Theta(z,\varepsilon))-X)^{-1}]
	=     \begin{bmatrix}
	g_{X,11}(\lambda,\varepsilon_0) &  g_{X,12}(\lambda,\varepsilon_0)\\
	g_{X,21}(\lambda,\varepsilon_0) & g_{X,22}(\lambda,\varepsilon_0)
	\end{bmatrix}.
	\end{equation}
	On the other hand, we have 
	\[
	G_{X+Y}(\Theta(z,\varepsilon))=
	\begin{bmatrix}
	g_{X+Y, 11}(z,\varepsilon) & g_{X+Y, 12}(z,\varepsilon)\\
	g_{X+Y, 21}(z,\varepsilon) & g_{X+Y, 22}(z,\varepsilon)
	\end{bmatrix}.
	\]
	Therefore the subordination relation $ G_{X+Y}(\Theta(z,\varepsilon))=G_X(\Omega_1(\Theta(z,\varepsilon)))$ is equivalent to
	\begin{equation}
	\label{eqn-in-proof-3.19-entries}
	g_{X+Y, 11}(z,\varepsilon)= g_{X,11}(\lambda,\varepsilon_0),
	\qquad
	g_{X+Y, 12}(z,\varepsilon)= g_{X,12}(\lambda,\varepsilon_0)
	\end{equation}
	for any $z\in\mathbb{C}$, $\varepsilon>0$, where 
	where $\lambda, \varepsilon_0$ are given by the relations \eqref{eqn-defn-varepsilon0-z-1} and \eqref{eqn-defn-varepsilon0-z-2}.
	
  It remains to show that $\varepsilon_0=w(\varepsilon;\lambda,t)$. We reminder the reader that 
	\[
	g_{X, 11}(\lambda,\varepsilon_0)=g_{X, 22}(\lambda,\varepsilon_0),
	\qquad
	g_{X, 12}(\lambda,\varepsilon_0)=\overline{g_{X, 21}(\lambda,\varepsilon_0)}.
	\]
	Now the relation \eqref{eqn-defn-varepsilon0-z-1} can be rewritten as
	\[
	\varepsilon_0=\varepsilon+{\i}t g_{X,22}(\lambda,\varepsilon_0)
	\]
	which is equivalent to
	\[
	\varepsilon_0=\varepsilon+t\cdot\varepsilon_0\phi\left(\big( (\lambda-{x})^*(\lambda-{x})+\varepsilon_0^2 
	\big)^{-1}\right).
	\]
	This can be further rewritten as
	\[
	\frac{\varepsilon_0-\varepsilon}{\varepsilon_0} \left(\int_0^\infty \frac{1}{\varepsilon_0^2+u^2}d\mu_{|{x}-\lambda\unit|}(u)\right)^{-1} =t.
	\]
	Hence, $\varepsilon_0=w(\varepsilon;\lambda,t)$ is the unique solution to the above equation following Definition \ref{defn:s-epsilon-t-00}. This finishes the proof. 
\end{proof}

For any $z\in\mathbb{C}$ and $\varepsilon>0$, if $\lambda$ is given by \eqref{eqn-z-lambda-subordination} which is rewritten as
\[
  \lambda=z-\gamma\cdot g_{X+Y, 21}(z,\varepsilon).
\] 
Then by the subordination relate \eqref{eqn:subordination-Cauchy-entries} and the proof of Theorem \ref{thm:sub-X0-gt-formula}, we have 
\[
   z=\lambda+\gamma\cdot g_{X+Y, 21}(z,\varepsilon)=
    \lambda+\gamma\cdot g_{X,21}(\lambda, w(\varepsilon;\lambda,t)).
\]
In Section \ref{section:regularized-pushforward-map}, we show that the relation $\lambda\mapsto z=\lambda+\gamma \cdot g_{X, 21}(\lambda, w(\varepsilon;\lambda,t))$ determines a homeomorphism of the complex plane. 

For any $d\in \mathbb{H}^+(M_2(\mathbb{C}))$, we denote 
\[
H(d)=d+R_Y(G_X(d)).
\]
Then by \eqref{eqn:Omega-1-formula-in-proof-0} for any $z\in\mathbb{C}$ and $\varepsilon>0$, the identity $H(\Omega_1(\Theta(z,\i \varepsilon)))=\Theta(z,\i \varepsilon)$
 holds tautologically. The following is a special case of \cite[Lemma 4.2]{BelinschiCapitaine2017jfa}. 
\begin{lemma}
\label{lemma:right-inverse-of-Omega1}
For any $d\in \mathbb{H}^+(M_2(\mathbb{C}))$ such that $\Im H(d)>0$ we have 
\begin{equation}
   \label{eqn:right-inverse-of-Omega1}
     \Omega_1(H(d))=d.
\end{equation}
\end{lemma}

\begin{corollary}
\label{cor:Omega-1-formula}
Given $\varepsilon>0$ and $\lambda\in\mathbb{C}$, let $z=\lambda+\gamma\cdot g_{X,21}(\lambda, w(\varepsilon;\lambda,t))$, then $\Omega_1(\Theta(z,{\i}\varepsilon))=\Theta(\lambda,{\i}w(\varepsilon;\lambda,t))$.
\end{corollary}
\begin{proof}
Denote $d=\Theta(\lambda,{\i}w(\varepsilon;\lambda,t))=\begin{bmatrix}
  {\i}w(\varepsilon;\lambda,t) & \lambda\\
  \overline{\lambda} & {\i}w(\varepsilon;\lambda,t)
\end{bmatrix}$. 
 By Lemma \ref{lemma:right-inverse-of-Omega1}, it suffices to show that
 $H(d)=\Theta(z,{\i}\varepsilon)=\begin{bmatrix}
   {\i}\varepsilon & z\\
   \overline{z} & {\i}\varepsilon
 \end{bmatrix}$. 
 Denote $w=w(\varepsilon;\lambda,t)$. Recalling the identity \eqref{eqn:R-Y-Theta}, 
 by the definition of $H$ and \eqref{eqn:condition-s-general-t-v2} in Definition \ref{defn:s-epsilon-t-00}, we have 
 \begin{align*}
  &\text{the imaginary part of (1,1)-entry of }\, H(d)\\
  &=w+t\cdot g_{X,22}(\lambda,w)\\
  &=w-tw\phi\big( ( (\lambda-{x})^*(\lambda-{x})+w^2)^{-1} \big)=\varepsilon.
 \end{align*}
The $(1,2)$-entry of $H(d)$ is exactly equal to $\lambda+\gamma\cdot g_{X,21}(\lambda, w)$, which is $z$.
\end{proof}

\subsection{The coupling Fuglede-Kadision determinants}

To help us remember entries of the Cauchy transform as in \eqref{eqn-Cauchy-X-Theta} are derivatives, we introduce the following notation. 
\begin{notation}
	\label{notation:partial-derivatives-lambda-c-x}
	For any $t>0$, let $c_t$ be a circular operator with variance $t$. For $\gamma\in\mathbb{C}$ such that $|\gamma|\leq t$, denote by $y=g_{t,\gamma}$ the twisted elliptic operator with parameters $t,\gamma$. For any $x\in \CM$, $\lambda\in\mathbb{C}$ and $\varepsilon\geq 0$, we set
	\begin{align*}
	S(x,\lambda, \varepsilon) 
	&=\log\Delta(  (x-\lambda\unit)^*(x-\lambda\unit)+\varepsilon^2\unit)\\
	&=\log\Delta(  (x-\lambda\unit)(x-\lambda\unit)^*+\varepsilon^2\unit).
	\end{align*}
	It is convenient to introduce the following notations
	\begin{align*}
	p_\lambda^{c, (t)}(\varepsilon) &=\frac{\partial S}{\partial \lambda}({x}+c_t, \lambda, \varepsilon)\\
	&=-\phi\left[({x}+c_t-\lambda\unit)^*\big( ({x}+c_t-\lambda\unit)({x}+c_t-\lambda\unit)^*+\varepsilon^2\unit\big)^{-1}\right]\nonumber\\
	p_{z}^{g, (t,\gamma)}( \varepsilon) &=\frac{\partial S}{\partial z}({x}+y, z, \varepsilon)\\
	&=-\phi\left[({x}+y-z\unit)^*\big( ({x}+y-z\unit)({x}+y-z\unit)^*+
	\varepsilon^2\unit\big)^{-1}\right] \nonumber\\
	p_\lambda^{(0)}(\varepsilon)&=\frac{\partial S}{\partial \lambda}({x}, \lambda, \varepsilon)=-\phi\left[({x}-\lambda\unit)^*\big( ({x}-\lambda\unit)({x}-\lambda\unit)^*+\varepsilon^2\unit\big)^{-1}\right]\nonumber\\
	q_\varepsilon^{c, (t)}(\lambda) &=\frac{1}{2}\frac{\partial S}{\partial \varepsilon}({x}+c_t, \lambda, \varepsilon)=\varepsilon\phi\left[\big( ({x}+c_t-\lambda\unit)^*({x}+c_t-\lambda\unit)+\varepsilon^2\unit \big)^{-1}\right]\nonumber\\
	q_\varepsilon^{g, (t,\gamma)}(\lambda) &=\frac{1}{2}\frac{\partial S}{\partial \varepsilon}({x}+y, \lambda, \varepsilon)=\varepsilon\phi\left[\big( ({x}+y-\lambda\unit)^*({x}+y-\lambda\unit)+\varepsilon^2\unit \big)^{-1}\right]\nonumber\\
	q_\varepsilon^{(0)}(\lambda) &=\frac{1}{2}\frac{\partial S}{\partial \varepsilon}({x}, \lambda, \varepsilon)=\varepsilon\phi\left[\big( ({x}-\lambda\unit)^*({x}-\lambda\unit)+\varepsilon^2\unit \big)^{-1}\right]\nonumber
	\end{align*}
	We also set
	\[
	p_{\overline{\lambda}}^{c, (t)}(\varepsilon)=\frac{\partial S}{\partial \overline{\lambda}}({x}+c_t, \lambda, \varepsilon),
	\quad
	p_{\overline{z}}^{g, (t,\gamma)}(\varepsilon)=\frac{\partial S}{\partial \overline{z}}({x}+y, z, \varepsilon),
	\quad
	p_{\overline{\lambda}}^{(0)}(\varepsilon)=\frac{\partial S}{\partial \overline{\lambda}}({x}, \lambda, \varepsilon).
	\] 
\end{notation}

We note that 
\[
p_{\overline{\lambda}}^{c, (t)}(\varepsilon)=\overline{ p_{{\lambda}}^{c, (t)}(\varepsilon)},
\quad
p_{\overline{z}}^{g, (t,\gamma)}(\varepsilon)
=\overline{p_{{z}}^{g, (t,\gamma)}(\varepsilon)},
\quad
p_{\overline{\lambda}}^{(0)}(\varepsilon)=\overline{p_{{\lambda}}^{(0)}(\varepsilon)}.
\]
and the Cauchy transform can be written as 
\begin{equation}
\label{eqn:sub-X0-Ct-p-lambdas}
G_{X+Y}(\Theta(\lambda,\varepsilon))=\mathbb{E} \begin{bmatrix}
(\Theta(\lambda,\varepsilon)-X-Y)^{-1}
\end{bmatrix}
=
\begin{bmatrix}
-{\i} q_\varepsilon^{(t,\gamma)}(\lambda) &  p_{\overline{\lambda}}^{(t,\gamma)}(\varepsilon)\\
p_{\lambda}^{(t,\gamma)}(\varepsilon) & -{\i} q_\varepsilon^{(t,\gamma)}(\lambda)
\end{bmatrix}.
\end{equation}
\begin{corollary}
	\label{thm:sub-comparison-Ct}
	The subordination relation \eqref{eqn:subordination-Cauchy-entries} is equivalent to
	\begin{align}
	\label{eqn:p-q-identity-x-t-gamma}
	q_{w(\varepsilon)}^{(0)}(\lambda)=q_{\varepsilon}^{g, (t,\gamma)}({z})=q_\varepsilon^{c, (t)}(\lambda), 
	\quad p_{\lambda}^{(0)}(w(\varepsilon))=p_{{z}}^{g,(t,\gamma)}(\varepsilon)=p_\lambda^{c,(t)}(\varepsilon).
	\end{align}
	where $w(\varepsilon)=w(\varepsilon;\lambda,t)$. In particular, we have
	\begin{align*}
  &\varepsilon \phi\bigg( \big( (\lambda-{x}-c_t)^*(\lambda-{x}-c_t) +\varepsilon^2 \big)^{-1} \bigg)\\
   &\qquad\qquad
    =w(\varepsilon;\lambda,t) \phi\bigg( \big( (\lambda-{x})^*(\lambda-{x}) +w(\varepsilon;\lambda,t)^2 \big)^{-1} \bigg)
\end{align*}
and
\begin{align*}
  &\ \phi\bigg((\lambda-x-c_t)^* \big( (\lambda-{x}-c_t)(\lambda-{x}-c_t)^* +\varepsilon^2 \big)^{-1} \bigg)\\
   &\qquad\qquad
    = \phi\bigg((\lambda-x)^* \big( (\lambda-{x})(\lambda-{x})^* +w(\varepsilon;\lambda,t)^2 \big)^{-1} \bigg).
\end{align*}
\end{corollary}
\begin{proof}
	We note that $c_t=g_{t,0}$. Then \eqref{eqn-z-lambda-subordination} reads $z=\lambda$ if $\gamma=0$. The result follows from \eqref{eqn:subordination-Cauchy-entries-1} and \eqref{eqn:subordination-Cauchy-entries-2}.
\end{proof}

The following proof was inspired by the proof of \cite[Lemma 4.14]{HaagerupSchultz2007}.
\begin{lemma}
	\label{lemma:main-FK-det-gt}
	Let $y=g_{t,\gamma}$ and ${x}$ be a random variable free from $y$. 
	For any $\lambda\in\mathbb{C}$ and $(\varepsilon, t)\in (0,\infty)\times(0,\infty)$, we have the coupling Fuglede-Kadison determinant formula
	\begin{align}
	\label{eqn:main-FK-det-V2-gt}
	\Delta \big( ({x}+y-z\unit)^* &({x}+y-z\unit))+\varepsilon^2\unit \big)
	={\Delta\big( ({x}-\lambda\unit)^*({x}-\lambda\unit)+w(\varepsilon)^2\unit \big)}\nonumber \\
	&\qquad\qquad \times \exp\left [\Re \bigg( \gamma\cdot (  p_\lambda^{(0)}(w(\varepsilon)) ) ^2\bigg)-\frac{(w(\varepsilon)-\varepsilon)^2}{t}   \right]
	\end{align}
	where $ z=\lambda+\gamma \cdot p_\lambda^{(0)}(w(\varepsilon)) $ and
	$w(\varepsilon)= w(\varepsilon;\lambda,t)$ is defined in Definition \ref{defn:s-epsilon-t-00}.
\end{lemma}
\begin{proof}
	Fix $\lambda\in\mathbb{C}$, then  $w(\varepsilon)= w(\varepsilon;\lambda, t)$ and $p_\lambda^{(0)}(w(\varepsilon)) $ are then completely determined by $\varepsilon$. 
	Recall tha $p_\lambda^{(0)}(w(\varepsilon))=p_z^{(t,\gamma)}(\varepsilon)$, we have 
	\[
	z=\lambda+\gamma\cdot p_\lambda^{(0)}(w(\varepsilon))=\lambda+p_z^{(t,\gamma)}(\varepsilon).
	\]
	Therefore, 
	\begin{align*}
	&\frac{d}{d\varepsilon} S({x}+y, z, \varepsilon)\\
	=&2q_\varepsilon^{(t,\gamma)}(z)+\left( \gamma \cdot p_{{z}}^{(t,\gamma)}(\varepsilon) \frac{d}{d\varepsilon}p_{{z}}^{(t,\gamma)}(\varepsilon)+\overline{\gamma}\cdot p_{\overline{{z}}}^{(t,\gamma)}(\varepsilon)
	\frac{d}{d\varepsilon}p_{\overline{{z}}}^{(t,\gamma)}(\varepsilon) \right)\\
	=& 2q_{w(\varepsilon)}^{(0)}(\lambda) + \bigg(\gamma\cdot p_\lambda^{(0)}(w(\varepsilon)) \frac{d}{d\varepsilon}p_\lambda^{(0)}(w(\varepsilon)) 
	+\overline{\gamma}\cdot p_{\overline{\lambda}}^{(0)}(w(\varepsilon)) \frac{d}{d\varepsilon}p_{\overline{\lambda}}^{(0)}(w(\varepsilon))  \bigg). 
	\end{align*}
   	Definition \ref{defn:s-epsilon-t-00} says that $w(\varepsilon)$ satisfies 
	\[
	\frac{w(\varepsilon)-\varepsilon}{t}= q_{w(\varepsilon)}^{(0)}(\lambda)=w(\varepsilon)\cdot\phi\bigg[ \big((\lambda-x)^*(\lambda-x)+w(\varepsilon)^2 \big)^{-1} \bigg].
	\]
	Recall that $\frac{1}{2}\frac{d}{d\varepsilon}S(x,\lambda,\varepsilon)=q_\varepsilon^{(0)}(\lambda)$. 
	We  then have 
	\begin{align*}
	& \int_{\varepsilon_0}^\varepsilon 
	q_{w(u)}^{(0)}(\lambda) du\\
	&=\int_{\varepsilon_0}^\varepsilon  q_{w(u)}^{(0)}(\lambda) \frac{d}{du} w(u) du+\int_{\varepsilon_0}^\varepsilon q_{w(u)}^{(0)}(\lambda) 
	\frac{d}{du}(u-w(u)) du\\
	&=\frac{S({x},\lambda,w(\varepsilon))-S({x},\lambda,w(\varepsilon_0))}{2}  -t\cdot\int_{\varepsilon_0}^\varepsilon
	q_{w(u)}^{(0)}(\lambda) \left(\frac{d}{du}q_{w(u)}^{(0)}(\lambda) \right)  \\
	&= \frac{S({x},\lambda,w(\varepsilon))-S({x},\lambda,w(\varepsilon_0))}{2}    -\frac{t}{2}\bigg( (q_{w(\varepsilon)}^{(0)}(\lambda) )^2 -( q_{w(\varepsilon_0)}^{(0)}(\lambda) )^2\bigg)
	\end{align*}    
	and
	\begin{align*}
	&\gamma \int_{\varepsilon_0}^\varepsilon 
	p_\lambda^{(0)}(w(u))\left( \frac{d}{du}p_\lambda^{(0)}(w(u)) \right)du
	+\overline{\gamma} \int_{\varepsilon_0}^\varepsilon     p_{\overline{\lambda}}^{(0)} (w(u))\left(\frac{d}{du}p_{\overline{\lambda}}^{(0)}(w(u)) \right) 
	\\
	=& \frac{\gamma}{2}\bigg(( p_\lambda^{(0)}(w(\varepsilon)))^2-(p_\lambda^{(0)}(w(\varepsilon_0))^2\bigg)
	+\frac{\overline{\gamma}}{2}\bigg((p_{\overline{\lambda}}^{(0)}(w(\varepsilon))^2-(p_{\overline{\lambda}}^{(0)}(w(\varepsilon_0))^2\bigg)
	\end{align*}
	Hence, there exists a constant $C$ such that 
	\begin{align*}
	&S(x+y,z,\varepsilon)-S(x,\lambda,w(\varepsilon))\\
	&=\frac{1}{2} \bigg( \gamma\cdot ( p_\lambda^{(0)}(w(\varepsilon)))^2+\overline{\gamma}\cdot (p_{\overline{\lambda}}^{(0)}(w(\varepsilon))^2-2t\cdot (q_{w(\varepsilon)}^{(0)}(\lambda) )^2  
	\bigg)+C\\
	&=\Re \bigg( \gamma\cdot (  p_\lambda^{(0)}(w(\varepsilon)) ) ^2\bigg)-t\cdot (q_{w(\varepsilon)}^{(0)}(\lambda) )^2  +C\\
	&=\Re \bigg( \gamma\cdot (  p_\lambda^{(0)}(w(\varepsilon)) ) ^2\bigg)-\frac{(w(\varepsilon)-\varepsilon)^2}{t}+C.
	\end{align*}
	
	It remains to show that $C=0$. We have
	\[
	S(x,\lambda,\varepsilon)=\int_0^\infty \log (u^2+\varepsilon^2) d\mu_{|{x}-\lambda\unit|}(u).
	\]
	Hence
	\begin{equation}
	\label{eqn:4.17-in-proof-limit}
	\lim_{\varepsilon \rightarrow \infty} ( S(x,\lambda,\varepsilon)-2\log \varepsilon)
	=\lim_{\varepsilon\rightarrow \infty} \log\left(1+\frac{u^2}{\varepsilon^2} \right)d\mu_{|{x}-\lambda\unit|}(u)=0.
	\end{equation}
	Observe that $\lim_{\varepsilon\rightarrow \infty}(w(\varepsilon)-\varepsilon)=0$. 
	Consequently,
	$\lim_{\varepsilon\rightarrow \infty} p_\lambda^{(0)}(w(\varepsilon))=0$.
	Hence
	\[
	\lim_{\varepsilon\rightarrow\infty}{z}=\lim_{\varepsilon\rightarrow\infty}(\lambda+t p_\lambda^{(0)}(w(\varepsilon))=\lambda
	\]
	and by a similar estimation as \eqref{eqn:4.17-in-proof-limit}
	\[
	\lim_{\varepsilon\rightarrow\infty} (S(x+y,z,\varepsilon)-2\log \varepsilon)=0
	\]
	Moreover, 
	\[
	\lim_{\varepsilon\rightarrow\infty} q_{w(\varepsilon)}^{(0)}(\lambda)=0.
	\]
	We conclude that $C$ must be zero. 
\end{proof}

\begin{theorem}
	\label{thm:main-FK-det-ct-0}
	For $\lambda\in\mathbb{C}$, set $\mu=\mu_{|{x}-\lambda\unit|}$ and let $w(0;\lambda,t)$ be as in  Definition \ref{defn:s-epsilon-t-00}.
	
	\begin{enumerate}
		\item If $\lambda\in\Xi_t$, then 
		\begin{equation}
		\label{eqn:main-FK-det-0-V2}
		\begin{aligned}
		\Delta \big({x}+c_t-\lambda\unit \big)^2
		&={\Delta\big( ({x}-\lambda\unit)^*({x}-\lambda\unit)+w(0; \lambda,t)^2 \unit \big)}\\ &\qquad\qquad\times{\exp\bigg(-\frac{(w(0;\lambda,t))^2}{t}\bigg)}.
		\end{aligned}
		\end{equation}

		\item If $\lambda\notin\Xi_t$, then 
		\[
		\Delta({x}+c_t-\lambda\unit)=\Delta({x}-\lambda\unit).
		\]
	\end{enumerate}
\end{theorem}
\begin{proof}
	The circular operator $c_t$ corresponds to $\gamma=0$ for a twisted elliptic operator. 
	If $\lambda\in\Xi_t$,
	then $w(0;\lambda,t)>0$. 
	Hence the first part follows from \eqref{eqn:main-FK-det-V2-gt} by letting $\varepsilon$ tend to zero. 
	For the second part, 
	note that 
	\[
	\Delta({x}+c_t-\lambda\unit)^2=\lim_{\varepsilon\rightarrow 0^+}  \Delta \big( ({x}+c_t-\lambda\unit)^*({x}+c_t-\lambda\unit)+\varepsilon^2\unit \big).
	\]
	By Lemma \ref{lemma:s-epsilon-t-limit}, we know if $\lambda\notin\Xi_t$, then $\lim_{\varepsilon\rightarrow 0^+} w(\varepsilon; \lambda,t)=0$. Hence by \eqref{eqn:main-FK-det-V2-gt}, we have
	\begin{align*}
	\Delta({x}+c_t-\lambda\unit)^2&=\lim_{\varepsilon\rightarrow 0^+}  \Delta \big( ({x}+c_t-\lambda\unit)^*({x}+c_t-\lambda\unit)+\varepsilon^2\unit \big)\\
	&=\lim_{\varepsilon\rightarrow 0^+} 
	{\Delta\big( ({x}-\lambda\unit)^*({x}-\lambda\unit)+w(\varepsilon)^2\unit \big)}\nonumber 
	\times \exp\left [-\frac{(w(\varepsilon)-\varepsilon)^2}{t}   \right]\\
	&=\Delta({x}-\lambda\unit)^2.
	\end{align*}
	This finishes the proof. 
\end{proof}

In light of Lemma \ref{lemma:main-FK-det-gt}, we would like to let $\varepsilon$ tend to zero for both sides of \eqref{eqn:main-FK-det-V2-gt}. Note that the the left hand side of \eqref{eqn:main-FK-det-V2-gt} is $\Delta \big( ({x}+y-z\unit)^* ({x}+y-z\unit))+\varepsilon^2\unit \big)$ where $z=\lambda+\gamma\cdot p_\lambda^{(0)}(w(\varepsilon;\lambda,t))$ also depends on $\varepsilon$. Hence, we need some regularity results to allow us to take the limit as we wish. To this end, we introduce the map $\Psi_{c,0}$ given by 
\[
\Psi_{c,0} (\lambda,\varepsilon)=(\lambda,w(\varepsilon;\lambda,t))
\]
where $w(\varepsilon;\lambda,t)$ was given in Definition \ref{defn:s-epsilon-t-00},
and the map $\Psi_{c,g}$ by 
\[
\Psi_{c,g}(\lambda,\varepsilon)=(z,\varepsilon)
\]
where 
\[
z=\lambda+\gamma\cdot p_\lambda^{(0)}(w(\varepsilon;\lambda,t)). 
\]

\begin{lemma}
	\label{lemma:Jacobian-circular-0}
	If $\lambda\in \Xi_t$, the Jacobian of $\Psi_{c,0}$ at $(\lambda,0)$ is invertible. 
\end{lemma}
\begin{proof}
	To show the Jacobian of $\Psi_{c,0}$ is invertible, it suffices to show that $\frac{\partial w(\varepsilon;\lambda,t)}{\partial \varepsilon}\neq 0$ at $\varepsilon=0$. Recall that $w(\varepsilon;\lambda,t)$ is the unique solution $s>0$ for 
	\[
	\int_0^\infty\frac{s}{s^2+u^2}d\mu_{|{x}-\lambda\unit|}(u)
	=\frac{s-\varepsilon}{t}.
	\]
	When $\varepsilon=0$, note that $w(0;\lambda,t)>0$ and we can rewrite \eqref{eqn:w-epsiton-0-identity} as
	\[
	\int_0^\infty\frac{1}{w(0;\lambda,t)^2+u^2}d\mu_{|{x}-\lambda\unit|}(u)=\frac{1}{t}.
	\]
	A direct calculation shows that $\frac{\partial w(\varepsilon;\lambda,t)}{\partial \varepsilon}\neq 0$ at $\varepsilon=0$
\end{proof}

\begin{lemma}
\label{lemma:convergence-p-lambda}
The function $p_\lambda^{(0)}(w(\varepsilon;\lambda,t))$ converges uniformly to $p_\lambda^{(0)}(w(0;\lambda,t))$ on $\mathbb{C}$ as $\varepsilon$ tends to zero.
\end{lemma}
\begin{proof}
For $\varepsilon_2>\varepsilon_1>0$, we set $w_i=w(\varepsilon_i;\lambda,t) (i=1,2)$. Recall that $w(\varepsilon;\lambda,t)>\varepsilon$. By the proof of Lemma \ref{lemma:sub-w-uniform-convergence}, we have $w_1<w_2$. 
By the resolvent identity, we have 
\begin{align*}
  & p_\lambda^{(0)}(w(\varepsilon_2;\lambda,t))-p_\lambda^{(0)}(w(\varepsilon_1;\lambda,t))\\
   =& \phi\bigg[ (\lambda\unit-{x})^*\big( (\lambda-{x})(\lambda-{x})^*+w_2^2 \unit\big)^{-1}\bigg]\\
   &\qquad - \phi\bigg[ (\lambda\unit-{x})^*\big((\lambda-{x})(\lambda-{x})^*+w_1^2 \unit\big)^{-1}\bigg]\\
   =&\phi(AH_1 H_2 B) (w_2-w_1),
\end{align*}
where
\begin{align*}
 A&=(\lambda\unit-{x})^*\big( (\lambda-{x})(\lambda-{x})^*+w_1^2 \unit\big)^{-1/2}\\
 B&=\big( (\lambda-{x})(\lambda-{x})^*+w_2^2 \unit\big)^{-1/2}
   (w_1+w_2)\\
  H_1&=\big( (\lambda-{x})(\lambda-{x})^*+w_1^2 \unit\big)^{-1/2}\\
 H_2&=\big( (\lambda-{x})(\lambda-{x})^*+w_2^2 \unit\big)^{-1/2}.
\end{align*}
Since $\varepsilon_1< w_1<w_2$, we have 
\[
   ||A||\leq 1, \qquad ||B||\leq 2.
  \]
We also have, for $i=1,2$, 
\begin{align*}
	 \phi(H_i^2)&=\phi( \big( (\lambda-{x})(\lambda-{x})^*+w_i^2 \unit\big)^{-1}  )\\
	  &\leq \phi( \big( (\lambda-{x})(\lambda-{x})^*+w(0; \lambda,t)^2 \unit\big)^{-1}  )  \leq \frac{1}{t},
\end{align*}
due to $w(0;\lambda,t)\leq w_i$. 
Therefore,
\begin{align*}
 |\phi(AH_1H_2B)|\leq 2\phi(H_1H_2)
   \leq 2\sqrt{ \phi(H_1^2) \phi(H_2^2) }\leq \frac{2}{t},
\end{align*}
which yields
\[
  |p_\lambda^{(0)}(w(\varepsilon_2;\lambda,t))-p_\lambda^{(0)}(w(\varepsilon_1;\lambda,t))|
   \leq \frac{2}{t}(w(\varepsilon_2;\lambda,t)-w(\varepsilon_1;\lambda,t)).
\]
We note that for $\lambda\notin \Xi_t$, then $w(0;\lambda,t)=0$ and in this case
the operator $\big( (\lambda-{x})(\lambda-{x})^*\big)^{-1}$ is regarded as an unbounded operator affiliated with $\mathcal{M}$. By Cauchy-Schwarz inequality, we have 
\begin{align*}
   &|p_\lambda^{(0)}(0)|^2=\left|\phi\bigg[ (\lambda\unit-{x})^*\big( (\lambda-{x})(\lambda-{x})^*\big)^{-1}\bigg]\right|^2\\
    &\leq \phi\bigg[\left|(\lambda\unit-{x})^*\big( (\lambda-{x})(\lambda-{x})^*\big)^{-1/2}\right|^2\bigg]
    \cdot \phi\bigg[ \big( (\lambda-{x})(\lambda-{x})^*\big)^{-1} \bigg]
    \leq \frac{1}{t}. 
\end{align*}
Hence, $p_\lambda^{(0)}(w(\varepsilon;\lambda,t))$ converges uniformly to $p_\lambda^{(0)}(w(0;\lambda,t))$ on $\mathbb{C}$ as $\varepsilon$ tends to zero,  thanks to the uniform convergence of $w(\varepsilon;\lambda,t)$ to $w(0;\lambda,t)$ proved in Lemma \ref{lemma:sub-w-uniform-convergence}.
\end{proof}

\begin{definition}
Define the map $\Phi_{t, \gamma}:\mathbb{C}\rightarrow\mathbb{C}$ by 
\begin{equation}
\label{eqn:def-Phi-t-gamma}
 	\Phi_{t,\gamma} (\lambda)= \lambda+\gamma\cdot p_\lambda^{(0)}( w(0;\lambda,t)  ),
	\qquad \lambda\in\mathbb{C}.
\end{equation}
\end{definition}

\begin{theorem}
	\label{thm:main-FK-det-x-t-gamma-0}
	Let $y=g_{t,\gamma}$ and ${x}$ be a random variable free from $x$. 
	Assume that $\Phi_{t, \gamma}$ is non-singular at $\lambda\in\Xi_t$ (the Jacobian of $\Phi_{t, \gamma}$ is invertible at $\lambda$). 
	Then, $(z,\varepsilon)\mapsto\Delta\big(|{x}+g_{t,\gamma}-z\unit|^2+\varepsilon^2\big)$ has a real analytic extension in a neighborhood of $(\Phi_{t,\gamma} (\lambda),0)$. Moreover, 
	\begin{equation}
	\label{eqn:relation-FK-det-c-g-0}
	\Delta({x}+g_{t,\gamma}-z\unit)^2
	=\Delta({x}+c_t-\lambda\unit)^2\exp \left[ \Re ( \gamma (p_\lambda^{(0)}(w(0;\lambda,t)))^2 ) \right]
	\end{equation}
	where $z=\Phi_{t,\gamma} (\lambda)$.
\end{theorem}
\begin{proof}
	Notice that the right hand side of \eqref{eqn:main-FK-det-0-V2} is a real analytic function of $\lambda$ by the analyticity of $w(0;\lambda,t)$ in Lemma \ref{lemma:s-epsilon-t-limit}. By Lemma \ref{lemma:Jacobian-circular-0}, the Jacobian of the map $(\lambda,\varepsilon)\mapsto\Psi_{c,0}(\lambda,\varepsilon)=(\lambda,w(\varepsilon;\lambda,t))$ is invertible. An application of inverse function theorem implies that the map 
	\[
	(\lambda,\varepsilon)\mapsto \Delta\big(|{x}+c_t-\lambda\unit|^2+\varepsilon^2 \big)
	\]
	has a real analytic extension in some neighborhood of $(\lambda,0)$. By choosing $\gamma=0$ and an arbitrary $\gamma$ in Lemma \ref{lemma:main-FK-det-gt}, we can obtain that 
	\begin{equation}
	\label{eqn:relation-FK-det-c-g-in-proof}
	\Delta\big(|{x}+g_{t,\gamma}-z\unit|^2+\varepsilon^2\big) 
	=\Delta\big(|{x}+c_t-\lambda\unit|^2+\varepsilon^2 \big)
	\exp \left[ \Re ( \gamma (p_\lambda^{(0)}(w(\varepsilon;\lambda,t)))^2 )\right],
	\end{equation}
	where 
	\[
	z=\lambda+\gamma\cdot p_\lambda^{(0)}(w(\varepsilon;\lambda,t)). 
	\]
	By Lemma \ref{lemma:convergence-p-lambda}, 
	\[
	\lim_{\varepsilon\rightarrow 0^+} p_\lambda^{(0)}(w(\varepsilon;\lambda,t))=p_\lambda^{(0)}(w(0;\lambda,t)).
	\]
	Hence, 
	$\lim_{\varepsilon\rightarrow 0^+}\Psi_{c,g}(\lambda,\varepsilon)=(\Phi_{t, \gamma}(\lambda),0)$. 
	Since $w(0;\lambda,t)>0$, then the assumption that $\Phi_{t,\gamma}$ is non-singular at $\lambda\in\Xi_t$ implies that 
	the Jacobian of the map $(\lambda,\varepsilon)\mapsto\Psi_{c,g}(\lambda,\varepsilon)=(z,\varepsilon)$ at $(\lambda,0)$ is invertible. 
	Hence, the analyticity of $\Delta\big(|{x}+c_t-\lambda\unit|^2+\varepsilon^2 \big)$ in some neighborhood of $(\lambda,0)$ implies that the map $(z,\varepsilon)\mapsto\Delta\big(|{x}+g_{t,\gamma}-z\unit|^2+\varepsilon^2\big)$ has a real analytic extension in a neighborhood of $(\Phi_{t,\gamma} (\lambda),0)$
	by inverse function theorem. Hence, \eqref{eqn:relation-FK-det-c-g-0} follows from \eqref{eqn:relation-FK-det-c-g-in-proof} by letting $\varepsilon$ tend to zero.
\end{proof}

\begin{corollary}
	\label{cor:FK-det-comparison-H}
	If $\lambda\in\Xi_t$, assume that $\Phi_{t,\gamma}$ is non-singular at $\lambda\in\Xi_t$, then
	\begin{align}
	\Delta \big( ({x}+&y-{z}\unit)^* ({x}+y-{z}\unit))\big)
	={\Delta\big( ({x}-\lambda\unit)^*({x}-\lambda\unit)+w(0; \lambda,t)^2\unit \big)} \nonumber\\
	\times &  
	\exp \left[   
	\frac{1}{2}\cdot
	\bigg( \gamma\cdot ( p_\lambda^{(0)}(w(0;\lambda,t)))^2+\overline{\gamma}\cdot (p_{\overline{\lambda}}^{(0)}(w(0;\lambda,t))^2\bigg)
	-\frac{ w(0;\lambda,t)^2}{t}	
	\right],
	\end{align}
	where $z=\lambda+\gamma\cdot p_\lambda^{(0)}(w(0;\lambda,t))$.
\end{corollary}

\section{Brown measure of addition with a circular operator}
\label{section:Brown-free-additiveBM}

In this section, we show that the Brown measure of ${x}+c_t$ has no atom and it is absolutely continuous with respect to Lebesgue measure with strictly positive and real analytic density in the open set $\Xi_t$, and the density formula can be expressed explicitly in terms of the function $w(0;\lambda,t)$.

\subsection{The density formula in the domain $\Xi_t$}
We first study the limits of $p_\lambda^{c,(t)}(\varepsilon)$ as $\varepsilon$ tends to zero.

\begin{lemma}
	\label{lemma:partial-lambda-c-t-limit}
	For any $t>0$.
	The function $\lambda\mapsto S({x}+c_t, \lambda,0)$ is a real analytic function for $\lambda\in\Xi_t$, and
	we have 
	\[
	p_\lambda^{c,(t)}(0)=p_\lambda^{(0)}(w(0;\lambda,t)),
	\]
\end{lemma}
\begin{proof}
	By Corollary \ref{thm:sub-comparison-Ct}, we have 
	\[
	p_\lambda^{c,(t)}(\varepsilon)=p_\lambda^{(0)}(w(\varepsilon;\lambda,t)), 
	\]
	Recall that, by Lemma \ref{lemma:s-epsilon-t-limit}, for $\lambda\in \Xi_t$, $\lim_{\varepsilon\rightarrow 0}w(\varepsilon;\lambda,t)=w(0;\lambda,t)\in (0,\infty)$. The result then follows by letting 
	$\varepsilon$ tend to zero.
\end{proof}

The following result generalizes \cite[Theorem 3.10]{HoZhong2020Brown} where ${x}$ is assumed that to be self-adjoint and \cite[Theorem 1.4]{Bordenave-Caputo-Chafai-cpam2014} where ${x}$ is assumed to be a Gaussian distributed normal operator (their techniques extend to the case when ${x}$ is a normal operator as in \cite{Bordenave-Capitaine-cpam2016}).
\begin{theorem}
	\label{thm:BrownFormula-x0-ct-general}
	The Brown measure is absolutely continuous with respect to Lebesgue measure in the open set $\Xi_t$.
	The density of the Brown measure at $\lambda\in \Xi_t$ is given by
	\begin{equation}
	\label{density-x0-ct-formula-2}
	\frac{1}{\pi} \left( \frac{1}{t}-  \frac{\partial}{\partial \overline{\lambda}} \bigg( \phi \big( {x}^* ( ({x}-\lambda\unit)({x}-\lambda\unit)^*+w(0; \lambda,t)^2\unit )^{-1} \big) \bigg) \right)
	\end{equation}
	where $w=w(0; \lambda,t)$ is determined by 
	\[
	\phi ( ({x}-\lambda\unit)^*({x}-\lambda\unit)+w(0; \lambda,t)^2\unit )^{-1} )=\frac{1}{t}.
	\]
	It can also be expressed as
	\begin{equation}
	\label{eqn:Brown-density-circular-positive}
	\frac{1}{\pi}\left( \frac{|\phi((\lambda\unit-{x})(h^{-1})^2)|^2}{\phi((h^{-1})^2)}
	+w(0;\lambda,t)^2\phi(h^{-1}k^{-1})  \right)
	\end{equation}
	where $h=h(\lambda,w(0;\lambda,t))$ and $k=k(\lambda,w(0;\lambda,t))$ for 
	\[
	h(\lambda,w)=(\lambda\unit-{x})^*(\lambda\unit-{x})+w^2
	\]
	and
	\[
	k(\lambda,w)=(\lambda\unit-{x})(\lambda\unit-{x})^*+w^2.
	\]
	In particular, the density of the Brown measure of ${x}+c_t$ is strictly positive in the set $\Xi_t$. 
\end{theorem}

After the present article was released to arXiv, in a joint work with Belinschi and Yin \cite[Section 7]{BelinschiYinZhong2021Brown}, we obtain the following strengthened result.

\begin{theorem}\cite[Theorem 7.10 and Lemma 7.11]{BelinschiYinZhong2021Brown}
\label{thm:main-theorem-from-BYZ2021}
The Brown measure $\mu_{x+c_t}$ is absolutely continuous with respect to Lebesgue measure on the complex plane. The density functions of both $\mu_{x+c_t,\varepsilon}$ and $\mu_{x+c_t}$ 
are bounded by $1/\pi{t}$. 
\end{theorem}

\begin{proof}[Proof of Theorem \ref{thm:BrownFormula-x0-ct-general}]
	For $\lambda\in\Xi_t$, we have $\phi \left[\big( ({x}-\lambda\unit)^*({x}-\lambda )\big)^{-1}\right]>\frac{1}{t}$, and
	\begin{align*}
	& \Delta \big( ({x}+c_t-\lambda\unit)^*({x}+c_t-\lambda\unit) \big)\\
	&\qquad=\Delta\big( ({x}-\lambda\unit)^*({x}-\lambda\unit) +w(0; \lambda,t)^2\unit \big)
	\cdot \exp\left[-t (q_{w(0;\lambda,t)}^{(0)}(\lambda))^2 \right]  	
	\end{align*}
	by Theorem \ref{thm:main-FK-det-ct-0}.
	In addition, $w(0; \lambda,t)\in (0,\infty)$ and $(a,b)\mapsto w(0; \lambda,t)$ for $\lambda=a+{\i}b$ is real analytic by Lemma \ref{lemma:s-epsilon-t-limit}. 
	Hence, $\lambda\mapsto\log\Delta({x}+c_t-\lambda\unit)$ is real analytic. 
	We put 
	\[
	g(\lambda)= \log \Delta \big( ({x}+c_t-\lambda\unit)^*({x}+c_t-\lambda\unit) \big)
	\]
	The Brown measure can be calculated as
	\[
	d\mu_{{x}+c_t}(\lambda)=\frac{1}{\pi}\frac{\partial^2}{\partial\overline{\lambda
		}\partial\lambda} g(\lambda).
	\]
	where the Laplacian can be calculated in the usual sense.

	By the identity $x^*(xx^*+\varepsilon\unit)^{-1}=(x^*x+\varepsilon\unit)^{-1}x^*$ and tracial property, then Lemma \ref{lemma:partial-lambda-c-t-limit} is equivalent to 
	\[
	\frac{\partial}{\partial\lambda} g(\lambda)=\phi\left(\frac{\partial h}{\partial\lambda} k^{-1}\right)=\phi\left(h^{-1}\frac{\partial h}{\partial\lambda}\right)
	\]
	where $\frac{\partial h}{\partial\lambda}=(\lambda\unit-{x})^*$. We can continue to take the derivative directly 
	\begin{align*}
	\frac{\partial^2}{\partial\overline{\lambda
		}\partial\lambda}g(\lambda)
	&=\frac{\partial}{\partial\overline{\lambda}}\phi\left(h^{-1}\frac{\partial h}{\partial\lambda}\right)\\
	&=\frac{\partial}{\partial\overline{\lambda}}\phi \big( (\lambda\unit-{x})^* ( ({x}-\lambda\unit)^*({x}-\lambda\unit)+w(0; \lambda,t)^2\unit )^{-1} \big)\\
	&=\frac{\partial}{\partial\overline{\lambda}}\phi \bigg( \frac{\overline\lambda}{t}-\phi \big( {x}^* ( ({x}-\lambda\unit)^*({x}-\lambda\unit)+w(0; \lambda,t)^2\unit )^{-1} \big) \bigg)\\
	&= \frac{1}{t}-  \frac{\partial}{\partial \overline{\lambda}} \bigg( \phi \big( {x}^* ( ({x}-\lambda\unit)^*({x}-\lambda\unit)+w(0; \lambda,t)^2\unit )^{-1} \big) \bigg).
	\end{align*}
	Then the first formula \eqref{density-x0-ct-formula-2} is established. 
	
	We adapt the calculation in \cite[Lemma 2.8]{HaagerupSchultz2007} to get another form of the density formula. By \cite[Lemma 3.2]{HaagerupT2005Annals}, since $g(\lambda)$ is a real analytic function of $\lambda$, we have
	\begin{align}
	\frac{\partial^2}{\partial\overline{\lambda
		}\partial\lambda}g(\lambda)
	&=\frac{\partial}{\partial\overline{\lambda}}\phi\left(h^{-1}\frac{\partial h}{\partial\lambda}\right)\nonumber\\
	&=\phi\left( -h^{-1} \left( \frac{\partial h}{\partial\overline{\lambda}}+2w \frac{\partial {w}}{\partial\overline{\lambda}}  \right) h^{-1} \frac{\partial h}{\partial\lambda}   + h^{-1} \frac{\partial^2 h}{\partial\overline{\lambda} \partial\lambda}  \right) \label{eqn:4.16-in-proof}\\
	&=\phi\left( -h^{-1}  \frac{\partial h}{\partial\overline{\lambda}} h^{-1} \frac{\partial h}{\partial\lambda}   + h^{-1} \frac{\partial^2 h}{\partial\overline{\lambda} \partial\lambda}  \right)
	-2w \frac{\partial {w}}{\partial\overline{\lambda}} 
	\phi\left( \frac{\partial h}{\partial {\lambda}} (h^{-1})^2  \right)\nonumber
	\end{align}
	and
	\[
	\frac{\partial h}{\partial\overline{\lambda}}= \frac{\partial h(\lambda,w)}{\partial\overline{\lambda}} =\lambda\unit -{x}; 
	\]
	and 
	\[
	\frac{\partial^2 h}{\partial\overline{\lambda} \partial\lambda} 
	=\frac{\partial}{\partial\overline{\lambda}} (\lambda\unit-{x})^*=1.
	\]
	
	We now apply the identity $x(x^*x+\varepsilon\unit)^{-1}=(xx^*+\varepsilon\unit)^{-1}x$
	to $x=\lambda\unit-{x}$ and $\varepsilon=w^2$, we find that
	\begin{align*}
	-\frac{\partial h}{\partial\overline{\lambda}} h^{-1} \frac{\partial h}{\partial\lambda}   +  \frac{\partial^2 h}{\partial\overline{\lambda} \partial\lambda}
	&=1-x(x^*x+w^2\unit)^{-1}x^*\\
	&=1-(xx^*+w^2\unit)^{-1}xx^*\\
	&=w^2 (xx^*+w^2\unit)^{-1}\\
	&=w^2 k^{-1}.
	\end{align*}
	Now, $w(0;\lambda,t)$ is determined by 
	\[
	\phi ( ({x}-\lambda\unit)^*({x}-\lambda\unit)+w(0; \lambda,t)^2\unit )^{-1} )=\frac{1}{t}
	\]
	which can be rewritten as
	\[
	\phi( h^{-1} )=\frac{1}{t}.
	\]
	Take implicit differentiation $\frac{\partial}{\partial\overline{\lambda}}$ and apply again \cite[Lemma 3.2]{HaagerupT2005Annals}, we then obtain
	\[
	\phi\left(h^{-1} \left( \frac{\partial h}{\partial\overline{\lambda}}+ 2w \frac{\partial {w}}{\partial \overline{\lambda}} \right)  h^{-1}   \right)=0,
	\]
	where $\frac{\partial h}{\partial\overline\lambda}=\lambda\unit-{x}$.
	This implies 
	\[
	- 2w \frac{\partial {w}}{\partial \overline{\lambda}}=\frac{ \phi \left( \frac{\partial h}{\partial\overline{\lambda} } (h^{-1})^2  \right) }{ \phi((h^{-1})^2)  }.
	\]
	By the tracial property, we observe that 
	\[
	\phi \left( \frac{\partial h}{\partial\overline{\lambda} } (h^{-1})^2  \right)
	=\overline{ \phi \left( \frac{\partial h}{\partial{\lambda} } (h^{-1})^2  \right)}.
	\]
	We therefore can continue to simplify \eqref{eqn:4.16-in-proof} as
	\begin{align*}
	\frac{\partial^2}{\partial\overline{\lambda
		}\partial\lambda}g(\lambda)
	&=\phi\left( -h^{-1}  \frac{\partial h}{\partial\overline{\lambda}} h^{-1} \frac{\partial h}{\partial\lambda}   + h^{-1} \frac{\partial^2 h}{\partial\overline{\lambda} \partial\lambda}  \right)
	-2w \frac{\partial {w}}{\partial\overline{\lambda}} 
	\phi\left( \frac{\partial h}{\partial {\lambda}} (h^{-1})^2  \right)\\
	&=w^2\phi(h^{-1} k^{-1})+
	\frac{ \left| \phi \left( \frac{\partial h}{\partial\overline{\lambda}  } (h^{-1})^2  \right) \right|^2}{ \phi((h^{-1})^2)  }\\
	&=w^2\phi(h^{-1/2} k^{-1}h^{-1/2})+
	\frac{ \left| \phi \left((\lambda\unit-{x}) (h^{-1})^2  \right) \right|^2}{ \phi((h^{-1})^2)  }>0
	\end{align*}
	for any $\lambda\in\Xi_t$.
	This finishes the proof. 
\end{proof}

\subsection{The support of the Brown measure of ${x}+c_t$}

Theorem \ref{thm:BrownFormula-x0-ct-general} does not tell us the Brown measure of ${x}+c_t$ outside the open set $\Xi_t$. We will show that the Brown measure of ${x}+c_t$ is supported in the closure of $\Xi_t$. 
\begin{lemma}
	\label{lemma:Brown-coincide-outside-Xi}
	The Brown measures of ${x}$ and ${x}+c_t$ coincide in the complement of the closure $\overline{\Xi_t}$. That is
	\[
	\mu_{{x}}|_{(\overline{\Xi_t})^c} = \mu_{{x}+c_t}|_{(\overline{\Xi_t})^c}.
	\]
	In particular, $\mu_{{x}+c_t}(\overline{\Xi_t})=1$ if and only if $\mu_{{x}}(\overline{\Xi_t})=1$.
\end{lemma}
\begin{proof}
	For $\lambda\in (\Xi_t)^c$, by Theorem \ref{thm:main-FK-det-ct-0}, we have 
	\[
	\log\Delta({x}-\lambda\unit)=\log\Delta({x}+c_t-\lambda\unit).
	\]
	Hence, 
	\[
	\int_\mathbb{C}\log |z-\lambda|d\mu_{{x}}(z)=\int_\mathbb{C}\log |z-\lambda|d\mu_{{x}+c_t}(z)
	\]
	for any $\lambda\in (\Xi_t)^c$. Then two Brown measures coincide in the open set $(\overline{\Xi_t})^c$ due
	to the Unicity Theorem of logarithmic potential (see \cite[Theorem 2.1 in Chapter II]{SaffTotik1997book}).
\end{proof}

We are grateful to Hari Bercovici for providing us a proof of the following result which is a refinement of an argument in \cite[Proposition 4.6]{HaagerupLarsen2000}.
\begin{lemma}
	\label{lemma:I-inverse2-bonded}
	Let $\mu$ be a finite Borel measure on $\mathbb{C}$. Define $I: \mathbb{C}\rightarrow (0,\infty]$ by
	\[
	I_\mu(\lambda)=\int_\mathbb{C}\frac{1}{\vert \lambda- z\vert^2}d\mu(z).
	\] 
	Then $I_\mu(\lambda)$ is infinite almost everywhere relative to $\mu$.  
\end{lemma}
\begin{proof}
Without loss of generality, we may assume that $\mu$ is compactly supported. 
For every $\lambda\in\mathbb{C}$ there exists $\lambda'\in\text{supp}(\mu)$
such that $I_\mu(\lambda)\le 4I_\mu(\lambda')$. Indeed, take $\lambda'$ to be the
closest point to $\lambda$, so
\[
|\lambda-\lambda'|\le|\lambda-z|,\quad z\in{\rm supp}(\mu).
\]
Then
\[
|\lambda'-z|\le|\lambda-\lambda'|+|\lambda-z|\le2|\lambda-z|,\quad z\in{\rm supp}(\mu),
\]
and this yields the desired inequality. Thus, if $I(\lambda)$ is
bounded on supp$(\mu)$ then it is bounded everywhere. 

If $\{\lambda\in\text{supp}(\mu):I(\lambda)<+\infty\}$ has positive
$\mu$ measure, we may further assume that there exists $C>0$ such that
	\[
	I(\lambda)<C, \qquad \lambda\in {\rm supp}(\mu).
	\]
Otherwise, we can replace $\mu$ by a restriction $\mu$ to some
compact set. Suppose this boundedness actually happens, and choose $R>0$ large
enough that $|z|+1<R$ for every $z\in{\rm supp}(\mu)$. We observe that
\[
+\infty>\int_{|\lambda|<R}I(\lambda)\,dm_2(\lambda)\ge\int_{{\rm supp}(\mu)}\left[\int_{|\lambda-z|<1}\frac{dm_2(\lambda)}{|\lambda-z|^{2}}\right]\,d\mu(z),
\]
where $m_2$ denotes Lebesgue measure on the complex plane. 
Since the integral inside the bracket is $+\infty$ for every $z$,
we conclude that $\mu$ must be the zero measure. Therefore, $I_\mu(\lambda)$ is infinite almost everywhere relative to $\mu$. 
\end{proof}

The following result  provides a new characterization of the support of the Brown measure of an arbitrary operator. The result also improves \cite[Theorem 1.2]{HallKemp2019} where a condition of local boundedness of $(T-\lambda\unit)^{-1}$ in $L^4$-norm is required. 
\begin{theorem}
	\label{thm:support-Brown-general}
	Let $\lambda\in \mathbb{C}$ and $T\in \mathcal{M}$. If $\phi(|T-\lambda|^{-2})<\infty$, then 
	\[
	\int_\mathbb{C}\frac{1}{|z-\lambda|^2}d\mu_T(z)\leq \phi(|T-\lambda|^{-2}). 
	\]
	If  $\phi(|T-\lambda|^{-2})<\infty$ in some neighborhood of $\lambda_0$, then $\lambda_0\notin \supp(\mu_T)$.
\end{theorem}
\begin{proof}
	Observe that for any $t\in (0,1)$, we have $2\log t>-\frac{1}{t^2}$. Hence, 
	\[
	\begin{aligned}
	2 \int_0^1 \log t d\mu_{|T-\lambda|}(t)
	&>-\int_0^1 \frac{1}{t^2}d\mu_{|T-\lambda|}(t)\\
	&>-\int_0^\infty \frac{1}{t^2}d\mu_{|T-\lambda|}(t)\\
	&=-\phi(|T-\lambda|^{-2})>-\infty. 
	\end{aligned}
	\]
	By \cite[Proposition 2.16]{HaagerupSchultz2007}, the operator $T-\lambda$ has an inverse $(T-\lambda)^{-1}$ that is possibly unbounded operator affiliated with $\mathcal{M}$ and the Brown measure $\mu_{(T-\lambda)^{-1}}$ of $(T-\lambda)^{-1}$ can be defined. Moreover, $\mu_{(T-\lambda)^{-1}}$ is the push-forward measure of $\mu_T$ via the map $z\mapsto (z-\lambda)^{-1}$. Hence, by \cite[Theorem 2.19]{HaagerupSchultz2007}, we obtain
	\begin{equation}
	\label{eqn:inequality-L2-inverse}
	\begin{aligned}
	\int_\mathbb{C}\frac{1}{|z-\lambda|^2}d\mu_T(z) & =\int_\mathbb{C} |z|^2 d\mu_{(T-\lambda)^{-1}}(z)\leq ||(T-\lambda)^{-1}||_2^2=\phi(|T-\lambda|^{-2}).   
	\end{aligned}
	\end{equation}
	If  $\phi(|T-\lambda|^{-2})<\infty$ in some neighborhood of $\lambda_0$.
	It then follows by Lemma \ref{lemma:I-inverse2-bonded} that $\mu_T(U)=0$ thanks to \eqref{eqn:inequality-L2-inverse}. Hence, $\lambda_0 \notin \supp(\mu_T)$. 
\end{proof}

\begin{theorem}
	\label{thm:support-x0-c_t-general}
	For any ${x}$ $*$-free from $c_t$, the Brown measure $\mu_{{x}+c_t}$ of ${x}+c_t$ has no atom and is supported in the closure $\overline{\Xi_t}$.
\end{theorem}
\begin{proof}
	Given $\lambda\in (\Xi_t)^c$, we have 
	\[
	\log\Delta({x}-\lambda\unit)=\log\Delta({x}+c_t-\lambda\unit),
	\]
	which yields that
	\[
	\int_0^\infty\log |t|d\mu_{|{x}-\lambda\unit|}(t)
	= \int_0^\infty\log |t|d\mu_{|{x}+c_t-\lambda\unit|}(t).
	\]
	In addition, we have $\phi(|{x}-\lambda\unit|^{-2})\leq 1/t$, and as in the proof of Theorem \ref{thm:support-Brown-general}, we have 
	\[
	\int_0^1 \log t\, d\mu_{|{x}-\lambda\unit|}(t)>-\infty.
	\]
	Hence,
	\[
	\int_0^1 \log t \, d\mu_{|{x}+c_t-\lambda\unit|}(t)>-\infty.
	\]
	It follows by \cite[Proposition 2.16]{HaagerupSchultz2007} that $\mu_{{x}+c_t-\lambda}(\{0\})=0$. Hence $\mu_{{x}+c_t}(\{\lambda\})=0$ for any $\lambda\in \mathbb{C}\backslash\Xi_t$, since $\mu_{{x}+c_t-\lambda\unit}$ is the translation of $\mu_{{x}+c_t}$ by $-\lambda$. Recall that $\mu_{{x}+c_t}$ is absolutely continuous in $\Xi_t$. It follows that $\mu_{{x}+c_t}$ has no atom in $\mathbb{C}$. 
	
	For any $\lambda\in (\Xi_t)^c$, there is some neighborhood $U$ of $\lambda$ such that $\phi(|{x}-z\unit|^{-2})\leq 1/t$ for $z\in U$. By Theorem \ref{thm:support-Brown-general}, $\mu_{{x}}(U)=0$. Hence, $\mu_{{x}}( (\Xi_t)^c )=0$.
	By Lemma \ref{lemma:Brown-coincide-outside-Xi}, the Brown measure of ${x}+c_t$ coincides with the Brown measure of ${x}$ within the open set $(\overline{\Xi_t})^c$. Hence, $\mu_{{x}+c_t}((\overline{\Xi_t})^c)=0$ and thus the Brown measure is supported in the closure $\overline{\Xi_t}$, thanks to Theorem \ref{thm:BrownFormula-x0-ct-general}.
\end{proof}

\section{Brown measure of addition with an elliptic operator}
\label{section:Brown-addition-elliptic}
In this section, we show the Brown measure of $\mu_{{x}+g_{t,\gamma}}$ is the push-forward measure of $\mu_{{x}+c_t}$ under the map $\Phi_{t,\gamma}$ defined in \eqref{eqn:def-Phi-t-gamma}, which is
\[
      \Phi_{t,\gamma} (\lambda)= \lambda+\gamma\cdot p_\lambda^{(0)}( w(0;\lambda,t)  ),
   \qquad \lambda\in\mathbb{C}.
\]
Recall that for $\lambda\notin \Xi_t$ we have $w(0;\lambda,t)=0$ and $p_\lambda^{(0)}(0)=\lim_{\varepsilon\rightarrow 0}p_\lambda^{(0)}(\varepsilon)$. The existence of this limit is due to the fact that $\phi(|\lambda-x|^{-2})\leq 1/t$ for $\lambda\notin \Xi_t$. 

We define the function $\Phi^{(\varepsilon)}_{t,\gamma}$ on $\mathbb{C}$ by
\begin{align*}
 \label{defn:Phi-t-gamma-epsitlon-v0}
 \Phi^{(\varepsilon)}_{t,\gamma} (\lambda)&= \lambda+\gamma\cdot \frac{\partial S({x}+c_t,\lambda,\varepsilon)}{\partial \lambda}\\
 &=\lambda+\gamma\cdot \phi\left( (\lambda-x-c_t)^*(((\lambda-x-c_t)(\lambda-x-c_t)^*+\varepsilon^2)^{-1}) \right)
\end{align*}
By the subordination relation \eqref{eqn:subordination-Cauchy-entries}, we can also write it as 
\begin{equation}
\label{defn:Phi-t-gamma-epsitlon}
\Phi^{(\varepsilon)}_{t,\gamma} (\lambda)= \lambda+\gamma\cdot p_\lambda^{(0)}( w(\varepsilon;\lambda,t)  ),
\qquad \lambda\in\mathbb{C}
\end{equation}
where 
\[
p_\lambda^{(0)}( w(\varepsilon;\lambda,t)  ) =
\phi\bigg[ (\lambda-{x})^*\big( (\lambda-{x})(\lambda-{x})^*+w(\varepsilon; \lambda,t)^2 \unit\big)^{-1}\bigg].
\]
By Notation \ref{notation:partial-derivatives-lambda-c-x}, these maps can also be expressed as
\[
  \Phi^{(\varepsilon)}_{t,\gamma}(\lambda)=\lambda+\gamma\cdot \frac{\partial S}{\partial \lambda}({x}, \lambda,w(\varepsilon;\lambda,t)), \qquad \lambda\in\mathbb{C}.
\]
For $\lambda\in\Xi_t$, since $w(0;\lambda,t)>0$, we have 
\[
  \Phi_{t,\gamma}(\lambda)=\lambda+\gamma\cdot \frac{\partial S}{\partial \lambda}({x}, \lambda,w(0;\lambda,t)),  \qquad \lambda\in\Xi_t.
\]
\begin{lemma}
	\label{lemma:regularization-Phi-uniform-convergence}
The function $\Phi^{(\varepsilon)}_{t,\gamma}$ converges to $\Phi_{t,\gamma}$ uniformly on $\mathbb{C}$ as $\varepsilon$ tends to zero. 
\end{lemma}
\begin{proof}
This follows from Lemma \ref{lemma:convergence-p-lambda} which says that $p_\lambda^{(0)}(w(\varepsilon;\lambda,t))$ converges uniformly to $p_\lambda^{(0)}(w(0;\lambda,t))$ as $\varepsilon$ tends to zero.
\end{proof}
Our strategy is to show first that the regularized Brown measure $\mu_{{x}+g_{t,\gamma},\varepsilon}$ is the push-forward measure of the regularized Brown measure of $\mu_{{x}+c_t,\varepsilon}$ under that map $\Phi^{(\varepsilon)}_{t,\gamma}$. We then show that $\Phi^{(\varepsilon)}_{t,\gamma}$ converges to $\Phi_{t,\gamma}$ uniformly on $\mathbb{C}$ as $\varepsilon$ tends to zero. 
\subsection{Regularized Brown measure and the regularized push-forward map}
\label{section:regularized-pushforward-map}

\begin{proposition}
	\label{prop:regularization-Phi-one2one-map}
	The map $\Phi^{(\varepsilon)}_{t,\gamma}$ is a $C^\infty$-self-diffeomorphism of $\mathbb{C}$. 
\end{proposition}
\begin{proof}
	It is clear that $\Phi^{(\varepsilon)}_{t,\gamma}$ is a real analytic map of $\mathbb{C}$. 
	Assume that $z=\Phi^{(\varepsilon)}_{t,\gamma}(\lambda_1)=\Phi^{(\varepsilon)}_{t,\gamma}(\lambda_2)$ for $\lambda_1, \lambda_2\in\Xi_t$. Then, by the subordination relation in Theorem \ref{thm:sub-X0-gt-formula}, we have, for $i=1,2$, 
	\[
	p_z^{g,(t,\gamma)}(\varepsilon)=p_{\lambda_i}^{c,(t)}(\varepsilon)=p_{\lambda_i}^{(0)}(w(\varepsilon;\lambda_i,t)),
	\]
	where $z=\Phi^{(\varepsilon)}_{t,\gamma}(\lambda_i) (i=1,2)$. 
	Rewrite the map $\Phi^{(\varepsilon)}_{t,\gamma}$ as 
	\[
	\Phi^{(\varepsilon)}_{t,\gamma}(\lambda)=\lambda+\gamma\cdot p_\lambda^{(0)}(w(\varepsilon;\lambda,t))=\lambda+ \gamma\cdot p_z^{g,(t,\gamma)}(\varepsilon).
	\]
	Then the condition $z=\Phi^{(\varepsilon)}_{t,\gamma}(\lambda_1)=\Phi^{(\varepsilon)}_{t,\gamma}(\lambda_2)$ yields that
	$\lambda_1=\lambda_2$. Hence, $\Phi^{(\varepsilon)}_{t,\gamma}$ is one-to-one in $\mathbb{C}$. 
	
	It follows from Corollary \ref{cor:Omega-1-formula} that the map $\Phi^{(\varepsilon)}_{t,\gamma}$ is also surjective. Its inverse map (see Theorem \ref{thm:sub-X0-gt-formula}) is 
	\[
	  \big(\Phi^{(\varepsilon)}_{t,\gamma}\big)^{-1}(z)=z-\gamma\cdot p_z^{g,(t,\gamma)}(\varepsilon),
	\]
	which is also a $C^\infty$ map of $\mathbb{C}$. 
\end{proof}

\begin{remark}
Proposition \ref{prop:regularization-Phi-one2one-map} also follows directly from  Corollary \ref{cor:Omega-1-formula}.
\end{remark}

\begin{theorem}
	\label{thm:main-push-forward-general-regularized}
	The regularized Brown measure $\mu_{{x}+g_{t,\gamma},\varepsilon}$ is the push-forward measure of the regularized Brown measure of $\mu_{{x}+c_t,\varepsilon}$ under the map $\Phi^{(\varepsilon)}_{t,\gamma}$. That is, 
	\[
	\mu_{{x}+g_{t,\gamma},\varepsilon}(\cdot)=\mu_{{x}+c_t,\varepsilon}\left( \big(\Phi^{(\varepsilon)}_{t,\gamma}\big)^{-1}(\cdot) \right).
	\]
\end{theorem}
\begin{proof}
	Let $\mathrm{\Gamma}\subset \mathbb{C}$ be a simply connected domain with piecewise smooth boundary. Since the regularized map $\Phi^{(\varepsilon)}_{t,\gamma}$ is a self-diffeomorphism of $\mathbb{C}$, it suffices to show that
	\begin{equation}
	\label{eqn:regular-push-forward-ct-x-t-gamma}
	\mu_{{x}+c_t,\varepsilon}(\mathrm{\Gamma})=\mu_{{x}+g_{t,\gamma},\varepsilon}(\Phi^{(\varepsilon)}_{t,\gamma}(\mathrm{\Gamma}))).
	\end{equation}
	The domain $\Phi^{(\varepsilon)}_{t,\gamma}(\mathrm{\Gamma}))$ is also a simply connected domain with piecewise smooth boundary. 
	For $\lambda\in\mathbb{C}$, set
	\[
	z= \Phi^{(\varepsilon)}_{t,\gamma} (\lambda)= \lambda+\gamma\cdot p_\lambda^{(0)}( w(\varepsilon;\lambda,t)  ).
	\]
	By Theorem \ref{thm:sub-X0-gt-formula} and Notation \ref{notation:partial-derivatives-lambda-c-x}, we have
	\[
	p_\lambda^{(0)}( w(\varepsilon;\lambda,t)  )=p_\lambda^{c,(t)}(\varepsilon)=\frac{\partial S({x}+c_t,\lambda,\varepsilon)}{\partial \lambda}=p_z^{g, (t,\gamma)}(\varepsilon) =\frac{\partial S({x}+g_{t,\gamma},z,\varepsilon)}{\partial z}.
	\] 
	Write $\lambda=\lambda_1+{\i}\lambda_2\in\mathbb{C}$ and $z=z_1+{\i}z_2 \in \mathbb{C}$, and we denote the vector fields
	\[
	P^c(\lambda_1, \lambda_2)=\Re (p_\lambda^{c,(t)}(\varepsilon)), 
	\qquad Q^c(\lambda_1, \lambda_2)=-\Im (p_\lambda^{c,(t)}(\varepsilon));
	\]
	and
	\[
	P^g(z_1, z_2)=\Re (p_z^{g,(t,\gamma)}(\varepsilon)), 
	\qquad Q^g(z_1, z_2)=-\Im (p_z^{g, (t,\gamma)}(\varepsilon)).
	\]
	We then have 
	\[
	p_\lambda^{c,(t)}(\varepsilon)
	=P^c(\lambda_1,\lambda_2)-iQ^c(\lambda_1,\lambda_2),
	\]
	and
	\[
	P^c(\lambda_1, \lambda_2)=P^g(z_1, z_2),
	\qquad Q^c(\lambda_1, \lambda_2)=Q^g(z_1, z_2).
	\]
	Let $\gamma=\gamma_1+{\i}\gamma_2$, then we have 
	\begin{equation}
	\label{eqn:formula-z1-z2}
	\begin{aligned}
	z_1&=\lambda_1+\gamma_1 P^c(\lambda_1, \lambda_2) +\gamma_2 Q^{c}(\lambda_1,\lambda_2);\\
	z_2&=\lambda_2+\gamma_2 P^c(\lambda_1, \lambda_2)-\gamma_1 Q^{c}(\lambda_1,\lambda_2).
	\end{aligned}
	\end{equation}
	Denote the differential form $\alpha$ in $\mathbb{C}$ as
	\begin{equation}
	\label{eqn:1-form-alpha}
	\alpha=-Q^g dz_1+P^g dz_2.
	\end{equation}
	Since $\Phi^{(\varepsilon)}_{t,\gamma}$ is one-to-one, we can change variables from $z$ to $\lambda$, which really means that we pull back the $1$-form $\alpha$ by the map $\Phi^{(\varepsilon)}_{t,\gamma}$. Letting $\beta$ be the pulled-back form, and also using formulas \eqref{eqn:formula-z1-z2} for $z_1, z_2$, we get 
	\[
	\begin{aligned}
	\beta&=-Q^c d(\lambda_1+\gamma_1 P^c+\gamma_2 Q^c)\\
	&+P^c d(\lambda_2+\gamma_2 P^c-\gamma_1 Q^c).
	\end{aligned}
	\]
	We can write this as
	\begin{equation}
	\label{eqn:1-form-beta}
	\begin{aligned}
	\beta&=-Q^cd\lambda_1+P^cd\lambda_2 -\gamma_1 (P^cdQ^c+Q^cdP^c)+\gamma_2 (P^cdP^c-Q^cdQ^c)\\
	&=-Q^cd\lambda_1+P^cd\lambda_2+d\bigg[ -\gamma_1 P^c Q^c + \frac{1}{2}\gamma_2((P^c)^2-(Q^c)^2) \bigg].
	\end{aligned}
	\end{equation}
	Hence, by Green's formula and the definition of $1$-forms $\alpha$ and $\beta$, we have 
	\begin{align*}
	&\mu_{{x}+g_{t,\gamma},\varepsilon}(\Phi^{(\varepsilon)}_{t,\gamma}(\mathrm{\Gamma}))\\
	&=\frac{1}{4\pi}\int_{\partial \Phi^{(\varepsilon)}_{t,\gamma}(\mathrm{\Gamma})} -\frac{\partial S({x}+g_{t,\gamma},z,\varepsilon)}{\partial z_2}dz_1+\frac{\partial S({x}+g_{t,\gamma},z,\varepsilon)}{\partial z_1}dz_2\\
	&=\frac{1}{2\pi}\int_{\partial \Phi^{(\varepsilon)}_{t,\gamma}(\mathrm{\Gamma})} \alpha 
	=\frac{1}{2\pi} \int_{\partial \mathrm{\Gamma}} \beta\\
	&=\frac{1}{2\pi} 
	\int_{\partial \mathrm{\Gamma}}-Q^c(\lambda_1, \lambda_2)  d\lambda_1+ P^c(\lambda_1, \lambda_2)  d\lambda_2\\
	&=\frac{1}{4\pi}\int_{\partial \mathrm{\Gamma}}-\frac{\partial S({x}+c_t,\lambda,\varepsilon)}{\partial\lambda_2}d\lambda_1+\frac{\partial S({x}+c_t,\lambda,\varepsilon)}{\partial\lambda_1}d\lambda_2\\
	&=\mu_{{x}+c_t,\varepsilon}(\mathrm{\Gamma}),
	\end{align*}
	where we used \eqref{eqn:1-form-beta} to deduce the fourth identity. 
\end{proof}

\subsection{Addition with an elliptic operator}
In this section, we show the Brown measure of $\mu_{{x}+g_{t,\gamma}}$ is the push-forward measure of $\mu_{{x}+c_t}$ under the map $\Phi_{t,\gamma}$.
 Hence, the following diagram commutes. 

\begin{center}
\begin{tikzpicture}
\matrix (m) [matrix of math nodes,row sep=3em,column sep=4em,minimum width=2em]
{
	\mu_{{x}+c_t,\varepsilon} & \mu_{{x}+g_{t,\gamma},\varepsilon} \\
	\mu_{{x}+c_t} & \mu_{{x}+g_{t,\gamma}} \\};
\path[-stealth]
(m-1-1) edge node [left] {$\varepsilon\rightarrow 0$} (m-2-1)
edge [double] node [below] {$\Phi_{t,\gamma}^{(\varepsilon)}$} (m-1-2)
(m-2-1.east|-m-2-2) edge [double] node [below] {{$\Phi_{t,\gamma}$}}
(m-2-2)
(m-1-2) edge node [right] {$\varepsilon\rightarrow 0$} (m-2-2)
(m-2-1);
\end{tikzpicture}
\end{center}

\begin{theorem}
\label{thm:main-push-forward-general}
The Brown measure of $\mu_{{x}+g_{t,\gamma}}$ is the push-forward measure of the Brown measure $\mu_{{x}+c_t}$ under the map $\Phi_{t,\gamma}$. 
\end{theorem}
\begin{proof}
It is known in Theorem \ref{thm:main-push-forward-general-regularized} that 
\[
   \int_\mathbb{C} F(u) d\mu_{{x}+g_{t,\gamma},\varepsilon}(u)=\int_\mathbb{C} F\circ \Phi^{(\varepsilon)}_{t,\gamma}(u)
     d\mu_{{x}+c_t,\varepsilon}(u)
\]
for any $F\in C_c^\infty(\mathbb{C})$. The regularized Brown measure converges to the Brwon measure weakly as $\varepsilon$ tends to zero. In addition, $\Phi^{(\varepsilon)}_{t,\gamma}$ converges to $\Phi_{t,\gamma}$ uniformly by Lemma \ref{lemma:sub-w-uniform-convergence}. It follows that,
for any $F\in C_c^\infty(\mathbb{C})$, 
\[
   \int_\mathbb{C} F(u) d\mu_{{x}+g_{t,\gamma}}(u)=\int_\mathbb{C} F\circ \Phi_{t,\gamma}(u)
     d\mu_{{x}+c_t}(u).
\]
Hence, $\mu_{{x}+g_{t,\gamma}}$ is the push-forward measure of $\mu_{{x}+c_t}$ under the map $\Phi_{t,\gamma}$. 
\end{proof}

\subsection{Some further properties of the push-forward map} In this section, we study the special case when the map $\Phi_{t,\gamma}$ is nonsingular.

\begin{lemma}
\label{lemma:FK-det-derivatives-connection-ct-gt}
Given $t>0$, $\gamma\in\mathbb{C}$ such that $|\gamma|\leq t$, and $\lambda\in\Xi_t$, assume that the {Jacobian} of $\Phi_{t,\gamma}$ is invertible at $\lambda$, then  the function $z\mapsto S({x}+g_{t,\gamma},z,0)$ is a real analytic function of $z$ in a neighborhood of $\Phi_{t,\gamma}(\lambda)$. Moreover, we have
\begin{equation}
\label{eqn:p-z-lambda-ct-vs-x}
  p_z^{g,(t,\gamma)}(0)=p_\lambda^{c,(t)}(0), \qquad  p_{\overline{z}}^{g,(t,\gamma)}(0)=p_{\overline\lambda}^{c,(t)}(0),
\end{equation}
where $z=\Phi_{t,\gamma}(\lambda)$. 

In particular, if the map $\lambda\mapsto \Phi_{t,\gamma} (\lambda)$ 
is non-singular at all $\lambda\in\Xi_t$, the functions $z\mapsto S({x}+g_{t,\gamma},z,0)$ is a real analytic function of $z$ for all $z\in\Phi_{t,\gamma}(\Xi_t)$, and the identities \eqref{eqn:p-z-lambda-ct-vs-x} hold for any $\lambda\in \Xi_t$. 
\end{lemma}
\begin{proof}
By Theorem \ref{thm:main-FK-det-x-t-gamma-0}, under the assumption that the Jacobian of $\Phi_{t, \gamma}$ at $\lambda$ is invertible, we know that the map $(z,\varepsilon)\mapsto\Delta\big(|{x}+g_{t,\gamma}-z\unit|^2+\varepsilon^2\big)$ has a real analytic extension in a neighborhood of $(\Phi_{t,\gamma} (\lambda),0)$. It follows that, we may take the limit as $\varepsilon$ goes to zero in the following identity from Corollary \ref{thm:sub-comparison-Ct}
\[
   p_{z(\varepsilon)}^{g,(t,\gamma)}(\varepsilon)=p_\lambda^{c,(t)}(\varepsilon),
\]
where $z(\varepsilon)=\Phi_{t, \gamma}^{(\varepsilon)}(\lambda)$. This yields, at $z=\Phi_{t, \gamma}(\lambda)$ we have 
\[
p_{z}^{g, (t,\gamma)}(0) -p_\lambda^{c,(t)}(0)=0.
\]
Similarly, we have, 
\[
p_{\overline{z}}^{g, (t,\gamma)}(0) -p_{\overline{\lambda}}^{c,(t)}(0)=0
\]
as well. The the above argument works for any $\lambda\in\Xi_t$ if the map $\Phi_{t,\gamma}$ is non-singular at any $\lambda\in\Xi_t$. This concludes the statement.
\end{proof}

\begin{remark}
It is interesting to compare Lemma \ref{lemma:FK-det-derivatives-connection-ct-gt} with Lemma \ref{lemma:partial-lambda-c-t-limit}. Choose $z=\lambda+\gamma p_\lambda^{(0)}(w(\varepsilon;\lambda,t))$. If the map $\Phi_{t,\gamma}$ associated with ${x}$ is singular,  it turns out $\lim_{\varepsilon\rightarrow 0^+} p_z^{g, (t,\gamma)}(\varepsilon)$ has a limit that depends on $\lambda$. Hence, the condition in Lemma \ref{lemma:FK-det-derivatives-connection-ct-gt} is necessary. See Example \ref{example-Phi-semicircular}.
\end{remark}

\begin{proposition}
\label{prop:non-singular-imply-1-1}
If the map $\lambda\mapsto \Phi_{t,\gamma} (\lambda)$ 
is non-singular at all $\lambda\in\Xi_t$, then the map $\lambda\mapsto\Phi_{t, \gamma}(\lambda)$ is one-to-one in $\Xi_t$. 
\end{proposition}
\begin{proof}
Assume that $z=\Phi_{t, \gamma}(\lambda_1)=\Phi_{t, \gamma}(\lambda_2)$ for $\lambda_1, \lambda_2\in\Xi_t$. Then, by Lemma \ref{lemma:FK-det-derivatives-connection-ct-gt}, we have, for $i=1,2$, 
\[
  p_z^{g,(t,\gamma)}(0)=p_{\lambda_i}^{c,(t)}(0)=p_{\lambda_i}^{(0)}(w(0;\lambda_i,t)).
\]
Rewrite the map $\Phi_{t, \gamma}$ as 
\[
  \Phi_{t, \gamma}(\lambda)=\lambda+p_\lambda^{(0)}(w(0;\lambda,t))=\lambda+ p_z^{g,(t,\gamma)}(0), 
\]
where $z=\Phi_{t, \gamma}(\lambda)$. 
Then the condition $z=\Phi_{t, \gamma}(\lambda_1)=\Phi_{t, \gamma}(\lambda_2)$ yields that
$\lambda_1=\lambda_2$. Hence, $\Phi_{t, \gamma}$ is one-to-one in $\Xi_t$. 
\end{proof}

\subsection{An intrinsic formulation of the push-forward map}
\label{section:Phit_new_formula}
We establish another formulation of the push-forward map $\Phi_{t,\gamma}$ without invovling the subordination function $w(\varepsilon;\lambda,t)$ and the boundary value $w(0;\lambda,t)$.
This is possible due to Theorem \ref{thm:main-theorem-from-BYZ2021} obtained in our joint work with Belinschi and Yin \cite{BelinschiYinZhong2021Brown}. We put some technical results in Appendix \ref{appendix-A}. 
The Cauchy transform of a probability measure $\mu$ on $\mathbb{C}$ is a function defined on $\mathbb{C}$ as following
\[
  G_\mu(\lambda)=\lim_{\varepsilon\downarrow 0}\int_{\vert \lambda-z\vert \geq \varepsilon}\frac{1}{\lambda-z}d\mu(z).
\] 
The integral above is finite at every $\lambda$ if there exists $s>1$ such that $\mu\left( \{ z: \vert z-\lambda\vert \leq r \} \right)\leq r^s$ for all $r>0$. 

The following is the main result in this Section. 
\begin{theorem}
 \label{thm:Phit_new}
 The map $\Phi_{t,\gamma}$ definied in \eqref{eqn:def-Phi-t-gamma} can be rewritten as
  \[
    \Phi_{t,\gamma}(\lambda)=\lambda+\gamma \cdot G_{\mu_{x+c_t}}(\lambda), \qquad 
      \lambda\in \mathbb{C}.
  \]
 For $\lambda$ not on the boundary set $\partial(\Xi_t)$ of $\Xi_t$, the map $\Phi_{t,\gamma}$ is rewritten as
 \[
    \Phi_{t,\gamma}(\lambda)=\lambda+\gamma\cdot \frac{\partial S}{\partial\lambda}(x+c_t,\lambda,0), \qquad \lambda\notin \partial(\Xi_t).
 \]
\end{theorem}

Similar to the approach in Section \ref{section:regularized-pushforward-map}, we first establish an analogous result for the regularized map $\Phi_{t,\gamma}^{(\varepsilon)}$. 
\begin{lemma}
 \label{lemma_new_regulariedPhi}
  The map $\Phi_{t,\gamma}^{(\varepsilon)}$ can be rewritten as
  \[
    \Phi_{t,\gamma}^{(\varepsilon)}(\lambda)=\lambda+\gamma\cdot \frac{\partial S}{\partial\lambda}(x+c_t,\lambda,\varepsilon)=\lambda+\gamma \cdot G_{\mu_{x+c_t,\varepsilon}}(\lambda)
  \]  
for any $\lambda\in\mathbb{C}$.
\end{lemma}
\begin{proof}
 By the definition of $\Phi_{t,\gamma}^{(\varepsilon)}$ given by \eqref{defn:Phi-t-gamma-epsitlon} and the subordination relation, we have
 \begin{align*}
p_\lambda^{(0)}( w(\varepsilon;\lambda,t)  ) &=
\phi\bigg[ (\lambda-{x})^*\big( (\lambda-{x})(\lambda-{x})^*+w(\varepsilon; \lambda,t)^2 \unit\big)^{-1}\bigg]\\
&=\phi\bigg[ (\lambda-{x}-c_t)^*\big( (\lambda-{x}-c_t)(\lambda-{x}-c_t)^*+\varepsilon^2 \unit\big)^{-1}\bigg]\\
&=\frac{\partial S}{\partial\lambda}(x+c_t,\lambda,\varepsilon).
 \end{align*}
 Note that the regularized Brown measure $\mu_{x+c_t,\varepsilon}$ is absolutely continuous and  its density is bounded on the complex plane thanks to Theorem \ref{thm:main-theorem-from-BYZ2021} (or \cite[Lemma 2.8]{HaagerupSchultz2007}). Then the second identity follows from Theorem \ref{thm:appendix-derivative-identity} in Appendix A. 
\end{proof}

\begin{lemma}\label{lemma_continuityCauchy}
The Cauchy transform of $\mu_{x+c_t}$ is continuous on $\mathbb{C}$. Moreover, 
\begin{equation}
 \label{eqn:5.7.limit}
 \lim_{\varepsilon\rightarrow 0}G_{\mu_{x+c_t,\varepsilon}}(\lambda)=G_{\mu_{x+c_t}}(\lambda)
\end{equation}
for any $\lambda\in\mathbb{C}$. 
\end{lemma}
\begin{proof}
The first claim follows from Lemma \ref{lemma:appendix_continuityCauchy} as $\mu_{x+c_t}$ satisfies the condition thanks to Theorem \ref{thm:main-theorem-from-BYZ2021}.

Note that the function $z\mapsto\frac{1}{\lambda-z}$ is locally integrable. Denote by 
$B_\rho(\lambda)$ the disk $B_\rho(\lambda)=\{w: |\lambda-w|<\rho\}$. Given $\rho<\text{min}\left\{\frac{1}{4}, \frac{1}{4|\lambda|}, \frac{1}{4\Vert x+c_t\Vert} \right\}$, we choose a continuous function $\chi$ defined on the complex plane with compact support such that:
(1) $\chi(z)\in [0,1]$ for all $z\in\mathbb{C}$;
(2) $\text{supp}(\chi)\subset B_{2/\rho}(\lambda)$; (3) $\chi(z)=1$ for $z\in B_{1/\rho}(\lambda)\backslash B_{2\rho}(\lambda)$;  (4) $\chi(z)=0$ for $z\in B_\rho(\lambda)$. Then
\[
  \lim_{\varepsilon\rightarrow 0}\int_{\mathbb{C}}
    \frac{\chi(z)}{\lambda-z}d\mu_{x+c_t,\varepsilon}(z)=
     \int_{\mathbb{C}}
    \frac{\chi(z)}{\lambda-z}d\mu_{x+c_t}(z).
\]
By Theorem \ref{thm:main-theorem-from-BYZ2021}, for $\mu$ being either $\mu_{x+c_t,\varepsilon}$ or $\mu_{x+c_t}$, we have 
\begin{align*}
  \left\vert \int_{B_{2\rho}(\lambda)}\frac{1-\chi(z)}{\lambda-z} d\mu(z) \right\vert 
   &\leq \frac{1}{\pi{t}} \int_{B_{2\rho}(\lambda)}\left\vert \frac{1}{\lambda-z}  \right\vert dm_2(z)\\
  & \leq \frac{1}{\pi{t}} \int_{0}^{2\pi}\int_{0}^{2\rho}drd\theta
   =\frac{4\rho}{t},
\end{align*}
where $m_2$ denotes Lebesgue measure on the complex plane. 

For $z$ with $\vert z\vert \geq \frac{1}{2\rho}$, then $z\notin \text{spec}(x+c_t)$, hence $z$ is not in the support of $\mu_{x+c_t}$ or equivalently $z\notin \overline{\Xi_t}$. Moreover, the density of $\mu_{x+c_t,\varepsilon}$ (see \cite[Lemma 2.8]{HaagerupSchultz2007}) can be estimated as
\begin{align*}
 & \frac{1}{4\pi}\varepsilon^2\phi \big ( (|z-x-c_t|^2+\varepsilon^2)^{-1} (|(z-x-c_t)^*|^2+\varepsilon)^{-1} \big) \\
   &\leq \frac{1}{4\pi}\phi \big ( (|z-x-c_t|^2+\varepsilon^2)^{-1} \big)\\
   & \leq \frac{1}{4\pi} \cdot \left(\Vert z-x-c_t \Vert\right)^{-2}\\
   &\leq \frac{1}{4\pi} \frac{1}{\left( \vert {z}\vert -\Vert x+c_t\Vert \right)^2}
   \leq \frac{1}{\pi \vert {z}\vert^2}
\end{align*}
where we used the assumption $\Vert x+c_t\Vert <\frac{1}{4\rho}\leq \frac{\vert{z}\vert}{2}$.
By the choice of $\rho$, if $z\in  \mathbb{C}\backslash B_{1/\rho}(\lambda)$, then
$\vert z\vert \geq \frac{1}{2\rho}$. Hence, there exists some $C>0$ such that
\begin{align*}
 \left | \int_{\mathbb{C}\backslash B_{1/\rho}(\lambda)}\frac{1-\chi(z)}{\lambda-z} d \mu_{x+c_t,\varepsilon} (z)  \right | \leq C\rho,
\end{align*}
for any $\lambda$ fixed. 

Putting the above results together, we obtain 
\[
  \lim_{\varepsilon\rightarrow {0}}\int_\mathbb{C}\frac{1}{\lambda-z}d\mu_{x+c_t,\varepsilon}(z)
   =\int_\mathbb{C}\frac{1}{\lambda-z}d\mu_{x+c_t}(z)
\]
for any $\lambda\in\mathbb{C}$. 
\end{proof}

We then study the derivative of $S(x+c_t,\lambda,0)$ outside $\overline{\Xi_t}$. 
\begin{proposition}\label{prop:5.13}
For any $\lambda\notin \overline{\Xi_t}$, we have 
\[
 \frac{\partial S}{\partial\lambda}(x+c_t,\lambda,0)=
  G_{\mu_{x+c_t}}(\lambda)=G_{\mu_x}(\lambda).
\]
\end{proposition}
\begin{proof}
By Theorem \ref{thm:main-FK-det-ct-0}, for $\lambda\notin \overline{\Xi_t}$ we have 
\begin{align*}
  \log\Delta(x+c_t-\lambda)&=\log\Delta(x-\lambda)\\
   &=\int_\mathbb{C}\log|\lambda-z|d\mu_x(z)
  =\int_\mathbb{C}\log|\lambda-z|d\mu_{x+c_t}(z).
\end{align*}
Taking the partial derivatives imply the desired identities. 
\end{proof}

We next study the derivative of $S(x+c_t,\lambda,0)$ in $\Xi_t$. 
\begin{proposition}\label{prop:5.14}
Given $t>0$ and $\lambda\in\Xi_t$, we have 
  \[
      \frac{\partial S}{\partial\lambda}(x+c_t,\lambda,0)=p_\lambda^{(0)}(w(0;\lambda,t))
      =G_{\mu_{x+c_t}}(\lambda).
  \]
\end{proposition}
\begin{proof}
By Theorem \ref{thm:main-FK-det-ct-0} and Notation \ref{notation:partial-derivatives-lambda-c-x}, for $\lambda\in\Xi_t$, we have 
\begin{align*}
  S(x+c_t,\lambda,0)&=\log\Delta(x+c_t-\lambda)^2\\
     &=\log\Delta(|x-\lambda|^2+w^2)-\frac{w^2}{t},
\end{align*}
where $w=w(0;\lambda,t)$ is defined in Definition \ref{defn:s-epsilon-t-00} by
\[
 	\int_0^\infty \frac{1}{w^2+u^2}d\mu_{|{x}-\lambda\unit|}(u)=\frac{1}{t}.
\]
We have 
\begin{align*}
 \frac{\partial\log\Delta(|x-\lambda|^2+w^2)}{\partial\lambda}
   &=\phi\large( (\lambda-x)^*(|\lambda-x|^2+w^2)^{-1} \large)\\
     &\qquad +2w\cdot \frac{\partial w}{\partial\lambda}\cdot\phi\large((|\lambda-x|^2+w^2)^{-1}\large)\\
    &= \phi\large( (\lambda-x)^*(|\lambda-x|^2+w^2)^{-1} \large)
     +\frac{2w}{t}\cdot \frac{\partial w}{\partial\lambda}.
\end{align*}
Cancelling the second term $\frac{\partial}{\partial\lambda}\left(\frac{w^2}{t}\right)$ yields 
\begin{align*}
  \frac{\partial S}{\partial\lambda}(x+c_t,\lambda,0)
  =\phi\large( (\lambda-x)^*(|\lambda-x|^2+w^2)^{-1} \large)
  =p_\lambda^{(0)}(w).
\end{align*}
This established the first identity. The second identity follows from Theorem \ref{thm:appendix-derivative-identity} in Appendix A and the fact that $\mu_{x+c_t}$ is absolutely continuous with respect to Lebesgue measure and its density function $h$ is bounded by the constant $1/\pi{t}$. 
\end{proof}

\begin{proof}[Proof of Theorem \ref{thm:Phit_new}]
By Lemma \ref{lemma_new_regulariedPhi} and Lemma \ref{lemma_continuityCauchy}, we have 
\begin{align*}
  \Phi_{t,\gamma}(\lambda)=\lim_{\varepsilon\rightarrow 0}
   \Phi_{t,\gamma}^{(\varepsilon)}(\lambda)
    =\lambda+\gamma\cdot G_{\mu_{x+c_t}}(\lambda),
\end{align*}
for any $\lambda\in\mathbb{C}$. 
The second formulation of $\Phi_{t,\gamma}$ is a combination of statements in Proposition \ref{prop:5.13} and Proposition \ref{prop:5.14}.
\end{proof}

\section{Addition with a self-adjoint operator}
\label{section:addition-selfadj}

The special case when ${x}$ is self-adjoint has drawn much attention in previous work \cite{HoHall2020Brown, Ho2020Brown, HoZhong2020Brown}. 
In this section, we apply our main result to generalize main results in those work. The generalization can be viewed as the addition analogue of recent work \cite{HallHo2021Brown} concerning free multiplicative Brownian motions.

\subsection{Subordination functions}
\label{section-sub-x0-self-adjoint}
Let ${x}$ be self-adjoint and $\mu=\mu_{{x}}$ be its spectral measure. We first study the subordination function as in Definition \ref{defn:s-epsilon-t-00} and Proposition \ref{prop:sub-identical-symmetric-measure}. 
The set $\Xi_t$ is expressed as
\begin{equation}
\label{eqn:Xi-t-self-adjoint}
  \Xi_t=\left\{\lambda=a+{\i}b: \int_\mathbb{R}\frac{1}{(u-a)^2+b^2}d\mu(u)>\frac{1}{t} \right\}.
\end{equation}
For $\lambda\in \Xi_t$ (which is equivalent to $t>\lambda_1(\mu)^2$), 
the condition subordination function $w(0;\lambda,t)$ is determined by the condition \eqref{eqn:condition-s-0-t}, that can be rewritten as 
\begin{equation}
 \label{eqn:condition-w-0-x0-self-adjoint}
    \int_\mathbb{R}\frac{1}{(u-a)^2+b^2+w(0; \lambda,t)^2} d\mu(u)=\frac{1}{t}.
\end{equation}
Following Biane's work \cite[Section 3]{Biane1997} on the spectral measure of ${x}+g_t$,  we set 
\[
   U_t=\left\{ a\in\mathbb{R}: \int_\mathbb{R}\frac{1}{(u-a)^2}d\mu(u)>\frac{1}{t} \right\},
\]
and define $v_t$ as follows
\begin{equation}
\label{defn:v-t}
  v_t(a) = \inf\left\{y>0: \int_{\mathbb{R}} \frac{1}{(a-x)^2+y^2}d\mu(x) \leq \frac{1}{t}\right\}, \quad a\in\R.
\end{equation}
We then set 
\[
  \Omega_t=\{ a+{\i}b: |b|>v_t(a) \}.
\]
It follows that $v_t(a)^2=b^2+w(0; \lambda,t)^2$ with $\lambda=a+{\i}b$ if $\lambda\in \Xi_t$; and $v_t(a)=0$ if $\{a+{\i}b: b\in\mathbb{R}\} \cap \Xi_t=\emptyset$. By regularity of $w(0;\lambda,t)$ in Lemma \ref{lemma:s-epsilon-t-limit}, we see that $v_t$ is a continuous function (see \cite{Biane1997} for the original definition and \cite[Section 2.3]{HoZhong2020Brown} for a review).

\begin{figure}[t]
	\label{figure1-circular}
	\begin{center}
		\begin{subfigure}[h]{0.75\linewidth}
			\includegraphics[width=\linewidth]{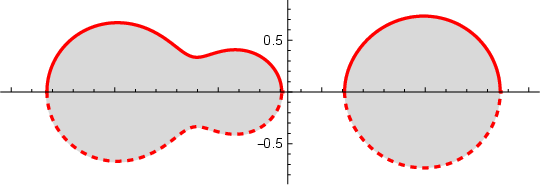}
			\caption{The domain $\Xi_t$ for $t=1$ and ${x}$ distributed as $0.4\delta_{-2}+0.1\delta_{-0.8}+0.5\delta_1$. The graph of $v_t$ is the solid curve above the $x$-axis.}
		\end{subfigure}
		\begin{subfigure}[h]{0.75\linewidth}
			\includegraphics[width=\linewidth]{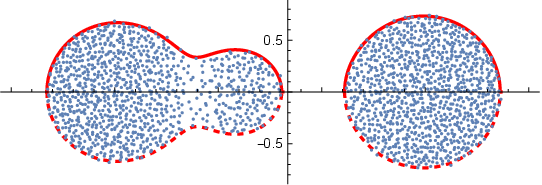}
			\caption{The $2000\times 2000$ random matrix simulation for the Brown measure of ${x}+c$.}
		\end{subfigure}	
	\end{center}
	\begin{center}
		\caption{The domain $\Xi_t$, the graph of $v_t$, and Brown measure simulation.}
	\end{center}
\end{figure}

\begin{proposition}
\label{prop:Xi-t-characterization-self-adjoint}
The subordination function $w(0;\lambda,t)$ as in Definition \ref{defn:s-epsilon-t-00} can be expressed as
\begin{align*}
    w(0;\lambda,t)=\begin{cases}
     \sqrt{v_t(a)^2-b^2}, &\qquad \text{for}\quad \lambda=a+{\i}b\in X_t;\\
     0, &\qquad \text{for}\quad \lambda=a+{\i}b\notin \Xi_t.
    \end{cases}
\end{align*}
In particular, $a\mapsto v_t(a)$ is a continuous function defined on $\mathbb{R}$, and is real analytic and $v_t(a)>0$ if $\{ a+{\i}b: b\in\mathbb{R} \}\cap \Xi_t \neq \emptyset$.

Moreover, we have $\Xi_t\cap \mathbb{R}=U_t$, and 
\[
  \Xi_t=\{ a+{\i}b\in\mathbb{C}: |b|<v_t(a) \}
\]
and
\[
   \overline{\Omega_t}=\{ a+{\i}b: |b|\geq v_t(a) \}=\big(\Xi_t \big)^c.
\]
\end{proposition}

\begin{definition}
For $t>0$ and $v_t$ as defined in \eqref{defn:v-t}, we set 
\begin{equation}
\label{defn:psi-t-a}
  \psi_t(a)=a+t\int_\mathbb{R}\frac{a-u}{(a-u)^2+v_t(a)^2}d\mu(u), \qquad a\in\mathbb{R},
\end{equation}
and 
\begin{equation}
 \label{defn:h-t-a}
 h_t(a)=t\int_\mathbb{R}\frac{a-u}{(a-u)^2+v_t(a)^2}d\mu(u), \qquad a\in\mathbb{R}.
\end{equation}
\end{definition}

The following result play a key role in the following discussion. It is taken from \cite[Lemma 5]{Biane1997} and \cite[Theorem 3.14]{HoZhong2020Brown}. 

\begin{proposition}
\label{prop:psi-t-derivative}
The function $a\mapsto \psi_t(a)$ is a homeomorphism of $\mathbb{R}$ onto $\mathbb{R}$. Moreover, if $v_t(a)>0$ or equivalently $\{a+{\i}b: b\in\mathbb{R}\}\cap \Xi_t\neq \emptyset$, then 
\begin{equation}\nonumber
   0<\frac{d\psi_t(a)}{da}<2.
\end{equation}
Consequently, if $v_t(a)>0$, then
\begin{equation} 
\label{eqn:h-t-derivative}
  -1<\frac{d\, h_t(a)}{da}<1.
\end{equation}
\end{proposition}

Finally, we need the following useful result taken from \cite{Biane1997}. 
\begin{lemma}\cite[Lemma 3]{Biane1997} and \cite[Proposition 3]{Biane1997}
\label{lemma:Cauchy-x0-sefladjoint-in-Omega-t}
The Cauchy transform of $G_\mu$ has a continuous extension to $\overline{\Omega_t}$ which is Lipschitz with Lipschitz constant $\leq 1/t$, and one has 
\[
   |G_\mu(z)|^2\leq \int_\mathbb{R}\frac{1}{|z-u|^2}d\mu(u)\leq \frac{1}{t}
\]
for $z\in \overline{\Omega_t}$. Moreover, the support of $\mu$ is contained in the closure of $U_t$. Consequently, $\supp(\mu)\subset \Xi_t$.
\end{lemma}

\subsection{The push-forward map}

Let ${x}$ be a self-adjoint operator and $\mu=\mu_{{x}}$ be its spectral measure. We study the map $\Phi_{t,\gamma}$ as defined in Section \ref{section:Brown-free-additiveBM}. 
We note that when $\gamma=t$, the operator $g_{t,\gamma}$ is a semicircular operator  $g_t$ with mean zero and variance $t$ and the Brown measure of ${x}+g_t$ is just the spectral measure of ${x}+g_t$. Hence, we will assume that 
$\gamma\neq t$ throughout this section. Let $\tau=t-\gamma$ and $\tau=\tau_1+i\tau_2$. Since $|\gamma|\leq t$, we see that $\tau_1> 0$ if $\gamma\neq t$. 

Using notations in Subsection \ref{section-sub-x0-self-adjoint}, we have 
\[
   \Phi_{t,\gamma}(\lambda)=\lambda+\gamma \int_\mathbb{R}\frac{\overline{\lambda}-u}{(u-a)^2+b^2+w(0;\lambda,t)^2}d\mu(u),\qquad \lambda\in\mathbb{C},
\]
and, if $\lambda\in\Xi_t$ then $\Phi_{t,\gamma}(\lambda)$ can be rewritten as 
\[
  \Phi_{t,\gamma}(\lambda)=\lambda+\gamma \int_\mathbb{R}\frac{\overline{\lambda}-u}{(u-a)^2+v_t(a)^2}d\mu(u), \qquad   \lambda\in \Xi_t.
\]

\begin{proposition}
\label{prop:Phi-t-gamma-x0-self-adjoint-formula}
For any $|\gamma|\leq t$ and $\gamma\neq t$, we have 
\begin{align*}
  \Phi_{t,\gamma}(\lambda)=
    \begin{cases}
      \psi_t(a)-\frac{\tau}{t}h_t(a)+i\frac{\tau b}{t}, \qquad & \text{for}\qquad 
          \lambda\in\Xi_t;\\
         \lambda+\gamma G_\mu(\lambda)  \qquad &\text{for}\qquad 
         \lambda\notin \Xi_t.
    \end{cases}
\end{align*}
\end{proposition}
\begin{proof}
We rewrite $\Phi_{t,\gamma}(\lambda)$ as 
\begin{align*}
     \Phi_{t,\gamma}(\lambda)&=\Phi_{t,t}(\lambda)-\tau  \int_\mathbb{R}\frac{\overline{\lambda}-u}{(u-a)^2+b^2+w(0;\lambda,t)^2}d\mu(u);
\end{align*}
and if $\lambda\in\Xi_t$, we can further rewrite it as 
\begin{align*}
       \Phi_{t,\gamma}(\lambda)  &=\Phi_{t,t}(\lambda)-\tau\int_\mathbb{R}\frac{\overline{\lambda}-u}{(u-a)^2+v_t(a)^2}d\mu(u)\\
       &=\Phi_{t,t}(\lambda)-\frac{\tau}{t}h_t(a)+i\frac{\tau b}{t},\qquad \qquad\qquad \lambda\in\Xi_t,
\end{align*}
where we used \eqref{eqn:condition-w-0-x0-self-adjoint} and the definition of $h_t$ in \eqref{defn:h-t-a}.
If $\lambda\in\Xi_t$, we have 
\begin{align*}
 \Phi_{t,t}(\lambda)&= \lambda+t\int_\mathbb{R}\frac{\overline{\lambda}-u}{(u-a)^2+b^2+w(0;\lambda,t)^2}d\mu(u)\\
 &=\psi_t(a)+i\left( b-t \int_\mathbb{R}\frac{b}{(u-a)^2+b^2+w(0;\lambda,t)^2}d\mu(u)\right)=\psi_t(a)
\end{align*}
where we used \eqref{eqn:condition-w-0-x0-self-adjoint} and the definition of $\psi_t$ in \eqref{defn:psi-t-a}. This proves the case when $\lambda\in\Xi_t$.

Recall that $\big(\Xi_t \big)^c= \overline{\Omega_t}$ and $\supp(\mu)\subset U_t\subset\Xi_t$ (recall Lemma \ref{lemma:Cauchy-x0-sefladjoint-in-Omega-t}). When $\lambda\notin \Xi_t$, then $w(0;\lambda,t)=0$. Hence, when $\lambda=a+{\i}b\notin \Xi_t$, 
\begin{align*}
  \Phi_{t,\gamma}(\lambda)&=\lambda+\gamma\int_\mathbb{R}\frac{\overline{\lambda}-u}{(u-a)^2+b^2}d\mu(u)\\
  &=\lambda+\gamma\int_\mathbb{R}\frac{1}{\lambda-u}d\mu(u)=\lambda+\gamma G_\mu(\lambda).
\end{align*}
This finishes the proof. 
\end{proof}

The following result says that $\Phi_{t,\gamma}$ is injective on $\big(\Xi_t \big)^c$. It can be viewed as a generalization of \cite[Lemma 4]{Biane1997} where the latter studies the inverse map of subordination function for the free additive convolution of ${x}+g_t$ (see also \cite[Section 2.4]{HoZhong2020Brown} and \cite[Proposition 5.5]{HallHo2021Brown}). We include a proof for convenience. 

\begin{proposition}
\label{prop:Phi-t-gamma-self-adjoint-outside}
The map $\Phi_{t,\gamma}(\lambda)$ is an injective map when restricted to $\big(\Xi_t \big)^c$ and is conformal when $\lambda$ is in the interior of $\big(\Xi_t \big)^c$.
\end{proposition}
\begin{proof}
If $\lambda\in \big(\Xi_t \big)^c$, by Proposition \ref{prop:Phi-t-gamma-x0-self-adjoint-formula}, we have $\Phi_{t,\gamma}(\lambda)=\lambda+\gamma G_\mu(\lambda)$. Hence, for $\alpha_1, \alpha_2 \in \big(\Xi_t \big)^c$, we have 
\begin{align*}
    \Phi_{t,\gamma}(\alpha_1)-\Phi_{t,\gamma}(\alpha_2)&=\alpha_1-\alpha_2+\gamma (G_\mu(\alpha_1)-G_\mu(\alpha_2))\\
    &=(\alpha_1-\alpha_2)\left( 1+\gamma \frac{G_\mu(\alpha_1)-G_\mu(\alpha_2)}{\alpha_1-\alpha_2} \right).
\end{align*}
Then, by Cauchy-Schwarz inequality, 
\begin{align*}
  \left| \frac{G_\mu(\alpha_1)-G_\mu(\alpha_2)}{\alpha_1-\alpha_2}\right|
    &=\left| \int_\mathbb{R} \frac{1}{(\alpha_1-u)(\alpha_2-u)} d\mu(u)\right|\\
    &\leq\left( \int_\mathbb{R}\frac{d\mu(u)}{|\alpha_1-u|^2}
      \int_\mathbb{R}\frac{d\mu(u)}{|\alpha_2-u|^2}\right)^{1/2}
     \leq \frac{1}{t}
\end{align*}
where we used Lemma \ref{lemma:Cauchy-x0-sefladjoint-in-Omega-t} in the final step. 
It can be shown that we can not have equality in the Cauchy-Schwarz inequality used above. Indeed, if equality holds, then
\[
  \left| \int_\mathbb{R} \frac{1}{(\alpha_1-u)(\alpha_2-u)} d\mu(u)\right|=\frac{1}{t}
\]
for some $\alpha_1, \alpha_2\in \big(\Xi_t \big)^c=\overline{\Omega_t}$. One can show that $\mu$ must be a Dirac measure (see \cite[Lemma 4]{Biane1997} for details). Therefore, 
\[
     \Phi_{t,\gamma}(\alpha_1)-\Phi_{t,\gamma}(\alpha_2)\neq 0
\]
for $\alpha_1\neq \alpha_2$ and $\alpha_1, \alpha_2 \in  \big(\Xi_t \big)^c$.
\end{proof}

\begin{lemma}
\label{lemma:simple-tau-t}
For $|\gamma|\leq t$ and $\gamma\neq t$, let $\tau=t-\gamma$ and $\tau=\tau_1+i\tau_2$, then we have $0< \frac{|\tau|^2}{t\tau_1}\leq 2 $.
\end{lemma}
\begin{proof}
The assumption implies directly that $\tau_1>0$. Write $\gamma=\gamma_1+i\gamma_2$, then $\tau_1=t-\gamma_1$ and $\tau_2=\gamma_2$. Hence
\begin{align*}
  |\tau|^2-2t\tau_1=(t-\gamma_1)^2+\gamma_2^2-2t(t-\gamma_1)
   =|\gamma|^2-t^2\leq 0.
\end{align*}
The result then follows.
\end{proof}

\begin{theorem}
\label{thm:det-Jacobian-Phi}
For any $|\gamma|\leq t$ and $\gamma\neq t$, 
the map $\Phi_{t,\gamma}$ is a smooth, injective map on $\Xi_t$. Moreover, the determinant of the Jacobian of $\Phi_{t,\gamma}$ is strictly positive at all $\lambda\in \Xi_t$ and can be expressed
\begin{equation}
\label{eqn:det-Jacobian-Phi-t-gamma-self-adjoint}
  \det(\emph{Jacobian}(\Phi_{t,\gamma}))(\lambda)=\frac{\tau_1}{t}\left[ 1+\left( 1-\frac{|\tau|^2}{t\tau_1} \right) \frac{d\,h_t(a)}{d a} \right], \qquad \lambda=a+{\i}b,
\end{equation}
where $\tau=\tau_1+i\tau_2$ and $h_t$ is defined in \ref{defn:h-t-a}.
\end{theorem}
\begin{proof}
By Proposition \ref{prop:Phi-t-gamma-x0-self-adjoint-formula}, for $\lambda\in \Xi_t$, 
we have
\begin{align*}
 \Phi_{t,\gamma}(\lambda)=  \psi_t(a)-\frac{\tau}{t}h_t(a)+i\frac{\tau b}{t}.
\end{align*}
Recall $\psi_t(a)=a+h_t(a)$. 
Hence, for $\Phi_{t,\gamma}(\lambda)=z_1+iz_2$, we deduce 
\begin{align}
\label{eqn:z-1-z-2}
 z_1=z_1(a,b)=a+\left(1-\frac{\tau_1}{t}\right) h_t(a)-\frac{\tau_2}{t}b, \quad
 z_2=z_2(a,b)=-\frac{\tau_2}{t}h_t(a)+\frac{\tau_1}{t}b.
\end{align}
We then have
\begin{align*}
 \text {Jacobian}(\Phi_{t,\gamma}(\lambda))=
        \begin{bmatrix}
          1+\left(1-\frac{\tau_1}{t}\right) h_t'(a) & -\frac{\tau_2}{t}\\
          -\frac{\tau_2}{t}h_t'(a) & \frac{\tau_1}{t}
       \end{bmatrix}.
\end{align*}
It follows that 
\[
   \det(\text{Jacobian}(\Phi_{t,\gamma}))(\lambda)=\frac{\tau_1}{t}\left[ 1+\left( 1-\frac{|\tau|^2}{t\tau_1} \right) \frac{d\,h_t(a)}{d a} \right].
\] 
By Proposition \ref{prop:psi-t-derivative}, $-1<h_t'(a)<1$, hence we have
\[
 -1< \left( 1-\frac{|\tau|^2}{t\tau_1} \right) \frac{d\,h_t(a)}{d a} <1
\]
thanks to Lemma \ref{lemma:simple-tau-t}. This implies that $\det(\text{Jacobian}(\Phi_{t,\gamma}))(\lambda)$ is independent of $\,b\,$ and $\det(\text{Jacobian}(\Phi_{t,\gamma}))(\lambda)>0$ for all $\lambda\in \Xi_t$.

The non-singular property of $\Phi_{t,\gamma}$ implies that $\Phi_{t,\gamma}$ is one-to-one in $\Xi_t$ by Proposition \ref{prop:non-singular-imply-1-1}. In this case, we can actually prove it directly. If $a_1+ib_1, a_2+ib_2\in \Xi_t$ and 
$\Phi_{t,\gamma}(a_1+ia_2)=\Phi_{t,\gamma}(b_1+ib_2)$. Then we have, by using \eqref{eqn:z-1-z-2},
\begin{align*}
  a_1-a_2+\left( 1-\frac{\tau_1}{t}\right) (h_t(a_1)-h_t(a_2)) &=-\frac{\tau_2}{t}(b_2-b_1)\\
  -\frac{\tau_2}{t}(h_t(a_1)-h_t(a_2))&=\frac{\tau_1}{t}(b_2-b_1),
\end{align*}
which yields (by canceling $b_2-b_1$), 
\[
  (a_1-a_2)+\left( 1-\frac{|\tau|^2}{t\tau_1}\right)(h_t(a_1)-h_t(a_2)).
\]
By the first part of the proof, we see that the function $a\mapsto \frac{\tau_1}{t}\left[ 1+\left( 1-\frac{|\tau|^2}{t\tau_1} \right) \frac{d\,h_t(a)}{d a} \right]$ is strictly positive. Consequently, $a_1=a_2$. We then deduce that $b_1=b_2$. 
The injectivity property is established. 
\end{proof}

\begin{corollary}
For any $|\gamma|\leq t$ and $\gamma\neq t$, 
the map $\Phi_{t,\gamma}$ is a homeomorphism of $\mathbb{C}$ to $\mathbb{C}$. 
\end{corollary}
\begin{proof}
We show that $\Phi_{t,\gamma} ( \Xi_t  )\cap \Phi_{t,\gamma}\left( (\Xi_t)^c \right)=\emptyset$. Recall that 
\[
\big(\Xi_t \big)^c=\{ a+{\i}b: |b|\geq v_t(a) \}=\overline{\Omega_t}.
\]
For $\lambda_0=a+{\i}b$ so that $b\geq v_t(a)$, if $\Phi_{t,\gamma}(\lambda_0)\in \Phi_{t,\gamma} ( \Xi_t  )$, consider the vertical half-line $\{ \lambda=a+{\i}d: d\geq b \}$  starting at $\lambda_0$. Since $\lim_{d\rightarrow \infty}\Phi_{t,\gamma}(a+{\i}d)=\infty$, we then see that there is another point $\lambda_1=a+{\i}b_1$ where $b_1>b$ so that $\Phi_{t,\gamma} (\lambda_1)\in \partial ( \Phi_{t,\gamma} ( \Xi_t  ) )$, the boundary of $\Phi_{t,\gamma} ( \Xi_t  )$.
On the other hand $\Phi_{t,\gamma}(a+{\i}v_t(b))\in \partial ( \Phi_{t,\gamma} ( \Xi_t  ) )$, and both $\lambda_1$ and $a+{\i}v_t(b)$ are in $(\Xi_t)^c$, this contradicts to Proposition \ref{prop:Phi-t-gamma-self-adjoint-outside}. Therefore, 
\[
\Phi_{t,\gamma} ( \Xi_t  )\cap \Phi_{t,\gamma}\left( (\Xi_t)^c \right)=\emptyset.
\]
Since $\Phi_{t,\gamma}$ is a continuous function on $\mathbb{C}$ and it maps a neighborhood of infinity to some neighborhood of infinity, 
one can then deduce that $\Phi_{t,\gamma}(\mathbb{C})=\mathbb{C}$ by some standard arguments in topology. We conclude that $\Phi_{t,\gamma}$ is a homeomorphism of $\mathbb{C}$ to $\mathbb{C}$ due to Proposition \ref{prop:Phi-t-gamma-self-adjoint-outside} and Theorem \ref{thm:det-Jacobian-Phi}.
\end{proof}

\subsection{Brown measure of addition with a self-adjoint operator}
\label{section:Brown-density-self-adjoint}

We apply results in Section \ref{section:Brown-free-additiveBM} to give a new proof for the Brown measure of ${x}+c_t$ and recover a result in our previous work with Ho \cite{HoZhong2020Brown}.
We then study the Brown measure of ${x}+g_{t,\gamma}$ that extends results by Hall and Ho \cite{HoHall2020Brown, Ho2020Brown} to all twisted elliptic operators. 
The result about ${x}+g_{t,\gamma}$ in this section is the additive counterparts of the results of Hall-Ho \cite{HallHo2021Brown}.

\begin{theorem}\cite[Theoorem 3.10]{HoZhong2020Brown}
\label{thm:density-x0-c-t-self-adjoint}
For $\lambda=a+{\i}b\in \Xi_t$, then the Brown measure of ${x}+c_t$ is absolutely continuous at $\lambda$ and the density at $\lambda$ is given by
\[
  d\mu_{{x}+c_t}(a+{\i}b)= \frac{1}{\pi t}\left( 1- \frac{t}{2}\frac{d}{da}\int_\mathbb{R}\frac{u}{(a-u)^2+v_t(a)^2}d\mu(u) \right) da\, db,
\]
which can also be expressed as
\begin{equation}
\label{density-Brown-x0-ct-derivative-psi}
    d\mu_{{x}+c_t}(a+{\i}b)=  \frac{1}{2\pi t} \frac{d\psi_t(a)}{da}da\, db.
\end{equation}
In particular, the density is constant along the vertical directions.
\end{theorem}
\begin{proof}
 One can deduce the density formula directly from the general formula in Theorem \ref{thm:BrownFormula-x0-ct-general}. We note that
 \begin{align*}
  \frac{\partial}{\partial \overline{\lambda}} \bigg( \phi \big( {x}^* ( ({x}-\lambda\unit) &({x}-\lambda\unit)^* +w(0; \lambda,t)^2 )^{-1} \big) \bigg) \\
   &=\frac{\partial}{\partial\overline{\lambda}} \bigg( \int_\mathbb{R} \frac{u}{(a-u)^2+b^2+w(0:\lambda,t)^2} d\mu(u)\bigg)\\
   &=\frac{1}{2}\frac{\partial}{\partial {a}} \bigg( \int_\mathbb{R} \frac{u}{(a-u)^2+v_t(a)^2}\mu(u) \bigg)
 \end{align*}
 where we used the fact that the integration is independent of $b$. 
 The first formula then follows from \eqref{density-x0-ct-formula-2}. Using the formula for $\psi_t$ in \eqref{defn:psi-t-a}, we have
\begin{align*}
	    \frac{d\psi_t(a)}{da}&=\frac{d}{da}\left( a+ t\int_\mathbb{R}\frac{a-u}{(a-u)^2+v_t(a)^2}d\mu(u)\right)\\
	&=2-t\frac{d}{da}\int_\mathbb{R}\frac{u}{(a-u)^2+v_t(a)^2}d\mu(u)
\end{align*}
 where we used \eqref{eqn:condition-w-0-x0-self-adjoint}. This establish the second formula. 
\end{proof}

\begin{proposition}\cite[Lemma 3.11 and Proposition 3.16]{HoZhong2020Brown} and \cite[Corollary 3 and Proposition 3]{Biane1997}
\label{prop:push-forward-semicircular}
The Brown measure $\mu_{{x}+c_t}$ concentrates in $\Xi_t$. In other words, $\mu_{{x}+c_t}\big(\Xi_t\big)=1$. The spectral measure of ${x}+g_t$ is the push-forward of the Brown measure $\mu_{{x}+c_t}$ under the map $\lambda \mapsto \Phi_{t,t}$, where 
\[
 \Phi_{t,t}(a+{\i}b)=\psi_t(a), \qquad\text{for}\qquad a\in U_t, |b|\leq v_t(a).
\]
Consequently, the spectral measure of ${x}+g_t$ is absolutely continuous and its density is given by
\[
p_t(\psi_t(a))=\frac{v_t(a)}{\pi t}.
\]
Moreover, the support of $\mu$ is contained in the closure of $U_t$.
\end{proposition}

\begin{proof}
Recall that $\Xi_t=\{a+{\i}b\in\mathbb{C}: a\in U_t, |b|<v_t(a) \}$, and $v_t>0$ if $a\in U_t$, and $v_t=0$ if $a\in\mathbb{R}\backslash U_t$. 
Since $\tau=0$ in this case, it follows by Proposition \ref{prop:Phi-t-gamma-x0-self-adjoint-formula} that
the push-forward map $\Phi_{t,t}$ from $\mu_{{x}+c_t}$ to $\mu_{{x}+g_t}$ can be calculated as 
\begin{align*}
	\Phi_{t,t}(\lambda)
	  &= \begin{cases}
	  \psi_t(a), \qquad & \text{for}\qquad 
	  \lambda\in\Xi_t,\\
	  \lambda+t G_\mu(\lambda), \qquad & \text{for}\qquad 
	  \lambda\notin\Xi_t,\\
	  \end{cases}
\end{align*}
where $\lambda=a+{\i}b$. 
Moreover, for $a\in U_t$, the integration of the density of $\mu_{{x}+c_t}$ over the vertical line segment $\{a+{\i}b: |b|<v_t(a) \}$ is $\frac{v_t(a)}{\pi t}\frac{d\psi_t(a)}{da}$ by the density formula \eqref{density-Brown-x0-ct-derivative-psi}. By the push-forward property from $\mu_{{x}+c_t}$ to $\mu_{{x}+g_t}$, it then follows that 
\[
 p_t(\psi_t(a))=\frac{v_t(a)}{\pi t}.
\]
The measure $\mu_{{x}+g_t}$ has no atom at ending points of components of $U_t$ because $\mu_{{x}+c_t}$ has no atom. Hence, $\mu_{{x}+g_t}$ is absolutely continuous. It follows that
$\mu_{{x}+c_t}(\Xi_t)=\mu_{{x}+g_t}(U_t)=1$. This finishes the proof. 
\end{proof}

\begin{lemma}
\label{lemma:delta-t-derivative}
Let $\delta(a)=\delta_{t,\gamma}(a)= a+\left( 1-\frac{|\tau|^2}{t\tau_1} \right) {h_t(a)}$, then $a\mapsto \delta(a)$ is a homeomorphism of $\mathbb{R}$ onto $\mathbb{R}$. Then, 
\[
   0<\delta'(a)<2
\]
for all $a\in U_t$ and $a\in \big(\overline{U_t} \big)^c$. Moreover, if $a\in U_t$,
\[
    \delta'(a)=\frac{t}{\tau_1} \cdot \det(\emph{Jacobian}(\Phi_{t,\gamma}))(\lambda).
\]
\end{lemma}
\begin{proof}
Using the characterization of $\Xi_t$ and $U_t$ in Proposition \ref{prop:Xi-t-characterization-self-adjoint}, it follows that the result holds for $a\in U_t$. Indeed, we have 
\[
   \delta'(a)=\frac{t}{\tau_1}\cdot  \det(\text{Jacobian}(\Phi_{t,\gamma}))(\lambda)
\]
for any $\lambda\in \Xi_t$ such that $\Re(\lambda)=a$. If $a\in \big(\overline{U_t} \big)^c$, by Lemma \ref{lemma:Cauchy-x0-sefladjoint-in-Omega-t}, then $\supp(\mu)\cap \big(\overline{U_t} \big)^c=\emptyset$. Hence, $v_t(a)=0$ and 
\[
   h_t(a)=t\int_\mathbb{R}\frac{1}{a-u}d\mu(u)=tG_\mu(a).
\]
It follows that $|h_t'(a)|\leq 1$ thanks to Lemma \ref{lemma:Cauchy-x0-sefladjoint-in-Omega-t}. Note that $a\mapsto h_t(a)$ is strictly convex on any open interval in $\big(\overline{U_t} \big)^c$. Hence, $h_t'(a)$ can not take local maximum in the open interval. We then conclude that $0<h_t'(a)<1$ if 
$a\in \big(\overline{U_t} \big)^c$. 
\end{proof}

\begin{figure}[H]
	\label{figure2-circular}
	\begin{center}
		\begin{subfigure}[h]{0.65\linewidth}
			\includegraphics[width=\linewidth]{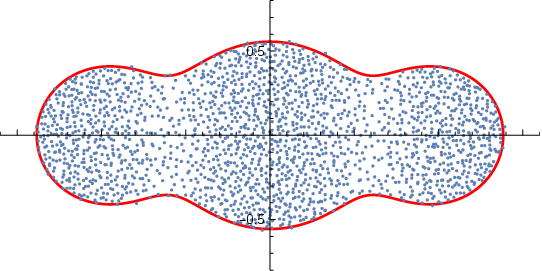}
		\end{subfigure}
\[
 \arrowdown{\huge ${\Phi_{t, \gamma}}$}
\]
		\begin{subfigure}[h]{0.65\linewidth}
			\includegraphics[width=\linewidth]{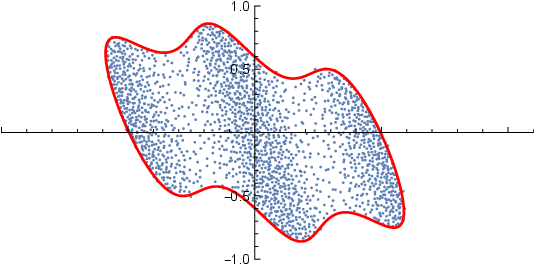}
		\end{subfigure}	
	\end{center}
	\begin{center}
		\caption{The Brown measure simulation (top) of ${x}+c_t$ for $t=0.5$, and
			the Brown measure simulation (bottom) of ${x}+g_{t,\gamma}$ for $t=0.5$ and $\gamma=0.25 + 0.25 i$, where ${x}$ is distributed as $0.25\delta_{-1}+0.5\delta_{0}+0.25\delta_1$.}
	\end{center}
\end{figure}

\begin{theorem}
\label{thm:dentisty-x0-x-t-gamma-main-self-adjoint}
Let ${x}$ be a self-adjoint operator that is free from $g_{t,\gamma}$. 
For any $|\gamma|\leq t$ with $\gamma\neq t$, the map $\Phi_{t,\gamma}$ is non-singular at any $\lambda\in\Xi_t$. The Brown measure of ${x}+g_{t,\gamma}$ is the push-forward map of the Brown measure of ${x}+c_t$ under the map $\Phi_{t,\gamma}$. 

Moreover, the Brown measure $\mu_{{x}+g_{t,\gamma}}$ takes the full measure on $\Phi_{t,\gamma}(\Xi_t)$ and the density is given by
\begin{equation}
   d\mu_{{x}+g_{t,\gamma}}(z)=\frac{1}{2\pi \tau_1}\frac{d\psi_t(a)}{d\delta(a)} dz_1 dz_2, \qquad z\in 
   \Phi_{t,\gamma}(\Xi_t)
\end{equation}
where  $z=z_1+{\i}z_2=\Phi_{t,\gamma}(a+{\i}b)$ and
\[
  \delta(a)= a+\left( 1-\frac{|\tau|^2}{t\tau_1} \right) {h_t(a)}.
\]
\end{theorem}
\begin{proof}
By Theorem \ref{thm:det-Jacobian-Phi}, we deduce that $\Phi_{t,\gamma}$ is non-singular at any $\lambda\in\Xi_t$.
Recall that if $a+{\i}b\in \Xi_t$, then $a\in U_t$. 
Using the density formula for ${x}+c_t$ in \eqref{density-Brown-x0-ct-derivative-psi} and determinant of the Jacobian of $\Phi_{t,\gamma}$ as in \eqref{eqn:det-Jacobian-Phi-t-gamma-self-adjoint}, it then follows that the density of ${x}+g_{t,\gamma}$ can be expressed as
\begin{align*}
     \frac{1}{2\pi t}\frac{d\psi_t(a)/da}{ \det(\text{Jacobian}(\Phi_{t,\gamma})(\lambda))}
       =\frac{1}{2\pi \tau_1}\frac{d\psi_t(a)/da}{d\delta(a)/da}=\frac{1}{2\pi \tau_1}\frac{d\psi_t(a)}{d\delta(a)}.
\end{align*}
The result is established. 
\end{proof}

\begin{example}\cite[Example 5.3]{BianeLehner2001}
\label{corollary:density-ellipse-x-t-gamma}
The Brown measure of $g_{t,\gamma}$ is the uniform measure in the rotated ellipse  with parametrization
\[
e^{{\i}\alpha} \left(\sqrt{t} e^{{\i}\theta}+\frac{|\gamma|}{\sqrt{t}} e^{-{\i}\theta}\right), 
\]
for $\theta\in [0,2\pi]$, where $\alpha$ is such that$\gamma=|\gamma| e^{2{\i}\alpha}$. 
\end{example}
\begin{proof}
Take ${x}=0$ and $\mu=\delta_0$. The formula for the set $\Xi_t$ is simplified as 
\[
  \Xi_t=\{\lambda=a+{\i}b: |\lambda|<\sqrt{t}\}.
\]
The condition determining $w(0;\lambda,t)$ is written as
\[
  \frac{1}{|\lambda|^2+w(0;\lambda,t)^2}=\frac{1}{t}, \quad |\lambda|<\sqrt{t}.
\]
Hence $w(0;\lambda,t)=\sqrt{(t-|\lambda|^2)_+}$. Then $\Phi_{t,\gamma}(\lambda)=\lambda+\gamma \frac{\overline{\lambda}}{|\lambda|^2+w(0;\lambda,t)^2}
    =\lambda+\frac{\gamma}{t}\overline{\lambda}$
when $|\lambda|<\sqrt{t}$. Write $\gamma=|\lambda|e^{{\i}\alpha}$. Then for any $0<r<\sqrt{t}$, the circle 
centered at origin with radius $r$ is mapped to the ellipse with parametrization
\[
   re^{{\i}\theta}+\frac{r|\gamma|}{t}e^{{\i}(\alpha-\theta)}.
\]
Moreover, for $|\lambda|<\sqrt{t}$ and $\lambda=a+{\i}b$, we have 
\[
  h_t(a)=t\frac{a}{|\lambda|^2+w(0;\lambda,t)^2}=a.
\]
Hence, the determinant of the Jacobian of the  map $\Phi_{t,\gamma}$ is the constant 
\begin{equation}
\nonumber
  \det(\text{Jacobian}(\Phi_{t,\gamma}))(\lambda)=\frac{\tau_1}{t}\left[ 1+\left( 1-\frac{|\tau|^2}{t\tau_1} \right) \right], \qquad \lambda=a+{\i}b,
\end{equation}
where $\tau_1+{\i}\tau_2=\tau=t-\gamma$. Recall that the density of the Brown measure of $c_t$ is the uniform measure on the circle $\{\lambda: |\lambda|<\sqrt{t} \}$. Therefore, the Brown measure of $g_{t,\gamma}$ is the uniform measure on the ellipse with parametrization
$\sqrt{t}e^{{\i}\theta}+|\gamma|e^{{\i}(\alpha-\theta)}/\sqrt{t}$ for $\theta\in [0,2\pi]$, where $\alpha$ is determined by $\gamma=|\gamma|e^{2{\i}\alpha}$. The result follows.  
\end{proof}

We now discuss some special cases which allow us to recover main results in \cite{HoHall2020Brown, Ho2020Brown}.

\begin{example}[The semicircular operator]
\label{example-Phi-semicircular}
If $\gamma=t$, the operator $g_{t,\gamma}$ is a semicircular operator $s_t$ with mean zero and variance $t$. Let $x$ be a self-adjoint operator that is freely independent from $s_t$ The push-forward map $\Phi_{t,t}$ is given in Proposition \ref{prop:push-forward-semicircular} (obtained in \cite{HoZhong2020Brown}). More precisely, the push-forward map sends each vertical line segments in $\Xi_t$ to a single point.

In this case, the injectivity of $\Phi_{t,t}$ fails dramatically, and the statement of Lemma \ref{lemma:FK-det-derivatives-connection-ct-gt} does not hold. Indeed, by Theorem \ref{thm:sub-X0-gt-formula}, for any $\varepsilon>0$, we have 
\[
  p_z^{g,(t,t)}(\varepsilon)=p_\lambda^{(0)}(w(\varepsilon))
\]
where $w(\varepsilon)=w(\varepsilon;\lambda,t)$ and $z=\Phi_{t,\gamma}(\lambda)$. By the Notation \ref{notation:partial-derivatives-lambda-c-x}, This can be rewritten as 
\begin{align*}
  -\phi&\left[({x}+g_t-z\unit)^*\big( ({x}+g_t-z\unit)({x}+g_t-z\unit)^*+
   \varepsilon^2\unit\big)^{-1}\right] \\
    &\qquad \qquad\qquad =-\phi\left[({x}-\lambda\unit)^*\big( ({x}-\lambda\unit)({x}-\lambda\unit)^*+
   w(\varepsilon)^2\unit\big)^{-1}\right].
\end{align*}
By Lemma \ref{lemma:s-epsilon-t-limit}, for any $\lambda\in\Xi_t$ (equivalently, $t>\lambda_1(\mu_{|{x}-\lambda\unit|})$), we have
\[
 \begin{aligned}
    \lim_{\varepsilon\rightarrow 0^+} p_\lambda^{(0)}(w(\varepsilon))&=-\phi\left[({x}-\lambda\unit)^*\big( ({x}-\lambda\unit)({x}-\lambda\unit)^*+
   w(0;\lambda,t)^2\unit\big)^{-1}\right]\\
    &=p_\lambda^{(0)}(w(0;\lambda,t)).
 \end{aligned}
\]
Hence, for any $\lambda\in\Xi_t$, by choosing $z=\lambda(t)=\lambda+\gamma p_\lambda^{(0)}(w(\varepsilon))$, we have
\begin{align*}
	  \lim_{\varepsilon\rightarrow 0^+}  p_z^{g,(t,t)}(\varepsilon)
	&=-\phi\left[({x}-\lambda\unit)^*\big( ({x}-\lambda\unit)({x}-\lambda\unit)^*+
	w(0;\lambda,t)^2\unit\big)^{-1}\right]\\
	&=p_\lambda^{(0)}(w(0;\lambda,t)).
\end{align*}
We note that $p_\lambda^{(0)}(w(0;\lambda,t))$ can be expressed as
\begin{align*}
  p_\lambda^{(0)}(w(0;\lambda,t))&=\int_\mathbb{R}\frac{\overline{\lambda}-u}{(a-u)^2+b^2+w(0;\lambda,t)^2}d\mu(u)\\
   &=\int_\mathbb{R}\frac{a-u}{(a-u)^2+v_t(a)^2}d\mu(u)-\frac{{\i}b}{t},
\end{align*}
where we used \eqref{eqn:condition-w-0-x0-self-adjoint} and the fact that $b^2+w(0;\lambda,t)^2=v_t(a)^2$ for $a+{\i}b\in\Xi_t$. 

To summarize, for $a\in\mathbb{R}$ fixed so that $v_t(a)>0$, then: (1) the limit $p_z^{g, (t,t)}(\varepsilon)$ as $\varepsilon$ tends to zero have different limit as long as $(z, \varepsilon)$ tends to $(\psi_t(a), 0)$ along different paths $(\lambda(t), \varepsilon)$ depending on $b$; (2) although the $\lim_{\varepsilon\rightarrow 0^+}\lambda(t)=\Phi_{t,t}(a+{\i}b)=\psi_t(a)$ for any $-v_t(a)\leq b\leq v_t(a)$, the limit of the partial derivative $ \lim_{\varepsilon\rightarrow 0^+}  p_z^{g,(t,t)}(\varepsilon)$ detects the value $b$ and 
\emph{remembers} where it came from. Namely, by looking at the limit $\lim_{\varepsilon\rightarrow 0^+}\lambda(t)=\Phi_{t,t}(a+{\i}b)$ we can not identify $b$ since $\Phi_{t,t}$ is not one-to-one, but the limit of the partial derivative $ \lim_{\varepsilon\rightarrow 0^+}  p_z^{g,(t,t)}(\varepsilon)$ does the work.
This phenomena is different from what Lemma \ref{lemma:FK-det-derivatives-connection-ct-gt} told us when $\Phi_{t,\gamma}$ is one-to-one.
\end{example}

\begin{example}[The imaginary multiple of a semicircular operator]
\label{example-Phi-i-semicircular}
If $\gamma=it$, the operator $g_{t,\gamma}$ has the same distribution as ${\i}s_t$. Let $x$ be a self-adjoint operator freely independent from $s_t$. In this case, $\tau=2t$, and for $\lambda\in\Xi_t$, we have
\begin{align*}
  \Phi_{t,it}(\lambda)&=\psi_t(a)-2h_t(a)+2{\i}b\\
    &=a-h_t(a)+2{\i}b\\
    &=t\int_\mathbb{R}\frac{u}{(u-a)^2+v_t(a)^2}d\mu(u) +2{\i}b
\end{align*}
and
\begin{align*}
   \delta_{t,it}(a)&=a-h_t(a)=t\int_\mathbb{R}\frac{u}{(u-a)^2+v_t(a)^2}d\mu(u).
\end{align*}
\end{example}

\begin{example}[The elliptic operator]
\label{example-Phi-elliptic}
If $\gamma=s\in \mathbb{R}$ with $-t<s<t$, the operator $g_{t,\gamma}$ has the same distribution as an elliptic operator. In this case, $\tau=t-s$, and for $\lambda\in\Xi_t$, we have
\begin{align*}
  \Phi_{t,s}(\lambda)&=\psi_t(a)-\frac{t-s}{t}h_t(a)+{\i}\frac{t-s}{t}b\\
    &=a+\frac{s}{t}h_t(a)+{\i}\frac{t-s}{t}b\\
  &=\left( a+s\int_\mathbb{R}\frac{a-u}{(a-u)^2+v_t(a)^2}d\mu(u) \right)+{\i}\frac{(t-s)b}{t},
\end{align*}
and
\begin{align*}
  \delta_{t,s}(a)=a+s\int_\mathbb{R}\frac{a-u}{(a-u)^2+v_t(a)^2}d\mu(u).
\end{align*}

We also note that $\Phi_{t,s}(\lambda)=z_1+{\i}z_2$ where 
\[
  z_1=\delta_{t,s}(a), \qquad z_2=\frac{(t-s)b}{t}
\]
and $\psi_t(a)$ can be written as
\[
  \psi_t(a)=z_1+\frac{t-s}{t}h_t(a)=z_1+(t-s)\int_\mathbb{R}\frac{a-u}{(a-u)^2+v_t(a)^2}d\mu(u).
\]
Therefore, the density of the Brown measure at its support $\overline{\Phi_{t,s}(\Xi_t)}$ is given by
\begin{align*}
&d\mu_{{x}+g_{t,s}}(z_1+{\i}z_2)\\
=&\frac{1}{2\pi (t-s)}\left( 1+(t-s)\frac{d}{dz_1}\int_\mathbb{R}\frac{a-u}{(a-u)^2+v_t(a)^2}d\mu(u) \right)dz_1dz_2,
\end{align*}
if $v_t(a)>0$.
Here we reminder the reader that the map $a\mapsto z_1=\delta_{t,s}(a)$ is a homeomorphism from $\mathbb{R}$ to $\mathbb{R}$. This recovers the main result in \cite{HoHall2020Brown, Ho2020Brown}
\end{example}

\begin{remark}
\label{remark:twisted-coordinate-x0-self-adjoint}
[The twisted $(\nu,\delta)-$coordinate.]
Fix $a\in\mathbb{R}$, as in the proof of Theorem \ref{thm:det-Jacobian-Phi} (see \eqref{eqn:z-1-z-2}), the map $b\mapsto\Phi_{t,\gamma}(a+{\i}b)$ is an affine transform of $b$. The density formula of ${x}+c_t$ in Theorem \ref{thm:density-x0-c-t-self-adjoint} is independent of $b$. Hence, the density formula of ${x}+g_{t,\gamma}$ must depend only on one parameter. It is indeed the case as in Theorem \ref{thm:dentisty-x0-x-t-gamma-main-self-adjoint}, where the density is expressed in terms of parameter $a$ coming from $\Xi_t$, the support of ${x}+c_t$.

We now describe an analogue of the formulation in the recent work \cite{HallHo2021Brown} where the authors study the free multiplicative Brownian motions as follows. For $z=z_1+{\i}z_2$, consider twisted $(\nu,\delta)-$coordinate determined by
\[
   a+{\i}b={\i}\tau \nu+\delta,
\]
where $\nu, \delta\in\mathbb{R}$. They can be written as 
\[
   \delta=a+\frac{\tau_2}{\tau_1}b, \qquad \nu=\frac{b}{\tau_1}.
\]
Using notations in the proof of Theorem \ref{thm:det-Jacobian-Phi} and the formula of $\Phi_{t,\gamma}(\lambda)$ as in \eqref{eqn:z-1-z-2}, for $\lambda=\lambda_1+i\lambda_2\in \Xi_t$, we have 
\begin{align*}
  \delta(\Phi_{t,\gamma}(\lambda))&=z_1+\frac{\tau_2}{\tau_1}z_2\\
    &=\left[a+\left(1-\frac{\tau_1}{t}\right) h_t(a)-\frac{\tau_2}{t}b\right]+\frac{\tau_2}{\tau_1}
       \left(-\frac{\tau_2}{t}h_t(a)+\frac{\tau_1}{t}b \right)\\
     &= a+\left( 1-\frac{|\tau|^2}{t\tau_1} \right) {h_t(a)}
     =\delta_{t,\gamma}(a),
\end{align*}
and
\[
   \nu(\Phi_{t,\gamma}(\lambda))=\frac{z_2}{\tau_1}=\frac{1}{t}\left[ b-\frac{\tau_2}{\tau_1}h_t(a) \right].
\]
We can then say that the Brown measure $\mu_{{x}+g_{t,\gamma}}$ is constant along the $\nu$-direction. 
\end{remark}


\section{Addition with an $R$-diagonal operator}
\label{section:addition-R-diag}

The family of $R$-diagonal operators was introduced by Nica-Speicher \cite{NicaSpeicher-Rdiag} which covers a large class of interesting operators in free probability theory. The circular operator and Haar unitary operator are special examples of $R$-diagonal operators. They have a number of remarkable symmetric properties. In a breakthrough paper, Haagerup--Larsen \cite{HaagerupLarsen2000} calculated the Brown measure of any bounded $R$-diagonal operator and they showed that the Brown measure of $T$ can be expressed in terms of the $S$-transform of the operator $T^*T$. This is the first nontrivial example of Brown measure formula after L. Brown introduced the definition in 1983. In \cite{HaagerupSchultz2007}, Haagerup-Schultz gave a second proof of the Brown measure formula of $R$-diagonal operators which also work for unbounded $R$-diagonal operators. In \cite{Zhong2021RdiagBrown}, we used subordination ideas and obtained simplification of technical arguments in Haagerup-Schultz's work. 

The Brown measure of an $R$-diagonal operator is the limit of the ESD of certain non-normal random matrix model called single ring theorem named after Feinberg and Zee (see \cite{GuionnetKZ-single-ring2011} and the earlier physics paper \cite{FeinbergZee1997-single-ring, FeinbergSZ2001-single-ring-2}).

In this section, we apply the results in Section \ref{section:Brown-free-additiveBM} and Section \ref{section:Brown-addition-elliptic} to study the Brown measures of $T+c_t$ and $T+g_{t,\gamma}$ where  $T$ and $\{c_t, g_{t,\gamma}\}$ are $*$-free. 
It is well-known that the sum of two $R$-diagonal operator is again $R$-diagonal. The new observation from \cite{Zhong2021RdiagBrown} will be used in Section \ref{section:Phi-t-gamma-Rdiag} to describe the push-forward map. We show that the Brown measure of $T+g_{t,\gamma}$ is supported in a deformed ring where the inner boundary is a circle and the outer boundary is an ellipse. The push-forward map $\Phi_{t, \gamma}$ sends a family of circles to a family of ellipses. 

\subsection{The work of Haagerup--Schultz and gradient functions}

Following \cite[Section 4]{HaagerupSchultz2007}, we introduce some auxiliary functions.
Let $\lambda\in\C\setminus\{0\}$, we set 
\begin{align}
	h_T(s)&= s\,\phi\big((T^* T+s^2\unit)^{-1}\big), \qquad & s>0,\\
	h_{T-\lambda\unit}(\varepsilon) &= \varepsilon\,\phi\big([(T-\lambda\unit)^*
	(T-\lambda\unit)+\varepsilon^2\unit]^{-1}\big), \qquad  &\varepsilon>0.
\end{align}
and
\[
  \lambda_1(T)=1/\sqrt{\phi((T^*T)^{-1})}, \qquad \lambda_2(T)=\sqrt{\phi(T^*T)}.
\]
It is known that the Brown measure of $T$ is supported in a single ring and $\lambda_1(T), \lambda_2(T)$ are inner and outer radii of this ring.

\begin{proposition}
	\cite[Definition 4.9]{HaagerupSchultz2007}
	\label{defn-for-s-t-Rdiag}
For any $\lambda\in\C\setminus\{0\}$, the equation
\begin{equation}\label{eqn:fixed-point-s-t-Rdiag}
(s-\varepsilon)^2-\frac{s-\varepsilon}{h_T(s)}+|\lambda|^2=0
\end{equation}
has a unique solution $s=s(|\lambda|,t)$ in the interval $(0,\infty)$ and 
\[
h_{T-\lambda\unit}(\varepsilon)=h_T(s(|\lambda|,\varepsilon)), \qquad \varepsilon>0. 
\]

	For $\lambda\in (\lambda_1(T),\lambda_2(T))$, the equation 
\begin{equation}
\label{eqn:fixed-point-s-t-0-Rdiag}
s^2-\frac{s}{h_T(s)}+|\lambda|^2=0
\end{equation}
has a unique solution  $s=s(|\lambda|,0)$ in the interval $(0,\infty)$. 
\end{proposition}

\begin{proposition}
	\label{prop:limit-s-lambda-epsilon-Rdiag}
	\cite[Lemma 4.10, Remark 4.11]{HaagerupSchultz2007}
	The function $(\lambda,\varepsilon)\mapsto s(\lambda,\varepsilon)$ is analytic in $(0,\infty)\times (0,\infty)$. Moreover, 
	\begin{align*}
		\lim_{
			\varepsilon\rightarrow 0}s(\lambda,\varepsilon)= \begin{cases}
			0, &\quad \text{if} \quad 0<\lambda\leq \lambda_1(T);\\
			s(\lambda,0), &\quad\text{if}\quad \lambda\in (\lambda_1(T),\lambda_2(T));\\
			+\infty,&\quad \text{if}, \quad \lambda\geq \lambda_2(T).
		\end{cases}
	\end{align*} 
	From now on, we denote by $s(\lambda,0)=\lim_{\varepsilon\rightarrow 0}s(\lambda,\varepsilon)$ for any $\lambda>0$. 
\end{proposition}

\begin{lemma}\label{lemma-delta-star-s-t} 
	\cite[Lemma 4.14]{HaagerupSchultz2007}
	Let $T$ be an $R$--diagonal element in $\CM$,
	let $\lambda\in\C\setminus\{0\}$, and let $\varepsilon>0$. We then have:
	\begin{equation}\label{eqn:delta-star-s-t}
		\Delta\big((T-\lambda\unit)^*(T-\lambda\unit)+\varepsilon^2\unit\big)=
		\frac{|\lambda|^2}{|\lambda|^2 + (s(|\lambda|,\varepsilon)-\varepsilon)^2} \Delta\big(T^*
		T + s(|\lambda|,\varepsilon)^2\unit\big).
	\end{equation}
\end{lemma}

\begin{theorem}\label{thm:FK-det-Rdiag-7.6}
	\cite[Theorem 4.15]{HaagerupSchultz2007}
	Let $T\in\CM$ be $R$--diagonal, we have 
	\begin{itemize}
		\item[(i)] If $\lambda_1(T)<|\lambda|<\lambda_2(T)$, then
		\begin{equation}
			\label{eqn:FK-det-Rdiag}
			\Delta(T-\lambda\unit)=
			\Bigg(\frac{|\lambda|^2}{|\lambda|^2+s(|\lambda|,0)^2}\,
			\Delta(T^* T +s(|\lambda|,0)^2\unit)\Bigg)^\frac12.
		\end{equation}
		\item[(ii)] If $|\lambda|\leq \lambda_1(T)$, then
		$\Delta(T-\lambda\unit)= \Delta(T)$.
		\item[(iii)] If $|\lambda|\geq \lambda_2(T)$, then
		$\Delta(T-\lambda\unit)=|\lambda|$.
	\end{itemize}
\end{theorem}

\begin{remark}
The reader might notice that formulas in this section are similar to our results in Section \ref{section-FK-det-subordination}. Indeed, circular operator is an $R$-diagonal operator. 
In \cite{Zhong2021RdiagBrown}, we used subordination ideas to give a simplified proof for Haagerup-Schultz's results and the Brown measure formula of $R$-diagonal operators. 
By choosing ${x}=0$, then Theorem \ref{thm:main-FK-det-ct-0} is a special case of Theorem \ref{thm:FK-det-Rdiag-7.6}. In a joint work with Bercovici \cite{BercoviciZhong2022}, we obtained a Fuglede-Kadison formula for operator $T+{x}$ where ${x}$ is an arbitrary operator $*$-free from $T$, which generalizes some free probability results obtained in \cite{DykemaSZ2017unbounded, HaagerupSchultz2009}.
\end{remark}

The following result can be deduced from the Fuglede-Kadison formulas in Lemma \ref{lemma-delta-star-s-t}
and Theorem \ref{thm:FK-det-Rdiag-7.6} and the defining equations for $s(|\lambda|,\varepsilon)$ and $s(|\lambda|,0)$. See \cite{Zhong2021RdiagBrown} for details. 
\begin{lemma}
	\cite[Section 4]{Zhong2021RdiagBrown}
	\label{lemma:formula-derivative-t-0}
	Let $T$ be an $R$-diagonal operator. 
	We have partial derivative formula
	 \begin{equation}
	 \phi( (\lambda\unit-T)^*[(\lambda \unit -T)(\lambda\unit- T)^*+\varepsilon^2]^{-1}) 
	   =\frac{\overline\lambda}{|\lambda|^2}\frac{(s(|\lambda|,\varepsilon)-\varepsilon)^2 }{(s(|\lambda|,\varepsilon)-\varepsilon)^2 +|\lambda|^2 },
	 \end{equation}
	 for any $\varepsilon>0$, and 
	\begin{align*}
	&\lim_{\varepsilon\downarrow {0}}\phi( (\lambda\unit-T)^*[(\lambda \unit -T)(\lambda\unit- T)^*+\varepsilon^2]^{-1}) \\
	=&
	\begin{cases}
	0, \quad & \text{for} \quad 0< |\lambda|\leq \lambda_1(T);\\
	\frac{\overline{\lambda}}{|\lambda|^2}\frac{s(|\lambda|,0)^2}{s(|\lambda|,0)^2+|\lambda|^2},\quad &\text{for} \quad \lambda_1(T)<|\lambda|<\lambda_2(T);\\
	\frac{\overline{\lambda}}{|\lambda|^2},\quad &\text{for}\quad |\lambda|\geq \lambda_2(T),
	\end{cases}
	\end{align*}
	where $s(|\lambda|,\varepsilon)$ and $s(|\lambda|,0)$ are defined in Definition \ref{defn-for-s-t-Rdiag}. Moreover, the Brown measure $\mu_T$ is the rotationally invariant probability measure such that
	\[
	\mu_T\{ z\in\mathbb{C}: |z|\leq r \}=\frac{s(r,0)^2}{s(r,0)^2+r^2}, \qquad 0<r<\infty
	\]
\end{lemma}

\subsection{The push-forward map and Brown measure}
\label{section:Phi-t-gamma-Rdiag}
We use Lemma \ref{lemma:formula-derivative-t-0} to study properties of the push-forward map from the Brown measure of $T+c_t$ to the Brown measure of $T+g_{t,\gamma}$. Observe that the Brown measure of $T+e^{{\i}\theta} g_{t,\gamma}$ is the same as the Brown measure of $e^{{\i}\theta}(T+g_{t,\gamma})$. In this section, we may restrict ourselves to the case $\gamma\in\mathbb{R}$, but the general result will follow easily. So we still keep using the complex notation $\gamma$ as before.

\begin{proposition}
\label{prop:Phi-t-gamma-formula-Rdiag}
	For ${x}=T$, the map $\Phi_{t,\gamma}$ is expressed as
\begin{align}
\label{eqn:Phi-t-gamma-Rdiag-formulas}
\Phi_{t,\gamma}(\lambda)=\lambda+\frac{\gamma}{\lambda}\frac{s(|\lambda|,0)^2}{s(|\lambda|,0)^2+|\lambda|^2}
=\lambda+\gamma\cdot\frac{\mu_{T+c_t}(\{z\in\mathbb{C} : |z|\leq |\lambda| \})}{\lambda},
\end{align}
for any $\lambda\in\mathbb{C}$, where for $\lambda\in\Xi_t$, $s=s(|\lambda|,0)$ is determined by
\[
s^2-\frac{s}{h(s)}+|\lambda|^2=0
\]
and $h(s)=s\cdot \phi( ( (T+c_t)^*(T+c_t)+s^2 )^{-1})$.

The regularized map $\Phi^{(\varepsilon)}_{t, \gamma}$ is expressed as
\begin{equation}
 \label{eqn:Phi-t-gamma-Rdiag-regularization}
\Phi^{(\varepsilon)}_{t,\gamma}(\lambda)=\lambda+\frac{\gamma}{\lambda}\frac{(s(|\lambda|,\varepsilon)-\varepsilon)^2}{(s(|\lambda|,\varepsilon)-\varepsilon)^2+|\lambda|^2},
\end{equation}
where $s=s(|\lambda|,\varepsilon)$ is determined by 
\[
(s-\varepsilon)^2-\frac{s-\varepsilon}{h(s)}+|\lambda|^2=0.
\]
\end{proposition}
\begin{proof}
We have the following equivalent definition of the regularized push-forward map
\[
  \Phi^{(\varepsilon)}_{t, \gamma}(\lambda)=\lambda+\frac{\partial}{\partial\lambda}S(T+c_t,\lambda,\varepsilon).
\]
The result follows by applying Lemma \ref{lemma:formula-derivative-t-0}.
\end{proof}

Since $T$ and $c_t$ are free to each other, we can calculate
\[
  \phi( (T+c_t)^*(T+c_t) )=\phi(T^*T)+\phi(c_t^*c_t)=\lambda_2(T)^2+t.
\]
Hence, $\lambda_2(T+c_t)=\sqrt{\lambda_2(T)^2+t }$. Although it is not important here, one can actually show that $\lambda_1(T+c_t)=\sqrt{(\lambda_1(T)^2-t)|_+}$ (\cite{BercoviciZhong2022}). 
When ${x}$ is an $R$-diagonal operator $T$, the set $\Xi_t$ is also the interior of the support of the Brown measure of $T+c_t$. That is,  
\[
  \Xi_t=\{ \lambda\in\mathbb{C}:  (\lambda_1(T)^2-t)|_+  < |\lambda|^2<\lambda_2(T)^2+t \}.
\]

\begin{figure}[t]
	\centering
	\begin{minipage}[c]{\textwidth}
		\begin{center}
			\begin{subfigure}[h]{0.4\linewidth}
				\includegraphics[width=\linewidth]{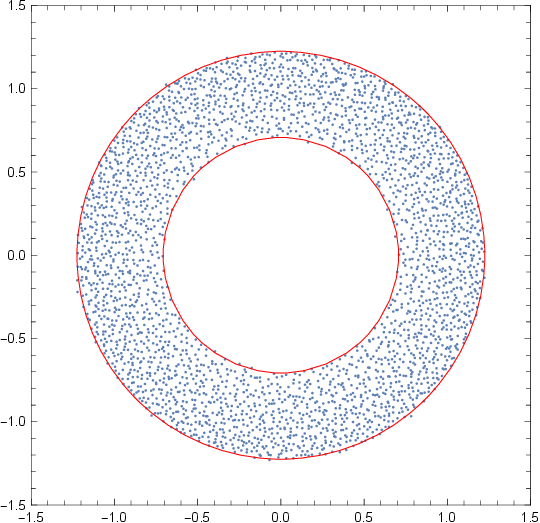}
			\end{subfigure}
			\begin{tikzpicture}
			\draw[-Latex,line width=4pt] (0,0) -- (1.5,0) node [pos=0.6, above] {\huge ${\Phi_{t, \gamma}}$};
			\end{tikzpicture}
			\begin{subfigure}[h]{0.45\linewidth}
				\includegraphics[width=\linewidth]{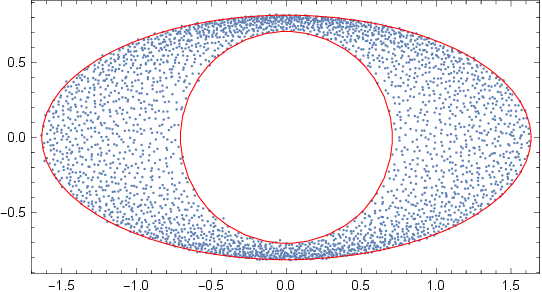}
			\end{subfigure}	
		\end{center}
		\caption{The $4000\times 4000$ random matrix simulation for Brown measures of $u+c_t$ and $u+g_t$ with $t=0.5$.}
		\label{figure3-circular}
	\end{minipage}
\end{figure}

\begin{theorem}
 \label{thm:Brown-addition-R-diag-elliptic-pushforward}
 Let $\gamma=|\gamma|e^{{\i}2\alpha}$ such that $|\gamma|\leq t$. 
 Set ${x}=T$ and use notations in Section \ref{section:Brown-addition-elliptic}. For any $r>0$, 
 the push-forward map $\Phi_{t,\gamma}$ sends the circle $\text{C}_r$ centered at the origin with radius $r$ to the ellipse
 \begin{align}
 \label{eqn:ellipse-Phi-t-gamma-T-ct}
 e^{{\i}\alpha}\left( re^{{\i}\theta} +\frac{e^{-{\i}\theta} |\gamma| m(r)}{r}\right), \quad 0\leq \theta\leq 2\pi, 
 \end{align}
 where 
 \[
 m(r)={\mu_{T+c_t}(\{z\in\mathbb{C} : |z|\leq r \})}.
 \]
 The semi-axes of the ellipse are 
 \[
 a(r)=r-|\gamma| m(r)/r, \qquad b(r)=r+|\gamma|m(r)/r. 
 \]
Moreover, both $a(r)$ and $b(r)$ are increasing functions of $r$ in the interval $(0,\infty)$.
When $|\gamma|<t$, then both $a(r)$ and $b(r)$ are strictly increasing and the map $\Phi_{t,\gamma}$ is a self-homeomorphism of $\mathbb{C}$.
 
In particular, the Brown measure of $T+g_{t,\gamma}$ is supported in the deformed ring where the inner boundary is the circle with radius $\lambda_1(T+c_t)$ and the outer boundary is the ellipse 
\[
  e^{{\i}\alpha}\left( \lambda_2 e^{{\i}\theta} +\frac{e^{-{\i}\theta} |\gamma|}{ \lambda_2} \right), \quad 0\leq \theta\leq 2\pi, 
\]
where $\lambda_2=\lambda_2(T+c_t)$ is the outer radii of the support of the Brown measure of $T+c_t$. 
\end{theorem}
\begin{proof}
We follow notations used in Proposition \ref{prop:Phi-t-gamma-formula-Rdiag}. 
 Recall that, by Lemma \ref{lemma:formula-derivative-t-0}, we have 
\[
m(r)={\mu_{T+c_t}(\{z\in\mathbb{C} : |z|\leq r \})}=\frac{s(r,0)^2}{s(r,0)^2+r^2}.
\]
We can rewrite $a(r), b(r)$ as
\[
a(r)=r-\frac{\gamma}{r}\frac{s(r,0)^2}{s(r,0)^2+r^2}, \qquad b(r)=r+\frac{\gamma}{r}\frac{s(r,0)^2}{s(r,0)^2+r^2}.
\]
For any $\varepsilon>0$, we set
\[
a(r,\varepsilon)= r-\frac{\gamma}{r}\frac{(s(r,\varepsilon)-\varepsilon)^2}{(s(r,\varepsilon)-\varepsilon)^2+r^2},
\qquad 
b(r,\varepsilon)=  r+\frac{\gamma}{r}\frac{(s(r,\varepsilon)-\varepsilon)^2}{(s(r,\varepsilon)-\varepsilon)^2+r^2}.
\]
Then, by \eqref{eqn:Phi-t-gamma-Rdiag-regularization}, we see that the regularized map $\Phi^{(\varepsilon)}_{t, \gamma}$ sends $C_r$ to an ellipse $E_{r,\varepsilon}$ centered at the origin, whose semi-axes  
$|a(r,\varepsilon)|, b(r,\varepsilon)$. 
 Since $a(r,\varepsilon)>0$ for large $r$ and $\Phi^{(\varepsilon)}_{t, \gamma}$ is a homeomorphism (see Proposition \ref{prop:regularization-Phi-one2one-map}), it follows that $a(r,\varepsilon)>0$ for any $r>0$. Moreover, it forces that the region enclosed by the ellipse $E_{r,\varepsilon}$ increases as $r$ increases. That is, 
\[
   a(r_1,\varepsilon)<a(r_2,\varepsilon), \qquad b(r_1,\varepsilon)<b(r_2,\varepsilon)
\]
for any $0<r_1<r_2$. By letting $\varepsilon$ go to zero and using Proposition \ref{prop:limit-s-lambda-epsilon-Rdiag} (or the fact that $\Phi^{(\varepsilon)}_{t, \gamma}$  converges uniformly to $\Phi_{t, \gamma}$ in $\mathbb{C}$ by Lemma \ref{lemma:regularization-Phi-uniform-convergence}), we deduce  
\[
   a(r_1)\leq a(r_2), \qquad b(r_1)\leq b(r_2), \qquad 0<r_1<r_2. 
\]

 Note that $a(r)=b(r)=r$ for $0<r\leq\lambda_1(T+c_t)$ and
\[
m(r)=\begin{cases}
0, \quad & \text{for}\quad 0\leq r\leq \lambda_1(T+c_t),\\
1, \quad & \text{for}\quad r> \lambda_2(T+c_t).
\end{cases}
\]	
Hence, for $r>\lambda_1(T+c_t)$, we have $a(r)>\lambda_1(T+c_t)$. Therefore, $a(r)>0$ for all $r$. 
 From \eqref{eqn:Phi-t-gamma-Rdiag-formulas}, we then see that $\Phi_{t, \gamma}$ maps the circle $C_r:=\{z: |z|=r \}$ to an ellipse $E_r$ with parametrization 
 \[
 e^{{\i}\alpha} \left( re^{{\i}\theta} +  \frac{|\gamma| e^{-{\i}\theta} m(r)}{r}  \right), \qquad 0\leq \theta\leq 2\pi,
 \]
 whose semi-axes are $a(r), b(r)$. 
 
 When $|\gamma|<t$, we write
\begin{equation}
    \label{eqn:a-r-in-proof}
     a(r)=r\left(1-\frac{|\gamma|}{t}\right)+\frac{|\gamma|}{t}\left( r- \frac{tm(r)}{r} \right).
\end{equation}
The previous discussion shows that the second term of the right hand side of \eqref{eqn:a-r-in-proof} is increasing. Therefore, $a'(r)>0$ for all $r\in (0,\infty)$. Similary, $b'(r)>0$ for all $r\in (0,\infty)$. This implies that $\Phi_{t,\gamma}$ is a self-homeomorphism of $\mathbb{C}$ that maps circles centered at the origin to a family of ellipses.
\end{proof}


\section{Example of explicit formulas}
\label{section:Brown-examples-non-selfadjoint}

In this section, we compute two examples of $x_0$ for the Brown measure of $x_0+c_t$ and $x_0+g_{t,\gamma}$. The first example is quasi-nilpotent DT operator, and the second example is Haar unitary operator.

\subsection{The quasi-nilpoten DT operator}

The quasi-nilpotent DT-operator $T$ was introduced by Dykema--Haagerup \cite{DykemaHaagerup2004-DT}. The operator played a key role in Brown measure of free random variables and invariant subspace problem in type $II_1$ factors \cite{DykemaHaagerup2001JFA, HaagerupSchultz2009}. It can be described as the limit in $*$-moments of random matrices of the form
\begin{align*}
  T^{(n)}=\begin{bmatrix}
     0 & t_{1,2}& \cdots & t_{1,n}\\
     0 & \ddots & \ddots &\vdots\\
     \vdots & \ddots &\ddots & t_{n-1,n}\\
     0 &\cdots & 0 & 0
  \end{bmatrix}
\end{align*}
where $\{ \Re(t_{i,j}), \Im (t_{i,j})\}$ is a set of $n(n-1)$ independent identically distributed Gaussian random variables with mean zero and variance $\frac{1}{2n}$. We refer the reader to \cite{DykemaHaagerup2001JFA, DykemaHaagerup2004-DT,HaagerupAagaard2004} for the construction of $T$ in a finite von Neumann algebra. 

\begin{proposition}
Given $\lambda\in\mathbb{C}$
Let $\sigma\in\mathbb{R}$ with $1+{|\lambda|^2}\sigma>0$ and $\mu^2=-\frac{e^\sigma}{\sigma
}(1+|\lambda|^2\sigma)$. We have 
\begin{align*}
  \begin{cases}
    &\phi(( (T-\lambda\unit)^*(T-\lambda\unit)+\mu^2)^{-1})=e^{-\sigma}-1\\
    &\phi( (T-\lambda\unit) ((T-\lambda\unit)^*(T-\lambda\unit)+\mu^2))^{-1} )=\lambda\sigma\\
    &\phi( (T-\lambda\unit)^* ((T-\lambda\unit)^*(T-\lambda\unit)+\mu^2))^{-1} )=\overline\lambda\sigma.
  \end{cases}
\end{align*}
\end{proposition}
\begin{proof}
It can be obtained by taking the expectation of (3.3) in \cite{HaagerupAagaard2004}.
See \cite[Pages 586-588]{HaagerupAagaard2004}.
\end{proof}

\begin{proposition}{\cite[Theorem 4.3]{HaagerupAagaard2004}}
For any $t>0$, 
the Brown measure of $T+c_t$ is the uniform measure on the closed disk $\overline{B}(0, \frac{1}{\sqrt{\log(1+1/t)}})$.
\end{proposition}
\begin{proof}
Choose $\sigma$ so that $e^{-\sigma}-1=\frac{1}{t}$. Then $\sigma=-\log(1+1/t)$. Consider the solution of $\varepsilon>0$ for the equation 
\[
   \phi( (T-\lambda\unit)^*(T-\lambda\unit)+\varepsilon^2)=\frac{1}{t},
\]
we then deduce that $w(0;\lambda,t)^2=\mu^2=-\frac{e^\sigma}{\sigma
}(1+|\lambda|^2\sigma)$ by applying the results in Section 5 for $x_0=T$. 
Hence, $w(0;\lambda,t)>0$ if and only if $(1+|\lambda|^2\sigma)>0$. Therefore,
\begin{align*}
    \Xi_t&=\left\{ \lambda: \phi(( (T-\lambda\unit)^*(T-\lambda\unit)+\mu^2)^{-1})>\frac{1}{t} \right\}\\
    &=\{\lambda: w(0;\lambda,t)>0 \}
   =\{\lambda: 1+|\lambda|^2\sigma>0 \}\\
   &=\left\{\lambda: |\lambda|< \frac{1}{\sqrt{\log(1+1/t)}}  \right\}.
\end{align*}
In addition, by the Proof of Theorem \ref{thm:BrownFormula-x0-ct-general}, we have 
\begin{align*}
  & \frac{\partial}{\partial \lambda} 
    \log  \Delta \big( (T+c_t-\lambda\unit)^*(T+c_t-\lambda\unit) \big)\\
   =&  p_\lambda^{(t)}(0) =p_\lambda^{(0)}(w(0;\lambda,t))\\
    =&-\phi \big( (T-\lambda\unit)^* ( (T-\lambda\unit)^*(x_0-\lambda\unit)+w(0; \lambda,t)^2 )^{-1} \big)\\
    =&-\phi \big( (T-\lambda\unit)^* ( (T-\lambda\unit)^*(x_0-\lambda\unit)+\mu^2 )^{-1} \big)
    =-\overline{\lambda}\sigma.
\end{align*}
Hence, the density of $T+c_t$ at $\lambda\in \Omega_t={B}(0, \frac{1}{\sqrt{\log(1+1/t)}})$ is given by
  \[
       d\mu_{T+c_t}=\frac{1}{\pi}\frac{\partial^2}{\partial \lambda^2} 
    \log  \Delta \big( (T+c_t-\lambda\unit)^*(T+c_t-\lambda\unit) \big)=-\frac{\sigma}{\pi}
      =\frac{\log(1+1/t)}{\pi}.
  \]
 It follows that the Brown measure of $T+c_t$ is equal to the uniform measure on the 
 $\overline{B}(0, \frac{1}{\sqrt{\log(1+1/t)}})$.
\end{proof}

We then calculate the push-forward map from $T+c_t$ to $T+g_{t,\gamma}$. We have
\[
 \Phi_{t,\gamma}(\lambda)= \lambda-\gamma \phi \big( (T-\lambda\unit)^* ( (T-\lambda\unit)^*(x_0-\lambda\unit)+w(0;\lambda,t)^2 )^{-1} \big)
  =\lambda-\gamma\overline{\lambda}\sigma,
\]
where $\sigma=-\log(1+1/t)$. Therefore,   
\[
   \Phi_{t,\gamma}(\lambda)=\lambda+\overline{\lambda}\gamma\cdot \log(1+1/t).
\]
Recall that $|\gamma|\leq t$. 
In particular, it follows that $\Phi_{t,\gamma}$ is non-singular in the disk $B(0, \frac{1}{\sqrt{\log(1+1/t)}})$. Hence, by the pushforward connection in Theorem \ref{thm:main-push-forward-general}, the Brown measure of $T+g_{t,\gamma}$ is supported in the ellipse with parametrization
\[
   \frac{1}{\sqrt{\log(1+1/t)}}e^{i\theta}+|\gamma| \sqrt{\log(1+1/t)} e^{i(\alpha-\theta)}, \qquad\theta\in[0,2\pi]
\]
where $\alpha=\arg \gamma$.
In fact, the Brown measure of $T+g_{t,\gamma}$ is the uniform measure on its support. 

\subsection{The Haar unitary}
\label{section:Brown-u-gt}

The Haar unitary operator $u$ is $R$-diagonal. Hence, results in Section \ref{section:addition-R-diag} applies to $u+g_{t,\gamma}$. 
The Brown measure of $u+c_t$ is $R$-diagonal and the density formula is determined by the $S$-transform of $(u+c_t)^*(u+c_t)$ \cite{HaagerupLarsen2000, HaagerupSchultz2007}.
However, it is not easy to obtain a precise formula for the Brown measure of $u+c_t$ by calculating the $S$-transform. See \cite[Section 3.1]{HKS-2010-tams} for a relevant calculation on the Cauchy transform of $c+\lambda u$ where $c$ is the standard circular operator with variance one and $\lambda\in \mathbb{C}\backslash\{0\}$. Hence, we do not apply Proposition \ref{prop:Phi-t-gamma-formula-Rdiag} directly.  We calculate the push-forward map from the Brown measure of $u+c_t$ to $u+g_{t,\gamma}$ using formulas in Section \ref{section:Brown-addition-elliptic} and as a by product we are able to calculate the Brown measure of $u+c_t$.

For $t>0$, following Definition \ref{defn:s-epsilon-t-00}, we would like to determine the choice of $s$ so that 
\[
\int_0^\infty\frac{1}{s^2+x^2}d\mu_{|u-\lambda\unit|}(x)=\frac{1}{t}.
\] 
To this end, we write down the following elementary results for convenience. 
Given $\lambda_0\neq 0$, 
consider the quadratic equation
\begin{equation}\label{eqn:quadratic-1}
(z-\lambda_0)(1-z\overline{\lambda_0})+s^2 z=0.
\end{equation}
Rewrite it as
\[
\overline{\lambda_0} z^2- (|\lambda_0|^2+1+s^2)z+\lambda_0=0.
\]
We set
\begin{equation}\label{eqn:delta}
\delta=|\lambda_0|^2+1+s^2. 
\end{equation}
Then \eqref{eqn:quadratic-1} has two solutions
\begin{align}
Z_1=\frac{\delta-\sqrt{\delta^2-4|\lambda_0|^2}}{2\overline{\lambda_0}},\quad
Z_2=\frac{\delta+\sqrt{\delta^2-4|\lambda_0|^2}}{2\overline{\lambda_0}}.\label{eqn:Z2}
\end{align}
Note that $Z_1Z_2=\lambda_0/\overline{\lambda_0}$, $|Z_1 Z_2|=1$,
and
\begin{equation}\label{eqn:delta-ineq}
\delta^2-4|\lambda_0|^2\geq (|\lambda_0|^2+1)^2-4|\lambda_0|^2=(|\lambda_0|^2-1)^2\geq 0
\end{equation}
for any $s$. Hence, for any $\lambda_0\neq 0$, we have $|Z_1|\leq 1$ and $|Z_2|\geq 1$. Moreover, 
when $s>0$, we have $|Z_1|<1$ and $|Z_2|>1$.

\begin{lemma}
	\label{lemma:7.3}
	For any $\lambda_0$, we have 
	\begin{equation}\label{eqn:p0-integration-0}
	\frac{1}{2\pi}\int_0^{2\pi}\frac{1}{|e^{i\theta}-\lambda_0|^2+s^2}d\theta=
	\frac{1}{\sqrt{\delta^2-4|\lambda_0|^2}}
	\end{equation}
	where $\delta=|\lambda_0|^2+1+s^2$ is given by (\ref{eqn:delta})
	and the integration is infinity if $\delta^2=4|\lambda_0|^2$. 
	In particular,  when $s=0$, 
	\begin{equation}\label{eqn:p0-integration-0-zero}
	\frac{1}{2\pi}\int_0^{2\pi}\frac{1}{|e^{i\theta}-\lambda_0|^2}d\theta=
	\begin{cases}
	\frac{1}{1-|\lambda_0|^2}\quad \text{if} \quad 0\leq |\lambda_0|<1,\\
	\frac{1}{|\lambda_0|^2-1}\quad \text{if} \quad |\lambda_0|>1,\\
	\infty \quad \text{if} \quad |\lambda_0|=1.
	\end{cases}
	\end{equation}
\end{lemma}
\begin{proof}
	We note that
	\[
	\frac{1}{|e^{i\theta}-\lambda_0|^2+s^2}=\frac{e^{i\theta}}{(e^{i\theta}-\lambda_0)(1-e^{i\theta}\overline{\lambda_0})+s^2e^{i\theta}}.
	\]
	Using the formula \eqref{eqn:Z2} for roots of the quadratic equation \eqref{eqn:quadratic-1}, we have 
	\begin{align*}
	\frac{1}{2\pi}\int_0^{2\pi}\frac{1}{|e^{i\theta}-\lambda_0|^2+s^2}d\theta&
	=\frac{1}{2\pi i}\int_{|z|=1}\frac{1}{(z-\lambda_0)(1-z\overline{\lambda_0})+s^2 z}dz\\
	&=\frac{1}{2\pi i}\int_{|z|=1}\frac{1}{-\overline{\lambda_0}(z-Z_1)(z-Z_2)}dz\\
	&=\frac{-\lambda_0}{|\lambda_0|^2(Z_1-Z_2)}\\
	&= \frac{1}{\sqrt{\delta^2-4|\lambda_0|^2}}
	\end{align*}
	provided that $\delta^2-4|\lambda_0|^2\neq 0$. When $s=0$ and $|\lambda_0|\neq 1$, the above calculation still works. 
	
	For the case that $\delta^2-4|\lambda_0|^2=0$. By using \eqref{eqn:delta-ineq}, we see that  $\delta^2-4|\lambda_0|^2=0$ if and only if $|\lambda_0|=1$ and $s=0$. It is clear the integration is infinity in this case.
\end{proof}

\begin{lemma}
	\label{lemma:7.4}
	For any $\lambda_0$, we have 
	\begin{equation}\label{eqn:p0-theta-1}
	\frac{1}{2\pi}\int_0^{2\pi}\frac{e^{i\theta}}{|e^{i\theta}-\lambda_0|^2+s^2}d\theta=
	\begin{cases}
	\frac{\lambda_0}{\sqrt{\delta^2-4|\lambda_0|^2}}\frac{\delta-\sqrt{\delta^2-4|\lambda_0|^2}}{2|\lambda_0|^2}, &\text{if}\quad \lambda_0\neq 0\\
	0, &\text{if}\quad \lambda_0=0.
	\end{cases}
	\end{equation}
	where $\delta=|\lambda_0|^2+1+s^2$ is given by (\ref{eqn:delta})
	and the integration is infinity if the denominator is zero. In particular, 
	when $u=0$.
	\begin{equation}\label{eqn:p0-theta-1-zero}
	\frac{1}{2\pi}\int_0^{2\pi}\frac{e^{i\theta}}{ |e^{i\theta}-\lambda_0|^2 }d\theta=
	\begin{cases}
	\frac{\lambda_0}{1-|\lambda_0|^2}\quad \text{if} \quad |\lambda_0|<1\\
	\frac{\lambda_0}{|\lambda_0|^2(|\lambda_0|^2-1)}\quad \text{if} \quad |\lambda_0|>1\\
	\infty \quad \text{if} \quad |\lambda_0|=1.
	\end{cases}
	\end{equation}
	
	Similarly, 
	\begin{equation}\label{eqn:p0-theta-2}
	\frac{1}{2\pi}\int_0^{2\pi}\frac{e^{-i\theta}}{|e^{i\theta}-\lambda_0|^2+s^2}d\theta=
	\begin{cases}
	\frac{\overline{\lambda_0}}{\sqrt{\delta^2-4|\lambda_0|^2}}\frac{\delta-\sqrt{\delta^2-4|\lambda_0|^2}}{2|\lambda_0|^2}, &\text{if}\quad \lambda_0\neq 0\\
	0, &\text{if}\quad \lambda_0=0.
	\end{cases}
	\end{equation}
	where $\delta=|\lambda_0|^2+1+s^2$ is given by (\ref{eqn:delta})
	and the integration is infinity if the denominator is zero. In particular, 
	when $u=0$.
	\begin{equation}\label{eqn:p0-theta-2-zero}
	\frac{1}{2\pi}\int_0^{2\pi}\frac{e^{-i\theta}}{ |e^{i\theta}-\lambda_0|^2 }d\theta=
	\begin{cases}
	\frac{\overline{\lambda_0}}{1-|\lambda_0|^2}\quad \text{if} \quad |\lambda_0|<1\\
	\frac{\overline{\lambda_0}}{|\lambda_0|^2(|\lambda_0|^2-1)}\quad \text{if} \quad |\lambda_0|>1\\
	\infty \quad \text{if} \quad |\lambda_0|=1.
	\end{cases}
	\end{equation}
\end{lemma}
\begin{proof}
	Recall that $|Z_1|<1$ and $|Z_2|>1$. 
	Using the formula \eqref{eqn:Z2} for roots of the quadratic equation \eqref{eqn:quadratic-1}, we have
	\begin{align*}
	\frac{1}{2\pi}\int_0^{2\pi}\frac{e^{i\theta}}{|e^{i\theta}-\lambda_0|^2+s^2}d\theta&
	=\frac{1}{2\pi i}\int_{|z|=1}\frac{z}{(z-\lambda_0)(1-z\overline{\lambda_0})+s^2 z}dz\\
	&=\frac{1}{2\pi i}\int_{|z|=1}\frac{z}{-\overline{\lambda_0}(z-Z_1)(z-Z_2)}dz\\
	&=\frac{-\lambda_0 Z_1}{|\lambda_0|^2(Z_1-Z_2)} \\
	&= \frac{\lambda_0}{\sqrt{\delta^2-4|\lambda_0|^2}}\frac{\delta-\sqrt{\delta^2-4|\lambda_0|^2}}{2|\lambda_0|^2}
	\end{align*}
	provided that $\delta^2-4|\lambda_0|^2\neq 0$ and $\lambda_0\neq 0$. If $\lambda_0=0$, the integration is clearly equal to zero. Note that $\delta^2-4|\lambda_0|^2= 0$ if and only if $|\lambda_0|=1$ and $s=0$ (see \eqref{eqn:delta-ineq}). The integration is infinity in this case. Hence \eqref{eqn:p0-theta-1} and \eqref{eqn:p0-theta-1-zero} are proved.
	
	Finally, when $s=0$, we have $\delta^2-4|\lambda_0|^2=(|\lambda_0|^2-1)^2$. Plugging this into \eqref{eqn:p0-theta-1}, we obtain \eqref{eqn:p0-integration-0-zero}.  
	By taking conjugation, we obtain \eqref{eqn:p0-theta-2} and \eqref{eqn:p0-theta-2-zero}. 
\end{proof}

	Hence, we have 
\begin{equation}
 \label{eqn:p-lambda-for-Haar}
	\begin{aligned}
		&p_\lambda^{(0)}(s)=\phi\bigg( (\lambda-u)^*[(\lambda-u)^*(\lambda-u)+s^2]^{-1} \bigg)\\
	=&\frac{1}{2\pi}\int_0^{2\pi} \frac{\overline{\lambda}-e^{-i\theta}}{|e^{i\theta}-\lambda|^2+s^2}d\theta\\
	=&\frac{\overline{\lambda}}{\sqrt{\delta^2-4|\lambda|^2}}
	-\frac{\overline{\lambda}}{\sqrt{\delta^2-4|\lambda|^2}}\frac{\delta-\sqrt{\delta^2-4|\lambda|^2}}{2|\lambda|^2},
	\end{aligned}
\end{equation}
		where $\delta=|\lambda|^2+1+s^2$.

A direct application of \eqref{eqn:p0-integration-0-zero} yields that 
\begin{equation}
\label{eqn:Xi-t-for-u-ct}
\Xi_t=\left\{ \lambda: \phi \left[\big( (u-\lambda\unit)^*(u-\lambda\unit) \big)^{-1}\right]>\frac{1}{t}  \right\}
  =\{ \lambda:   \lambda_1<|\lambda|<\lambda_2 \}
\end{equation}
where $\lambda_1=\sqrt{(1-t)_+}$ and $\lambda_2=\sqrt{1+t}$. 
\begin{proposition}
The support of the Brown measure of $u+c_t$ is the single ring
\[
   \overline{\Xi_t}=\{ \lambda:   \lambda_1\leq|\lambda|\leq\lambda_2 \}.
\]
Moreover, the Brown measure is absolutely continuous  and the density of $u+c_t$ is strictly positive on its support. 
\end{proposition}
\begin{proof}
	By \eqref{eqn:Xi-t-for-u-ct} and Theorem \ref{thm:support-x0-c_t-general}, it follows that the support of the Brown measure of $u+c_t$ is $\overline{\Xi_t}$. One can also compare it with the general result for the Brown measure of $R$-diagonal operators \cite{HaagerupLarsen2000}. 
\end{proof}

We can now calculate the subordination function $w(0;\lambda, t)$ in this context. For any $\lambda\in\Xi_t$, following Definition \ref{defn:s-epsilon-t-00}, $w(0;\lambda, t)$ is determined by the condition
\[
   \phi \left[\big( (u-\lambda\unit)^*(u-\lambda\unit) +w(0;\lambda, t)^2\big)^{-1}\right]=
    \frac{1}{t}.
\]
By \eqref{eqn:p0-integration-0}, this condition is equivalent to 
\[
  \delta^2-4|\lambda|^2=  ( |\lambda|^2+1+w(0;\lambda, t)^2)^2-4|\lambda|^2=t^2.
\]
It follows that
\begin{equation}
  w(0;\lambda, t)^2=\sqrt{4|\lambda|^2+t^2}-(|\lambda|^2+1), \qquad  \lambda_1<|\lambda|<\lambda_2.
\end{equation}

We now describe the map $\Phi_{t,\gamma}$ more closely. For $\lambda\in \Xi_t$, we set
$s=|\lambda|^2-1$,
where $\lambda_1^2-1<s<\lambda_2^2-1$. 
Let $\delta=|\lambda|^2+1+w(0;\lambda,t)^2$. 
Since $\delta^2-4|\lambda|^2=t^2$, we now have 
\begin{align*}
	\frac{\delta-\sqrt{\delta^2-4|\lambda|^2}}{2|\lambda|^2}
	    =\frac{\delta-t}{2|\lambda|^2}
	    =\frac{\sqrt{4(1+s)+t^2}-t}{2(1+s)}.
\end{align*}
Therefore, for $\lambda_1^2-1<s<\lambda_2^2-1$ and $\lambda$ such that $|\lambda|^2=s+1$, by \eqref{eqn:p-lambda-for-Haar}, we have, 
\begin{align}
\label{eqn:Phi-HaarU-gt}
 \Phi_{t,\gamma}(\lambda)&=\lambda+\gamma p_\lambda^{(0)}(w(0;\lambda,t))\nonumber\\
   &=\lambda+ \gamma\left( \frac{\overline{\lambda}}{\sqrt{\delta^2-4|\lambda|^2}}
    -\frac{\overline{\lambda}}{\sqrt{\delta^2-4|\lambda|^2}}\frac{\delta-\sqrt{\delta^2-4|\lambda|^2}}{2|\lambda|^2}
\right)
\end{align}
 where $\delta=|\lambda|^2+1+w(0;\lambda,t)^2$. 
 
 We first look at the case when $\gamma=t$.
 Since $\delta^2-4|\lambda|^2=t^2$, for $\lambda= r e^{i\theta}$ where $r=\sqrt{s+1}$, $\Phi_{t,t}$
 as in \eqref{eqn:Phi-HaarU-gt} can be further simplified as 
 \begin{align}
 \label{eqn:J-t-HaarU}
  \Phi_{t,t}( r e^{i\theta})
  &=r e^{i\theta}+t re^{-i\theta}\left( \frac{1}{t} +\frac{1}{2r^2}- \frac{\sqrt{4r^2+t^2}}{2r^2t} \right)
   \\ \label{eqn:8.17}
  &=r e^{i\theta}+r e^{-i\theta}\left( 1- \frac{\sqrt{4(1+s)+t^2}-t}{2(1+s)} \right)\nonumber\\
    &=a(s)\cos (\theta)+ib(s) \sin(\theta)
 \end{align}
 where, $\lambda_1^2-1<s<\lambda_2^2-1$, and
 \begin{align}
 \label{eqn:a-b-t-J-t-HaarU}
  \begin{cases}
   a(s)=\left(2\sqrt{1+s}- \frac{\sqrt{4(1+s)+t^2}-t}{2\sqrt{1+s}}\right);\\
   b(s)=\frac{\sqrt{4(1+s)+t^2}-t}{2\sqrt{1+s}}.
  \end{cases}
 \end{align}
In particular, if $\lambda_1^2-1=s$, we have
$1+s=\lambda_1^2=(1-t)_+$ and 
 $4(1+s)+t^2=(4(1-t)_+)+t^2$. Then for $s=((1-t)_+)-1$, 
 \[
  \begin{cases}
   a(s)= a^-=\sqrt{(1-t)_+}\\
   b(s)= b^{-}=\sqrt{(1-t)_+}
  \end{cases}.
\]
Similarly, when $\lambda_2^2-1=s$, we have $s=t$ and
\[
   \begin{cases}
   a(s)= a^+=\frac{2t+1}{\sqrt{t+1}}\\
   b(s)=b^+=\frac{1}{\sqrt{1+t}}\\
  \end{cases}.
\]

\begin{proposition}
\label{prop:HaarU-a-b-s-monotonicity}
For $\lambda_1^2-1<s<\lambda_2^2-1$, the circle $C_s:=\{\lambda: |\lambda|^2=1+s \}$ is mapped to the ellipse centered at the origin with semi-axes $a(t)$ and $b(t)$ as in \eqref{eqn:a-b-t-J-t-HaarU}. Moreover, the function $s\mapsto a(s)$ is a strictly increasing function from $(\lambda_1^2-1, \lambda_2^2-1)$ onto $(a^-, a^+)$; and the 
 the function $s\mapsto b(s)$ is a strictly increasing function from $(\lambda_1^2-1, \lambda_2^2-1)$ onto $(b^{-}, b^+)$.
\end{proposition}
\begin{proof}
We have 
\[
  a'(s)=\frac{\frac{t^2}{\sqrt{4 s+t^2+4}}+4 s-t+4}{4 (s+1)^{3/2}},
  \qquad \text{and} \qquad 
   b'(s)=\frac{t \left(\sqrt{4 s+t^2+4}-t\right)}{4 (s+1)^{3/2} \sqrt{4 s+t^2+4}}.
\]
Hence $b'(s)>0$. To show $a'(s)>0$, we set
\[
  f(s)={\frac{t^2}{\sqrt{4 s+t^2+4}}+4 s-t+4}
\]
and calculate its derivatives 
\[
  f'(s)=4-\frac{2 t^2}{\left(4 s+t^2+4\right)^{3/2}}
\]
and 
\[
  f''(s)=\frac{12 t^2}{\left(4 s+t^2+4\right)^{5/2}}>0.
\]
Hence $f'(s)>f'(\lambda_1^2-1)=f'((1-t)_+)>0$. Consequently, $f$ is increasing. 
We then check that 
\[
  f((1-t)_+)>0.
\]
It follows that $a'(s)>0$ for $\lambda_1^2-1<s<\lambda_2^2-1$.
\end{proof}

\begin{proposition}
	\label{prop:Phi-non-sigular-Haar}
The map $\Phi_{t,t}$ is non-singular at any $\lambda\in \{ z:   \lambda_1<|z|<\lambda_2 \}$. 
\end{proposition}
\begin{proof}
	Choose the coordinate $(s,\theta)$ so that $x=r\cos\theta$ and $y=r\sin\theta$, where $r^2=s+1$ and $\lambda_1^2-1<s<\lambda_2^2-1$. 
The map $(s,\theta)\mapsto \Phi_{t,t}(re^{i\theta})=a(s)\cos (\theta)+ib(s) \sin(\theta)$ has the Jacobian given by
\[
 D\Phi_{t,t}(s,\theta)=\begin{bmatrix}
   a'(s)\cos\theta &-a(s)\sin\theta\\
   b'(s)\sin\theta & b(s)\cos\theta
 \end{bmatrix}.
\]
By Proposition \ref{prop:HaarU-a-b-s-monotonicity}, we have $a'(s)>0, b'(s)>0$, we hence have 
$\det(D\Phi_{t,t})>0$. 
\end{proof}

 \begin{theorem}
 \label{thm:Brown-HaarUnitary-c-t}
 The Brown measure of $u+c_t$ is given by
 \begin{equation}
   \mu_{u+c_t}(\{\lambda: |\lambda|\leq r)\})=
    r^2\left( \frac{1}{t} +\frac{1}{2r^2}- \frac{\sqrt{4r^2+t^2}}{2r^2t} \right),
 \end{equation}
 where $\sqrt{(1-t)_+}\leq r\leq \sqrt{1+t}$.
 \end{theorem}
 \begin{proof}
	By Lemma \ref{lemma:formula-derivative-t-0} (\cite[Corollary 4.4]{Zhong2021RdiagBrown}), for $\sqrt{(1-t)_+}\leq |\lambda|\leq \sqrt{1+t}$,
	\begin{equation*}
	\begin{aligned}
	p_\lambda^{c,(t)}(0)&=\frac{\partial S}{\partial \lambda}(u+c_t,\lambda,0)\\
	&=\frac{\partial}{\partial\lambda} \log\Delta((u+c_t-\lambda\unit)^*(u+c_t-\lambda\unit))\\
	&=\frac{\mu_{u+c_t}(\{\lambda: |\lambda|\leq r)\})}{\lambda},
	\end{aligned}
	\end{equation*}
	which yields $\Phi_{t,t}(\lambda)=\lambda+\frac{t}{\lambda}  \mu_{u+c_t}(\{\lambda: |\lambda|\leq r)\})$. Hence, the result follows by comparing this general formula with the explicit formula of $\Phi_{t,t}(\lambda)$ as in \eqref{eqn:8.17}.
\end{proof}

\begin{theorem}
\label{thm:Brown-HaarUnitary-x-t-gamma}
The support of the Brown measure of $u+g_t$ is the deformed single ring where the inner boundary is the circle centered at the origin with radius $\sqrt{(1-t)_+}$ and the outer boundary is the ellipse centered at the origin with semi-axes 
$\frac{2t+1}{\sqrt{t+1}}$ and $\frac{1}{\sqrt{1+t}}$.
The Brown measure is absolutely continuous and its density is strictly positive in the support. 

 Moreover, the Brown measure of $u+g_t$ is the push-forward map of the Brown measure of $u+c_t$ under the map $\Phi_{t,t}$ defined as
 \[
   \Phi_{t,t}(re^{i\theta})=a(s)\cos (\theta)+ib(s) \sin(\theta),\qquad r=\sqrt{s+1},
 \]
 where $ ((1-t)_+)-1<s<t$ and $a(s), b(s)$ are given by \eqref{eqn:a-b-t-J-t-HaarU}.

\end{theorem}
\begin{proof}
It follows from Proposition \ref{prop:Phi-non-sigular-Haar} that the map $\Phi_{t,t}$ is, non-singular, one-to-one and onto from the single ring  $\{\lambda: \lambda_1<|\lambda|<\lambda_2 \}$, the support of the Brown measure of $u+c_t$, to the deformed ring. Theorem \ref{thm:main-push-forward-general} applies and the push-forward map is given by \eqref{eqn:J-t-HaarU} and \eqref{eqn:a-b-t-J-t-HaarU}. Since the Brown measure of $u+c_t$ is strictly positive on its support, the result then follows. 
\end{proof}

We now fix $\gamma=|\gamma|e^{2i\alpha}$ such that $|\gamma|\leq t$ and $\lambda_1^2-1<s<\lambda_2^2-1$ and $r^2=|\lambda|^2=s+1$,  then similar to the case $\gamma=t$ as in \eqref{eqn:Phi-HaarU-gt}, we have
 \begin{align}
 \label{eqn:Phi-t-gamma-HaarU}
  \Phi_{t,\gamma}( r e^{i\theta})&=r e^{i\theta}+r e^{i(2\alpha-\theta)}\frac{|\gamma|}{t}\left( 1- \frac{\sqrt{4(1+s)+t^2}-t}{2(1+s)} \right)\nonumber\\
   &=e^{i\alpha}\left[ \sqrt{s+1}e^{i(\theta-\alpha)}+\frac{|\gamma|}{t} e^{i(\alpha-\theta)}\left( \sqrt{s+1}- \frac{\sqrt{4(1+s)+t^2}-t}{2\sqrt{1+s}} \right)\right],
 \end{align}
which is an ellipse twisted by angle $\alpha$ with semi-axes
\begin{align}
\label{eqn:Haar-x-t-gamma-semi-axes}
 \begin{cases}
    f(s)&=\frac{t-|\gamma|}{t}\sqrt{s+1}+\frac{|\gamma|}{t} a(s),\\
   g(s)&=\frac{t-|\gamma|}{t}\sqrt{s+1}+\frac{|\gamma|}{t}b(s),
 \end{cases}
\end{align}
where $a(s), b(s)$ are given by \eqref{eqn:a-b-t-J-t-HaarU}. It follows by Proposition \ref{prop:HaarU-a-b-s-monotonicity} that both $f$ and $g$ are increasing function of $s$ in the interval $(\lambda_1^2-1, \lambda_2^2-1)$. 

Then we deduce the following result in the same way as Theorem \ref{thm:Brown-HaarUnitary-x-t-gamma}, as a special case of Theorem \ref{thm:Brown-addition-R-diag-elliptic-pushforward}. We leave details for interested readers. 
\begin{theorem}
\label{thm:Brown-Haar-x-general}
The support of the Brown measure of $u+g_{t,\gamma}$ is the deformed single ring where the inner boundary is the circle centered at the origin with radius $\sqrt{(1-t)_+}$, and the outer boundary is the ellipse twisted by angle $\alpha$, centered at the origin, with semi-axes $\frac{t+|\gamma|+1}{\sqrt{t+1}}$ and $\frac{t-|\gamma|+1}{\sqrt{t+1}}$.

Moreover, the Brown measure of $u+g_{t,\gamma}$ is the push-forward map of the Brown measure of $u+c_t$ under the map $\Phi_{t,\gamma}$ defined as
\[
  \Phi_{t,\gamma}(re^{i\theta})=e^{i\alpha} \left( f(s)\cos(\theta-\alpha)+i g(s) \sin(\theta-\alpha) \right),
  \qquad r=\sqrt{s+1} 
\]
where $((1-t)_+)-1<s<t$, $f(s), g(s)$ are defined by \eqref{eqn:Haar-x-t-gamma-semi-axes}.
\end{theorem}

%


\bigskip

\bigskip
\noindent
{\bf Acknowledgment.} 
The author would like to express his deep gratitude to Hari Bercovici and Zhi Yin for many valuable discussions and for careful reading of earlier versions of this paper. The author wants to thank Serban Belinschi for his comments on Hermitian reduction method and subordination functions. Special thanks to Brian Hall for encouraging the author to develop the results in Section \ref{section:Phit_new_formula} and for his many insightful suggestions.
The author also wants to thank Ken Dykema, Ching-Wei Ho, Dang-Zheng Liu, Weihua Liu, Alexandru Nica, Zhuang Niu, Edward Saff, Piotr  \'{S}niady, and Alex Solynin for discussions on the topics of Brown measures and potential theory during the course of the investigation. 
The author is grateful to the anonymous reviewers for their thorough reading of the manuscript and their valuable insights and suggestions.

\appendix
\section{The distributional derivative of the logarithmic potential}
 \label{appendix-A}
We establish some elementary results about computing the distributional derivative
of the logarithmic potential of a regularized Brown measure. The results in this Appendix are useful in Section \ref{section:Phit_new_formula}.

\begin{lemma}
\label{parts.thm}
Suppose $f$ is a locally integrable function on
$\mathbb{R}^{n},$ with $n\geq1.$ Suppose also that $f$ is continuously
differentiable on $\mathbb{R}^{n}\setminus\{0\}$ and that $\partial f/\partial
x_{j}$ is locally integrable on $\mathbb{R}^{n}$ (not just on $\mathbb{R}%
^{n}\setminus\{0\}\,$). Finally, suppose that
\begin{equation}
\lim_{\mathbf{x}\rightarrow0}\left\vert \mathbf{x}\right\vert ^{n-1}%
f(\mathbf{x})=0.\label{condition}%
\end{equation}
Then $\partial f/\partial x_{j}$ in the distribution sense is computed by
integrating a test function against the classical partial derivative, and the
one point where the classical derivative is not defined. In the case $n=1,$
the condition (\ref{condition}) can be replaced by the condition that
$\lim_{x\rightarrow0}f(x)$ exists. 
\end{lemma}

\begin{proof}
The distributional derivative is computed by looking at the integral%
\[
-\int_{\mathbb{R}^{n}}f(\mathbf{x})\frac{\partial\chi}{\partial x_{j}%
}(\mathbf{x})~d^{n}\mathbf{x},
\]
where $\chi$ is a test function. Since $f$ is locally integrable, the above
integral exists and can be computed as
\[
-\lim_{\varepsilon\rightarrow0}\int_{\mathbb{R}^{n}\setminus D_{\varepsilon}%
}f(\mathbf{x})\frac{\partial\chi}{\partial x_{j}}(\mathbf{x})~d^{n}\mathbf{x},
\]
where $D_{\varepsilon}$ is the ball of radius $\varepsilon$ centered at the
origin. Once a ball around the origin has been cut out, everything is $C^{1}$
and we can integrate by parts to get
\begin{equation}
-\int_{\mathbb{R}^{n}\setminus D_{\varepsilon}}f(\mathbf{x})\frac{\partial
\chi}{\partial x_{j}}(\mathbf{x})~d^{n}\mathbf{x=-}\int_{S_{\varepsilon}%
}f(\mathbf{x})\chi(\mathbf{x)~}v_{j}~d\Gamma\mathbf{+}\int_{\mathbb{R}%
^{n}\setminus D_{\varepsilon}}\frac{\partial f}{\partial x_{j}}(\mathbf{x}%
)\chi(\mathbf{x})~d^{n}\mathbf{x},\label{parts}%
\end{equation}
where $S_{\varepsilon}$ is the sphere of radius $\varepsilon,$ $v_{j}$ is the
$j$th component of the outward (with respect to $\mathbb{R}^{n}\setminus
D_{\varepsilon}$) unit normal, and $\Gamma$ is the surface area measure on the
sphere $S_{\varepsilon}.$ 

Now, $\left\vert v_{j}\right\vert \leq1$ and $\chi$ is bounded, so the first
(boundary) term on the right-hand side of (\ref{parts}) can be bounded using
the area of the sphere, which is proportional to $\varepsilon^{n-1}$:%
\[
\left\vert \int_{S_{\varepsilon}}f(\mathbf{x})\chi(\mathbf{x)~}v_{j}%
~d\Gamma\right\vert \leq C\varepsilon^{n-1}\sup_{S_{\varepsilon}}\left\vert
f(\mathbf{x})\right\vert ,
\]
where the right-hand side of the above inequality tends to zero by the
assumption (\ref{condition}). Letting $\varepsilon\rightarrow0$ in
(\ref{parts}) using the local integrability of $\partial f/\partial x_{j}$
therefore gives%
\[
-\int_{\mathbb{R}^{n}}f(\mathbf{x})\frac{\partial\chi}{\partial x_{j}%
}(\mathbf{x})~d^{n}\mathbf{x=}\int_{\mathbb{R}^{n}}\frac{\partial f}{\partial
x_{j}}(\mathbf{x})\chi(\mathbf{x})~d^{n}\mathbf{x},
\]
which is what we wanted to show. 

In the case $n=1,$ if $\lim_{x\rightarrow0}f(x)$ exists, we can add a constant
to $f$ (without changing the distributional derivative) so that the limit is
zero and (\ref{condition}) will hold with $n=1.$ 
\end{proof}

\begin{example}
If $n=2$ and $f(\mathbf{x})=\log\left\vert \mathbf{x}\right\vert ,$ then
(\ref{condition}) holds. Furthermore,
\[
\frac{\partial f}{\partial x_{j}}=\frac{1}{\left\vert \mathbf{x}\right\vert
}\frac{\partial}{\partial x_{j}}\left\vert \mathbf{x}\right\vert =\frac
{1}{\left\vert \mathbf{x}\right\vert }\frac{x_{j}}{\left\vert \mathbf{x}%
\right\vert },
\]
which is a locally integrable near the origin. 
\end{example}

We now turn to the logarithmic potential $S_{\mu}$ and Cauchy transform $G_{\mu}$ of a
compactly supported probability measure on the complex plane $\mathbb{C}$, defined as%
\begin{align*}
S_{\mu}(z)  & =\int_{\mathbb{C}}\log(\left\vert z-w\right\vert ^{2})~d\mu(w),\\
G_{\mu}(z)  & =\lim_{\varepsilon\downarrow 0}\int_{\vert z-w\vert\geq \varepsilon}\frac{1}{z-w}~d\mu(w).
\end{align*}
We will work on probability measures $\mu$ satisfying
\begin{equation} 
  \label{eqn-condition.A3}
  \mu\left( \{ w: \vert z-w\vert \leq r \} \right)\leq r^s
\end{equation}
for some $s>1$. 
By an elementary argument using Fubini's theorem and the local
integrability of the functions $\log\left(  \left\vert z\right\vert
^{2}\right)  $ and $1/z$, we can show (1) that the integrals defining
$S_{\mu}(z)$ and $G_{\mu}(z)$ are finite at every $z$, 
and (2) that $S_{\mu}$ and $G_{\mu}$ are locally integrable functions on the complex plane.

Recall that the regularized Brown measure $\mu_{x,\varepsilon}$ of $x\in\mathcal{M}$ is determined by
  \[
      \int_\mathbb{C}\log\vert \lambda-z\vert d\mu_{x,\varepsilon}(z)=
      \frac{1}{2}\log\Delta(\vert x-\lambda\vert^2+\varepsilon^2).
  \]
Then, by \cite[Lemma 2.8]{HaagerupSchultz2007}, the regularized measure $\mu_{x,\varepsilon}$ is absolutely continuous, and its density (up to a constant $1/4\pi$) is 
\begin{align}
\label{eqn_4_order}
 &\frac{\partial^2}{\partial_{\bar\lambda}\partial_\lambda}\log\Delta(\vert x-\lambda\vert^2+\varepsilon^2)\\
   &\qquad=\varepsilon^2\phi\Huge( (\vert x-\lambda\vert^2+\varepsilon^2)^{-1} (\vert (x-\lambda)^*\vert^2+\varepsilon^2)^{-1}\Huge)\nonumber
   \leq 4 \varepsilon^2/|\lambda|^4,  
\end{align}
provided that $|\lambda|>2\Vert x\Vert$.

\begin{lemma}
 \label{lemma:bounded1}
For any $R>0$ and $\varepsilon>0$, we have 
\[
      \int_\mathbb{C}\int_{|z|\leq R}\vert \log|w-z| \vert d^2z d\mu_{x,\varepsilon}(w)<\infty. 
   \]
\end{lemma}
\begin{proof}
 We denote 
 \[
   I_1(w)=\int_{|z|\leq R}\vert \log(|w-z|)\vert d^2z
 \]
 and change the variable $u=w-z$ to get 
 \[
   I_1(w)=\int_{|u-w|\leq R}\vert \log(|u|)\vert d^2u.
 \]
We then consider two cases: $|w|\leq R+1$ and $|w|>R+1$. 

If $|w|\leq R+1$, we note that $|u|\leq |u-w|+|w|\leq 2R+1$. Thus, we have $\{u: |u-w|\leq R\}\subset \{u: |u|\leq 2R+1\}$ and in this case 
\[
  I_1(w)\leq C:=\int_{|u|\leq 2R+1}\vert \log(|u|) \vert d^2u.
\]
If $|w|>R+1$, then $|u|\geq |w|-|u-w|\geq 1$ when $|u-w|\leq R$. Thus, $\log (|u|)\geq 0$ when $|u-w|\leq R$. In addition, note that $|u|\leq |u-w|+|w|\leq R+|w|$, we then get
\[
  I_1(w)\leq \pi R^2\log(R+|w|).
\]

Therefore, we have
\[
 I_1(w)\leq \begin{cases}
      C    &\text{for} \quad |w|\leq R+1,\\
      \pi R^2\log(R+|w|) &\text{for} \quad |w|>R+1.
   \end{cases}
\]
By \eqref{eqn_4_order}, the density of $d\mu_{x+c_t,\varepsilon}$ decays of the order $1/|w|^4$ for large $|w|$. It follows that $\int_\mathbb{C}I_1(w) d\mu(w)<\infty$. 
\end{proof}

\begin{corollary}\label{cor3_distderivative}
If $\mu_{x,\varepsilon}$ is the regularized Brown measure of $x\in\mathcal{M}$ with parameter $\varepsilon>0$, then the
distributional $z$-derivative of the log potential is the Cauchy transform. 
\end{corollary}
\begin{proof} 
 Suppose that $\chi$ is a test function supported in the ball $\{z: |z|\leq R\}$. In addition, we may suppose that $|\chi(z)|\leq 1$ and 
$\left|\frac{\partial\chi(z)}{\partial z}\right |\leq 1$. 
We denote
\[
   I_2(w)=\int_{|z|\leq R}\frac{1}{|z-w|}d^2z.
\]
Note that $I_2(0)=2\pi R$.
Then we have 
\[
 I_2(w)\leq \begin{cases}
    2\pi & \text{for}\quad |z-w|\leq 1,\\
    \pi R^2 & \text{for} \quad |z-w|>1.
  \end{cases}
\]
It follows that
\begin{equation}
  \label{eqn:0020}
 \int_\mathbb{C}\int_{|z|\leq R}\frac{1}{|z-w|}d^2zd\mu_{x,\varepsilon}(w)<\infty.
\end{equation}

Lemma \ref{lemma:bounded1} and Equation \eqref{eqn:0020} guarantee that we can apply Fubini's theorem in the following calculation:
\begin{align*}
& -\int_{{|z|\leq R}}\frac{\partial\chi}{\partial z}(z)\int_{\mathbb{C}}\log(\left\vert
z-w\right\vert ^{2})~d\mu(w)~d^{2}z\\
& =-\int_{\mathbb{C}}\int_{{|z|\leq R}}\frac{\partial\chi}{\partial z}(z)\log(\left\vert
z-w\right\vert ^{2})~d^{2}z~d\mu(w)\\
& =\int_{\mathbb{C}}\int_{{|z|\leq R}}\chi(z)\frac{\partial}{\partial z}\log(\left\vert
z-w\right\vert ^{2})~d^{2}z~d\mu(w)\\
& =\int_{\mathbb{C}}\int_{|z|\leq R}\chi(z)\frac{1}{z-w}~d^{2}z~d\mu(w)\\
& =\int_{|z|\leq R}\chi(z)\int_{\mathbb{C}}\frac{1}{z-w}~d\mu(w)~d^{2}z,
\end{align*}
where we used Lemma \ref{parts.thm} to derive the second identity. 
This finishes the verification.
\end{proof}

We next study the continuity of the Cauchy transform of a probability measure on $\mathbb{C}$ under some regularity conditions. 

\begin{proposition}\label{lemma:appendix_continuityCauchy}
Given a probability measure $\mu$ on the complex plane $\mathbb{C}$, suppose that $\mu$ is absolutely continuous with respect to Lebesgue measure and its density function $h$ satisfying $||h||_\infty<C$, then the Cauchy transform of $\mu$ is a uniformly continuous function on $\mathbb{C}$. 
\end{proposition}
\begin{proof}
Let $B_\rho(\lambda)$ be the disk $B_\rho(\lambda)=\{ w: |\lambda-w|<\rho \}$. Then
\begin{equation}
  \left|\int_{B_\rho(\lambda)}\frac{1}{\lambda-w}d\mu(w)\right|\leq
  C\int_{B_\rho(\lambda)}\left|\frac{1}{\lambda-w}\right|d^2w\leq C\int_0^{2\pi}\int_0^\rho \frac{1}{r}rdrd\theta=2\pi {C}\rho.
\end{equation}
Given $\varepsilon>0$, if $|\lambda_1-\lambda_2|<\min\{\varepsilon, \varepsilon^3\}$, then for any $w\notin B_\varepsilon(\lambda_1)\cup B_\varepsilon(\lambda_2)$, we note that
\[
 \left| \frac{1}{\lambda_1-w}-\frac{1}{\lambda_2-w} \right|
  =\left|\frac{\lambda_1-\lambda_2}{(\lambda_1-w)(\lambda_2-w)}\right|<\varepsilon.
\]
The condition $|\lambda_1-\lambda_2|<\min\{\varepsilon, \varepsilon^3\}$ ensure that $B_\varepsilon(\lambda_1)\cup B_\varepsilon(\lambda_2) \subset B_{2\varepsilon}(\lambda_i)$ for $i=1,2$. 
We hence have
\begin{align*}
   & \int_{B_\varepsilon(\lambda_1)\cup B_\varepsilon(\lambda_2)} \left|\frac{1}{\lambda_1-w}-\frac{1}{\lambda_2-w}\right| d\mu(w)\\
   &\leq \int_{B_{2\varepsilon}(\lambda_1)}\left|\frac{1}{\lambda_1-w}\right|d\mu(w)+
   \int_{B_{2\varepsilon}(\lambda_2)}\left|\frac{1}{\lambda_2-w}\right|d\mu(w)
   \leq 8\pi{C}\varepsilon, 
\end{align*}
and 
\begin{align*}
  \int_{\mathbb{C}\backslash(B_\varepsilon(\lambda_1)\cup B_\varepsilon(\lambda_2))} \left|\frac{1}{\lambda_1-w}-\frac{1}{\lambda_2-w}\right| d\mu(w)< \varepsilon.
\end{align*}
Therefore, $|C_\mu(\lambda_1)-C_\mu(\lambda_2)|<(8\pi{C}+1)\varepsilon$ provided that $|\lambda_1-\lambda_2|<\min\{\varepsilon, \varepsilon^3\}$. 
\end{proof}

\begin{theorem}
\label{thm:appendix-derivative-identity}
If $\mu$ is a compactly supported probability measure on the plane, suppose that $\mu$ is absolutely continuous with respect to Lebesgue measure and its density function $h$ satisfying $||h||_\infty<C$, and $S_\mu$ has a continuous derivative on some open domain $\mathcal{O}$, then the $z$-derivative of the log potential is the Cauchy transform of $\mu$. 

In particular, for any $z\in\mathbb{C}$ we have
\begin{equation}
  \label{eqn:identity6}
    \frac{\partial S_{\mu_{x+c_t,\varepsilon}}(z)}{\partial z}=\int_\mathbb{C}\frac{1}{z-w}d\mu_{x+c_t,\varepsilon}(w).
\end{equation}
\end{theorem}
\begin{proof}
Since the distributional derivative of a continuously differentiable function concide with the usual derivative almost every where, by Corollary \ref{cor3_distderivative}, $\frac{\partial S_\mu(z)}{\partial z}=C_\mu(z)$ almost everywhere. Proposition \ref{lemma:appendix_continuityCauchy} implies that $C_\mu$ is continuous on $\mathbb{C}$. If $\frac{\partial S_\mu(z)}{\partial z}$ is continuous on $\mathcal{O}$, then it follows that the identity $\frac{\partial S_\mu(z)}{\partial z}=C_\mu(z)$ holds for every $z\in \mathcal{O}$. 

It is shown that $S_{\mu_{x+c_t,\varepsilon}}$ has a continuous derivative on $\mathbb{C}$ and ${\mu_{x+c_t,\varepsilon}}$ is absolutely continuous with respect to Lebesgue measure and its density function is bounded by $1/4\pi {t}$ (see Theorem \ref{thm:main-theorem-from-BYZ2021} or \cite[Lemma 7.11]{BelinschiYinZhong2021Brown}). Hence, \eqref{eqn:identity6} holds for all $z\in\mathbb{C}$. 
\end{proof}

\bibliographystyle{abbrvnat}
\bibliography{BrownMeasure}

\end{document}